\documentclass{article}
\usepackage{graphicx}
\usepackage[english]{babel} 
\usepackage{xcolor}
\usepackage{xcolor,colortbl}
\usepackage[document]{ragged2e}
\usepackage{amsmath}
\usepackage{booktabs} 
\usepackage[bottom=1.8cm,top=1.3cm,left=2.6cm,right=2.6cm,includehead,includefoot]{geometry}
\usepackage{abstract} 
\usepackage{titlesec} 
\usepackage{fancyhdr} 
\usepackage{titling} 
\usepackage{hyperref} 
\usepackage{amssymb}
\usepackage{xspace}
\usepackage[utf8]{inputenc}
\usepackage{newunicodechar}
 \usepackage{relsize}

\newunicodechar{￥}{\textyen}
\DeclareTextCommandDefault{\textyen}{%
  \vphantom{Y}%
  {\ooalign{Y\cr\hidewidth\yenbars\hidewidth\cr}}%
}

\newcommand{\yenbars}{%
  \vbox{
     \hrule height.1ex width.4em
     \kern.15ex
     \hrule height.1ex width.4em
     \kern.3ex
  }%
}
\newcommand\nd{\textsuperscript{nd}\xspace}
\newcommand\rd{\textsuperscript{rd}\xspace}
\newcommand\st{\textsuperscript{st}\xspace}
\newcommand\nth{\textsuperscript{th}\xspace} 
\usepackage{listings}
\usepackage{listings}
\definecolor{Gray}{gray}{0.95}
\pretitle{\begin{center}\justifying\large\bfseries} 
\posttitle{\end{center}} 
\title{\hspace*{3cm} Riemann’s hypothesis is rejected by definition.\\\\\\
Riemann's hypothesis is rejected by definition, because $\boldsymbol{\zeta(s)}$ is not be equal by definition to the particular sum, which it assumes to be equal. $\boldsymbol{R(s)=1 / 2}$ holds only for the zeros of $\boldsymbol{\zeta(s)=0}$ and for the zeros of certain related functions. However, it does not hold for certain special generalized functions of $\boldsymbol{\zeta()}$, such \newline \hspace*{3.3cm}  the Zeta Hurwitz functions and their sums.} 
\author{%
N.Mantzakouras\\ 
\normalsize \href{mailto:nikmatza@gmail.com}{e-mail: nikmatza@gmail.com}\\ 
\normalsize Athens - Greece  \\ 
}
\date{}
\renewcommand{\maketitlehookd}{%
\begin{abstract}
\noindent 
\justifying
\color{gray}
The Riemann zeta function is one of the most Euler's important and fascinating functions in mathematics. Summarizing and analyzing other similar functions to $\zeta(s)$, we obtain zeros with $\operatorname{Re}(s)$ in the interval $(0,1)$ for the generalized case ${ }_{\text{ }\text{ }G} \zeta( m \cdot s , p )={ }_{ m }^{\alpha, \beta} \sigma( s )=\sum_{ k =1}^{\infty} \frac{1}{(\alpha \cdot k +\beta)^{m \cdot s}}$, where $p =\frac{\alpha+\beta}{\alpha},\{\alpha, \beta \in R , \alpha \neq 0\}, s \in C , m \in Z ^{+}$ with 2 general cases (if $\alpha<0$, $m=1$) out of 4 that is a generalized sum which is partially related to $\zeta (s)$ and finally for the Davenport Heilbronn case, which has zeros of a periodic Dirichlet series with different critical lines, within the interval $(0,1)$. But there is the \textbf{Paradox in the Riemann Hypothesis} for which we have non-trivial zeros that zero a sum that is different from what the function $\zeta()$ represents and not what it should be. A strong pre-hypothesis that raises the question that \textbf{the Riemann Hypothesis}, should rejected by definition or its definition should be reconsidered on another basis. But this paradox is negated when we change the sign of the series and make \textbf{the series negative}. That is, the relationship $\zeta(s)=\sum_{n=1}^{\infty} \frac{1}{n^s}$ is not valid, but the correct relationship is $\zeta(s)=\sum_{n=1}^{\infty} \frac{1}{(-n)^s}$. This simple mistake had not been discovered until now. The proof is in chapter 6.7. This does not follow directly, but because the negative sum represented by the function $\zeta( s )$ is the sum of two series that are each zeroed separately. Also if we generalize the hypothesis to $\zeta(q * s)=0, q>1 / 2$ in $R^{+}$ then we have infinite critical lines that originate from and are related to the baseline critical line of $\operatorname{Re}(s)=1 / 2$, but are not identical to it and lie within the interval $(0,1)$. Finally, it is possible to have other constructed cases with relevance of the function $\zeta(z)$ in addition to those considered. Therefore, whether or not the Riemann Hypothesis holds depends on whether or not we accept the fact that the Riemann Hypothesis is or is not fully stated, in the functional equation $\zeta(s)=0$ or directly in the sum it represents, or in more general forms of functions $\zeta()$ i.e. we are talking more about forms of Zeta Hurwitz functions.
\end{abstract}
}
\begin{document}
\maketitle
\justifying
\begin{flushleft}
\section*{Part I. The pseudo-truth of the Riemann hypothesis in relation to the function $\boldsymbol{\zeta()}$ and some related functions}
\textbf{The function $\boldsymbol{\zeta(z)}$ is 1-1 on the critical strip.}\\
\vspace*{5pt}
\textbf{\textit{1.Theorem}}
\end{flushleft}
\normalsize
\textbf{The Hypothesis of Riemann focuses on the point where we must prove that if $\boldsymbol{s = Re(z)}+$ $\boldsymbol{ + Im(z)^{*}I}$,\\\\
I) The functions $\boldsymbol{\zeta(z)}$ and $\boldsymbol{\zeta(1-z)}$ are 1-1 on the critical strip. \\
II) The common roots of the equations $\boldsymbol{\zeta(z)-}$ $\boldsymbol{\zeta(1-z) = 0}$  they have $\boldsymbol{Re(z)=Re(1-z)=1/2}$ within the interval (0,1) and determine unique position, which is called critical line.}
\newpage\noindent
\textbf{I.1. Proof:}\\
The functions $\zeta (z)$ and $\zeta (1-z)$ are 1-1 on the critical strip. For this we need to analyze when and where the exponential function $n^{z}$  is 1-1 when $n\in Z^{+}$. \\
The exponential function $n^{z}$:\\
I) is 1-1 in each of these strips defined by the intervals $2\pi k/\ln(n)$ and $2\pi(k+1)/\ln(n)$ “where $k\in Z$. For this we must prove two cases. If $A\subset C \wedge f: A\rightarrow C$ then if $f(z_{1})=f(z_{2})$ then $z_{1}=z_{2}$ that is.. $x_{1}=x_{2}\wedge y_{1}=y_{2}$. Indeed, if $z_{1}=x_{1}+i\cdot y_{1}$, $z_{2}=x_{2}+i\cdot y_{2}$ are two points within into a such strip  such that $e^{z_{1}}=e^{z_{2}}$ then,
\[n^{x_{1}+y_{1}\cdot i }=n^{x_{1}+y_{1}\cdot i }\Rightarrow\left|n^{z_{1}}\right|=\left|n^{z_{2}}\right|\Rightarrow\left|n^{x_{1}}n^{y_{1}\cdot i}\right|=\left|n^{x_{2}}n^{y_{2}\cdot i}\right|\Rightarrow x_{1}=x_{2}\]
and since $x_{1}=x_{2}$, the relation $e^{z_{1}\cdot \ln (n)}=e^{z_{2}\cdot \ln (n)}$ gives $e^{i\cdot y_{1}\cdot \ln (n)}=e^{i\cdot y_{2}\cdot \ln (n)}$, so $y_{1}-y_{2}=2k\pi/\ln(n)$, is an integer multiple of $2\pi$. But $z_{1}$ and $z_{2}$ belong to the strip, so $\left|y_{1}-y_{2}\right|<2\pi/\ln(n)$. That is, the difference $y_{1}-y_{2}$ is at the same time a multiple of $2\pi/\ln(n)$ and at an absolute value of less than $2\pi/\ln(n)$. The only case that this is true is when $y_{1}=y_{2}$. We finally conclude that $z_{1}=z_{2}$, so $n^{2}$ is 1-1 in the strip $\{z: 2k\pi/\ln(n)\leq Im(z)<2(k+1)\pi/\ln(n)\}.$ ''lines down closed – open up”. We also notice that $n^{z}$ is on $C-\{0\}$. Because if $w\neq 0$ and we put $z=\ln |w|/\ln(n)+iargw/\ln(n)$ then,
\[n^{z}=n^{\ln |w|+iargw/\ln(n)}=n^{\ln|w|}n^{iargw/\ln(n)}=|w|\cdot n^{iargw/\ln(n)}=w\]
\begin{figure}[h!]
\centering
\includegraphics[scale=0.6]{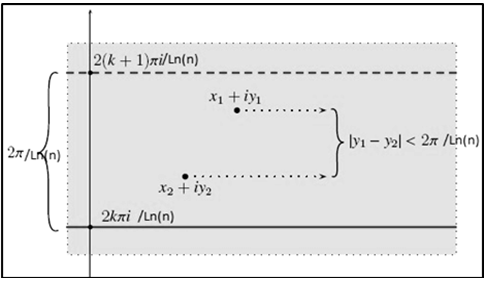}
\caption{The $z_{1}$ and $z_{2}$ belong to the lane width $2\pi/\ln(n)$}
\end{figure}
\flushleft
\justifying
\textbf{Formation of the strips 1-1 for $\boldsymbol{\zeta(s)}$ and $\boldsymbol{\zeta(1-s)}$.}\\\\
\textbf{a1.)} If we accept the non-trivial zeroes on critical strip of the Riemann Zeta Function $\zeta(s)$ as $s_{1}=\sigma_{1}+it_{1}$ and $s_{2}=\sigma_{2}+it_{2}$ with $|s_{2}|>|s_{1}|$, and if we suppose that the real coordinates  $\sigma_{1}$, $\sigma_{2}$ of each non-trivial zero of the Riemann Zeta, $[1, 2, 3, 4]$ function $\zeta(s)$ correspond to two imaginary coordinates $t_{1}$ and $t_{2}$, then, we have the following equations group:
\[\zeta (\sigma_{1}+i\cdot t_{1})=\dfrac{1}{1^{\sigma_{1}+i\cdot t_{1}}}+\dfrac{1}{2^{\sigma_{1}+i\cdot t_{1}}}+\dfrac{1}{3^{\sigma_{1}+i\cdot t_{1}}}+...+\dfrac{1}{n^{\sigma_{1}+i\cdot t_{1}}}+...=0\]
\[\zeta (\sigma_{2}+i\cdot t_{2})=\dfrac{1}{1^{\sigma_{2}+i\cdot t_{2}}}+\dfrac{1}{2^{\sigma_{2}+i\cdot t_{2}}}+\dfrac{1}{3^{\sigma_{2}+i\cdot t_{2}}}+...+\dfrac{1}{n^{\sigma_{2}+i\cdot t_{2}}}+...=0\]
\\
Taking the first equation and deducting the second, we obtain:
\[\zeta (\sigma_{1}+i\cdot t_{1})-\zeta (\sigma_{2}+i\cdot t_{2})=\sum\limits_{n=1}^{\infty}\left(\dfrac{1}{n^{\sigma_{1}+i\cdot t_{1}}}-\dfrac{1}{n^{\sigma_{2}+i\cdot t_{2}}}\right)=\]
\[\sum\limits_{n=1}^{\infty}\dfrac{n^{\sigma_{2}+i\cdot t_{2}}-n^{\sigma_{1}+i\cdot t_{1}}}{n^{\sigma_{2}+i\cdot t_{2}}\cdot n^{\sigma_{1}+i\cdot t_{1}}}=\sum\limits_{n=1}^{\infty}\dfrac{n^{\sigma_{2}}e^{i\cdot t_{2}\ln(n)}-n^{\sigma_{1}}e^{i\cdot t_{1}\ln(n)}}{n^{\sigma_{2}+i\cdot t_{2}}\cdot n^{\sigma_{1}+i\cdot t_{1}}}=0\]\\
From the previous relation, we conclude that if 
\[n^{\sigma_{2}+i\cdot t_{2}}\cdot n^{\sigma_{1}+i\cdot t_{1}}\neq 0,\text{ } n>1\text{ then } \sigma_{2}=\sigma_{1}\wedge t_{2}=t_{1}\pm \dfrac{2k\pi}{\ln(n)}, (k=1,2,...). \]
That is to say $t_{1}$ and $t_{2}$ can take any value, but according to the previous relation. So from (Fig.1), that means $\zeta (s)$ is 1-1 on the lane of critical Strip in Intervals such as defined, therefore and on the critical Line. So $n^{-z}$ is 1-1 in the strip and $\{z:2k\pi/\ln(n)\leq Im(z)<2(k+1)\pi/\ln(n)\}$.''lines down closed – open up”.\\\\
\textbf{a2.)} If we do the same work with $|s_{2}|>|s_{1}|$ for the case $\zeta(1-s)$ will we have:
\[\zeta (1-\sigma_{1}-i\cdot t_{1})-\zeta (1-\sigma_{2}-i\cdot t_{2})=\sum\limits_{n=1}^{\infty}\left(\dfrac{1}{n^{1-\sigma_{1}-i\cdot t_{1}}}-\dfrac{1}{n^{1-\sigma_{2}-i\cdot t_{2}}}\right)=\]
\[\sum\limits_{n=1}^{\infty}\dfrac{n^{1-\sigma_{2}-i\cdot t_{2}}-n^{1-\sigma_{1}-i\cdot t_{1}}}{n^{1-\sigma_{2}-i\cdot t_{2}}\cdot n^{1-\sigma_{1}-i\cdot t_{1}}}=\sum\limits_{n=1}^{\infty}\dfrac{n^{1-\sigma_{2}}e^{-i\cdot t_{2}\ln(n)}-n^{1-\sigma_{1}}e^{-i\cdot t_{1}\ln(n)}}{n^{1-\sigma_{2}-i\cdot t_{2}}\cdot n^{1-\sigma_{1}-i\cdot t_{1}}}=0\]
From the previous relation, we see that if
\[n^{1-\sigma_{2}-i\cdot t_{2}}\cdot n^{1-\sigma_{1}-i\cdot t_{1}}\neq 0, n>1\]
\[\& (\{1-\sigma_{2}=1-\sigma_{1}\}\vee \sigma_{2}=\sigma_{1})\wedge t_{2}=t_{1}\pm\dfrac{2k\pi}{\ln(n)}, k=(1,2,...)\]
So, from (Fig.1), it follows that $\zeta (1-s)$ is 1-1 on the  critical Strip, So $n^{-(1-s)}$ is 1-1 in the strip and $\{s:2k\pi/\ln(n)\leq Im(s)\leq 2(k+1)\pi/\ln(n)\}$.''lines down closed – open up”\\\\
As we have seen before, the complex exponential form $n^{s}$ but and $n^{-s}$ is also periodic with a period of $2\pi/\ln(n)$. For example $n^{s+2\cdot \pi\cdot i}=$ $n^{s}\cdot n^{2\cdot\pi\cdot i/\log(n)}=n^{s}$ is repeated in all the horizontal strips with $2\pi/\ln(n)$ on (Fig.2) below.\\
\begin{figure}[h!]
\centering
\includegraphics[scale=0.65]{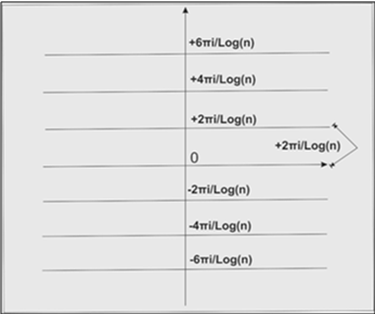}
\caption{All the ''bottom - closed upper “horizontal strips with $2\pi/\ln(n)$ for $n^{-s}$, $n\geq 1$, $n\in Z$.}
\end{figure}\\\\
If we are asked to find the strips that are 1-1 of the function $\zeta(s)=\sum\limits_{i=1}^{n}\dfrac{1}{n^{s}}$ it will be the union of the strips formed with period $2\pi k/\ln(n)$ where $k,n\in Z$. This results from the analysis we made for the cases a1, a2. Therefore, if we analytically assume that $n=2,3,...$ then the exponential complex function will be 1-1 in each of the strips, i.e. in each one set of the form le
\[S_{n}=\{s:2k\pi/\ln(n)\leq Im(s)<2(k+1)\pi/\ln(n)\}, k\in Z\]
and therefore in each one subset of the unified expanded set $S=\bigcup\limits_{i=2}^{\infty}S_{i}$. The largest strip obviously contains all the rest it is wide, and has width $2\pi/\ln 2$ i.e. with $n = 2$.\\\\\\
\textbf{I.2. Overlapping contiguous Strips.}\\
If we start from the upper bound of the strip with $n = 2$ and $k=1$, Fig.3 the upper ones will descend to zero by 
increasing n with a width difference the interval $\delta = 2\pi/\ln(n-1) - 2\pi/\ln(n)$. In addition, we will prove that the strip $i=n\geq 2$ primary (k=1) that generated in ascending order, the upper limit of that created of its double width enters the zone of the previous of $(i =n-1, k = 1)$ also primary. With this logic we can find for which n we have overlapping strips in relation to the increase of n on the imaginary axis. Below we find for which n this dimension applies which is useful to distinguish them but also to make them visible schematically.
\begin{figure}[h!]
\centering
\includegraphics[scale=0.7]{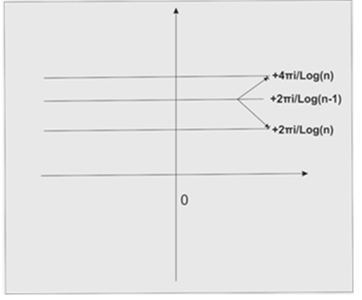}
\caption{The upper ones will go down to zero by increasing n}
\end{figure}
\flushleft 
\justifying
Apply $2\pi/\ln(n-1)<2\cdot 2\cdot\pi/\ln n \Rightarrow n<(n-1)^{2}$ $\Rightarrow n^{2}-3\cdot n +1 >0\Rightarrow n>\dfrac{3+\sqrt{5}}{2}$ and because we talked about integers then  $n\geq 3$. The same phenomenon is also created at higher level strips.\\\\
\textbf{I.3.} The $\boxed{\zeta (z)=\sum\limits_{i=1}^{n}\dfrac{1}{n^{z}}}$ it is directly 1-1 from the comparison of 2 points $z_{1}, z_{2}$.\\\\
\textbf{Proof:}\\
From the equality\\
\[\zeta (z_{1})=\sum\limits_{i=1}^{n}\dfrac{1}{n^{z_{1}}}=\sum\limits_{i=1}^{n}\dfrac{1}{n^{z_{2}}}=\zeta (z_{2})\Rightarrow 1+2^{-(x_{1}+y_{1}\cdot i)}+3^{-(x_{1}+y_{1}\cdot i)}+...+n^{-(x_{1}+y_{1}\cdot i)}=\zeta (z_{2})=\]

\[=1+2^{-(x_{2}+y_{2}\cdot i)}+3^{-(x_{2}+y_{2}\cdot i)}+...+n^{-(x_{2}+y_{2}\cdot i)}\Rightarrow\]
with comparison of similar terms the result arises.
\[\boxed{x_{1}=x_{2}\wedge y_{1}=y_{2}-\dfrac{2\cdot \pi\cdot k}{\ln(n)}, k\in Z}\]
After to the analysis of cases before, will be arises 1-1 as shown in sections \{I$_{\text{1}}$.a$_{\text{1}}$, I$_{\text{1}}$.a$_{\text{2}}$\} pages 1\& 2. Also according to Lagrange's generalized theory, each of the roots according to the obvious relation $\zeta (z)=\zeta (1-z)$ will result from the generalized theorem of Lagrange. From the equation that results from 3 Functional Equations, with condition of common roots. We take $\zeta (z)=\zeta (1-z)$ and for first approach we then follow the 3 first terms on each $\zeta$-equation we will have in the Regular form $\zeta (z)$, $\zeta (1-z)$ with the analysis below.\\
\vspace*{0.5\baselineskip}
\[\begin{array}{l}
1^{-z}+2^{-z}+3^{-z}=1^{-(1-z)}+2^{-(1-z)}+3^{-(1-z)}\Rightarrow 2^{-z}+3^{-z}=\dfrac{1}{2}2^{z}+\dfrac{1}{3}3^{z} \Rightarrow \\
2^{-z}+3^{-z}-2^{-(1-z)}-3^{-(1-z)}=0\Rightarrow 6^{z}+\dfrac{3}{2}4^{z}-3-3\left(\dfrac{2}{3}\right)^{z}=0
\end{array}\]
With the replacement
\[6^{z}=y\Rightarrow z\cdot Log(6)=Log(y)+2\cdot k\cdot \pi\cdot i\Rightarrow z=\dfrac{Log(y)+2\cdot k\cdot \pi\cdot i}{Log(6)}\]
\[\boxed{z_{in}=\dfrac{Log(y)+2\cdot k\cdot \pi\cdot i}{Log(6)}}\]
we consider it as initial value for the solution of Transcendental equation with the method of Lagrange inversion, then: 
\\\\
$\resizebox{\hsize}{!}{$\boxed{z=\dfrac{Log(y)+2\cdot k\pi\cdot i }{Log(6)}+\sum\limits_{w=1}^{\infty}\dfrac{(-1)^{w}}{\Gamma (w+1)}\dfrac{d^{w-1}}{dy^{w-1}}\left(\left(\dfrac{1}{y\cdot Log(6)}\right)\cdot\left(\dfrac{3}{2}\cdot 4^{\dfrac{Log(y)+2\cdot k\cdot\pi\cdot i}{Log(6)}}-3\cdot\left(\dfrac{2}{3}\right)^{\dfrac{Log(y)+2\cdot k\cdot\pi\cdot i}{Log(6)}}\right)^{w}\right)_{y\rightarrow 3}} \text{ (IL)}$}$\\\\
Some solutions for $k\in Z$ (Table 1)
\begin{table}[h!]
\centering
\begin{tabular}{|cccc|}
 \hline
 \rowcolor{Gray}
 1\st   0.4985184869351  & -/+ &3.6072014300530 I,&k=1\\
  \rowcolor{Gray}
 2\nd  0.4997077392995 & -/+ &10.329023684711  I,&k=3\\
   \rowcolor{Gray}
 3\rd  0.5003426560276& -/+ & 13.990255078243  I,&k=4\\
   \rowcolor{Gray}
 4\nth 0.4993087825728 & -/+ & 17.585326315022 I,&k=5\\
    \rowcolor{Gray}
 5\nth 0.5026747758005& -/+ &  21.304192633721 I,&k=6\\[1ex]  
 \hline
\end{tabular}
\caption{6 Roots of equation $\zeta (z)=\zeta (1-z)$}
\end{table}\\
And so on for the sequence of infinity. Now we see that we approach the real roots and its complexities. With the additional method Newton or the Bisection method we approach more rapidly to them. And these roots will be contained in corresponding strips that we have previously defined and in closer relation with their position on the critical line.\\
By this logic for the previous approximation equation in relation to its solutions of, $\zeta(z)=0$ the strips which determine the imaginary roots of the equation, will be between the intervals $2\cdot \pi \cdot k/\ln 6$, $k\in Z$ as defined by the analysis of the generalized Theorem of Lagrange. The strips of the imaginary part of the roots, as they appear for 6 consecutive intervals of the approximate equation with method of Lagrange that we have solved above and we will have the results on the (Table 2).
\begin{table}[h!]
\centering
\textbf{Strips – Imaginary Roots:}\\
\begin{tabular}{|ll|}
 \hline
 \rowcolor{Gray}
 1\st  (7.013424977055179  &3.506712488527589)\\ 
  \rowcolor{Gray}
 2\nd  (10.52013746558276  & 7.013424977055179 \\
  \rowcolor{Gray}
 3\rd  (14.02684995411035  &10.52013746558276) \\
  \rowcolor{Gray}
 4\nth (17.53356244263794  & 14.02684995411035) \\
  \rowcolor{Gray}
 5\nth (21.04027493116553  &17.53356244263794) \\
  \rowcolor{Gray}
 6\nth (24.54698741969312 &21.04027493116553)\\
 \hline
\end{tabular}
\caption{The strips of the imaginary part of the roots}
\end{table}\\
It is clear the first root $3,607..i$ is located in 1\st of strip $(7.01-3.506)$, the second $10.3290..i$ in the 2\st strip $(7.013-10.52)$ etc. Now the root $13.99..i$ which approximates the first root of $\zeta(z) = \zeta(1-z) = 0$ i.e. the $14.1347..i$, is in the 4\st strip $(14.02-17.533)$ with a lower limit the $14.026..i$ a value that is very close to the one  required for the approximate and exclusive (1\st root) when $\zeta(z) = 0$ of the Riemann Hypothesis.\\\\
\textbf{Corollary 1.}\\
The only roots of the Zeta function not included in the set $\{ z \in \mathbb{C}: 0 \leq Re(z)\leq 1 \}$ are the points $-2, -4, -6..$ [11. p47]\\\\
\textbf{Proof}\\
On the functional equation:\\
$\zeta(1-z)= 2 \cdot \zeta(z) \cdot \Gamma (z) \cdot \cos( \pi / 2\cdot z) \cdot(2\pi)^{-z}$ we know that for $Re(z)>1$ the functions $\zeta(z)$ and $\Gamma(z)$ do not equal zero (Proof-Th.3.p6). We also find that for $Re (z)> 1$ the  $1 -Re (z) <0$ and putting $u = 1-z$ we will find all the roots of $\zeta(u)$ for $Re(u)<0$. Therefore we will have that $\zeta(u)$ will be zero, where it is zeroed the i.e. $z = 3,5,7$. Then the roots of $\zeta(u)$ for $Re(u)<0$ will be the points $-2, -4, -6...$ and will be all the roots of the function $\zeta$ out of the strip $\{ z \in \mathbb{C}: 0 \leq Re(z) \leq 1\}$.\\\\
\textbf{Theorem-Helpful 1}\\
The Riemann's Z-function has no roots on the lines $Re(z)=1$ and $Re(z)=0$.\\
The proof is detailed in the book [11. p50-51].\\\\
\textbf{II. The common roots of the equations $\boldsymbol{\zeta(z)-\zeta(1-z)=0}$  they have $\boldsymbol{Re(z)= Re(1-z)=1/2}$ within the interval $\boldsymbol{(0,1)}$  if moreover apply $\boldsymbol{\zeta(z)=0}$. (Refer to Theorem 3, p.8)\\\\
Proof:} Let us assume $z$ to be such that for complex $z=x_0+iy_0$, $0<Re(z)<1\wedge Im(z)\neq0$ and $\zeta(z)=\zeta(1-z)=0$. According to the two equations, they must apply to both, because they are equal to zero that: $\zeta(z)=\zeta(\overline{z})=0 \wedge \zeta(1-z)=\zeta(\overline{1-z})=0$. But from [Theorem 1, II, a1, a2] the $\zeta(z)$ and $\zeta(1-z)$ are 1-1 on the critical Strip. If suppose generally that $\zeta(x_0 \pm y_0 \cdot i)=\zeta(1-x_{0}' \pm y_0' \cdot i)=0$ then:
\[ y_0 = y_0' \wedge x_0=1-x_0',\]
but because are the 1-1 then we will apply two cases for complex $z$ i.e.
\[x_0'=1-x_0,\text{ } x_0 \neq x_0',\text{ } y_0=y_0' \vee x_0'=1-x_0,\text{ } x_0=x_0'=\dfrac{1}{2}, \text{ } y_0=y_0'\]
\{Corollary 1, Theorem Helpful 1\}, because we suppose $\zeta(z)=\zeta{1-z}=0$.] We conclude that $0<x_0<1$. Therefore we have three cases:\\\\
\textbf{II.1.} If $\boxed{x_0=x_0'=\dfrac{1}{2}}$ we apply the obvious i.e.
$\zeta\left(1 / 2 \pm y_{0} \cdot i\right)=\zeta\left(1-1 / 2 \pm y_{0}' \cdot i\right)=\zeta\left(1 / 2 \pm y_{0}' \cdot i\right)=0 \wedge$
$\wedge {y}_{0}={y}_{0}'.$ Which fully meets the requirements of
hypothesis!\\
\flushleft
\textbf{II.2.} If $\boxed{0<{x}_{0} \neq {x}_{0}' \neq \dfrac{1}{2}<1}$ in this case if $x_{0}<x_{0}' \wedge x_{0}+x_{0}'=1$ or $x_{0}>x_{0}' \wedge x_{0}+x_{0}'=1$ then
for the Functional equation $\zeta(z)-\zeta(1-z)=0$ simultaneously apply:
\[\zeta\left(x_{0}+y_{0} \cdot i\right)=\zeta\left(x_{0}-y_{0} \cdot i\right)=0 \wedge \zeta\left(x_{0}'+y_{0}' \cdot i\right)=\zeta\left({x}_{0}^{\prime}-{y}_{0}^{\prime} {i}\right)=0.\] So let's assume that:\\
\vspace*{0.5\baselineskip}
1\st. $x_{0}<x_{0}' \wedge x_{0}+x_{0}'=1 \Rightarrow x_0=\dfrac{1}{2}-a \wedge x_{0}'=\dfrac{1}{2} + b , a \neq b$. But apply ${x}_{0}+{x}_{0}^{\prime}=1 \Rightarrow 1+b-a=1 \Rightarrow$ $b=a$. That is, $a, b$ is symmetrical about of $\dfrac{1}{2}$. But the $\zeta(z)$ and $\zeta(1-z)$ are 1 - 1 on the critical Strip then apply the relation:
\[ \zeta\left(\dfrac{1}{2}-a \pm {y}_{0} \cdot {i}\right) = \zeta\left(\dfrac{1}{2}+a \pm y_0' \cdot i\right) \wedge y_0 = y_0' \Rightarrow \text{must } \dfrac{1}{2}-a=\dfrac{1}{2}+a \Rightarrow 2a=0 \Rightarrow a=0\]
therefore $x_{0}=x_{0}'$\\
\vspace*{0.5\baselineskip}
2\nd. $x_{0}>x_{0}' \wedge x_{0}+x_{0}'=1 \Rightarrow x_{0}=\dfrac{1}{2}+a \wedge x_{0}^{\prime}=\dfrac{1}{2}-a$, as before. Then because the $\zeta(z)$ and $\zeta(1-z)$ are 1-1 on the critical Strip will apply:
\[\zeta\left(\dfrac{1}{2}+a \pm {y}_{0} \cdot {i}\right)=\zeta\left(\dfrac{1}{2}-a \pm {y}_{0}^{\prime} \cdot {i}\right) \wedge {y}_{0}={y}_{0}^{\prime}\] that
$\dfrac{1}{2}+a=\dfrac{1}{2}-a \Rightarrow 2 a=0 \Rightarrow a=0$
therefore ${x}_{0}={x}_{0}^{\prime}$. We see that in all acceptable cases
is true that "The non-trivial zeros of $\zeta({z})-\zeta(1-{z})=0$ have real part which is equal to $1/2$ within the interval $(0,1)$. In the end, we see still one partial case that appears in the roots of the equation $\zeta(z)$ $-\zeta(1-z)=0$.\\
\vspace*{\baselineskip}
\flushleft
\justifying
\textbf{II.3.} If $\boxed{1<x_{0} \wedge x_{0}'<0}$ furthermore, in this case with that is, it is symmetrical about ${x}_{0}+{x}_{0}'=1$ then for the functional equations we apply:
\[\zeta\left(x_{0}-y_{0} \cdot i\right)=\zeta\left(1-x_{0}-y_{0} \cdot i\right) \wedge \zeta\left(x_{0}+y_{0} \cdot i\right)=\zeta\left(1-x_{0}+y_{0} i\right) \wedge y_{0}=y_{0}'\] 
and for the three cases it is valid $y_{0}=0 \vee y_{0} \neq 0$. This case applies only when $\zeta(z)=\zeta(1-z) \neq 0 \wedge x_{0}=$ $=1+a$, $x_{0}'=-a$, and not for the equality of $\zeta$ function with zero. But in order to have common roots of the two functions $\zeta(z)$ and $\zeta(1-z)$, «because they are equal to zero and equal to each other», then first of all, it must be i.e $x_{0}=x_{0}'=1-x_{0}$ and also because the $\zeta(\mathrm{z})$ is 1-1 and for each case applies $(y_{0}=0$ and $y_{0}'=0)$ or $\left(y_{0} \neq 0\right.$ and $\left.y_{0}'\neq 0\right)$, therefore we will apply $x_{0}=x_{0}'=\dfrac{1}{2} \wedge y_{0}=y_{0}^{\prime}$ and with regard to the three cases \{II.1, II.2, II.3\} the cases \{II.2, II.3\} cannot happen, because they will have to be within $(0,1)$ and will
therefore be rejected.\\\\
\emph{
Therefore if $\zeta(z)=\zeta(1-z)=0$ then because "these two equations of zeta function are 1-1 on the lane of critical strip" as shown in theorem 1, is "sufficient condition" that all non - trivial zeros are on the critical line $Re(z)=1/2$. That means that the real part of $z$ of $\zeta(z)=0$, equals to $1/2$". So Theorem 1 \& Corollary 1}, they have been proved and also helps Theorem 3 below.\\\\\\
\textbf{2. Theorem 2}\\
For the non- trivial zeroes of the Riemann Zeta Function $\zeta(s)$ apply\\\\
There exists an upper-lower bound of  $Re(s)$ of the Riemann Zeta Function $\zeta(s)$ and more specifically in the closed interval $\left[\dfrac{\ln 2}{\ln 2 \pi}, \dfrac{\ln \pi}{\ln 2 \pi}\right]$. The non-trivial zeroes of the Riemann Zeta function $\zeta(s)$ of the upper-lower bound are distributed symmetrically on the straight line $Re(s)=1/2$. The average value of the upper lower bound of $Re(s)=1/2$.\\\\
\textbf{Proof:}\\
Here, we formulate two of the functional equations from E - q. Set
\[
\def\arraystretch{1.5}\begin{array}{ll}
\zeta(1-s) / \zeta(s)=2(2 \pi)^{-s} \cos(\pi s / 2) \Gamma(s), &{Re}(s)>0\\
\zeta(s) / \zeta(1-s)=2(2 \pi)^{s-1} \sin(\pi s / 2) \Gamma(1-s), &Re(s)<1
\end{array}\]
We look at each one equation individually in order to identify the set of values that we want each time.\\\\
\textbf{2.1.} For the first equation and for real values with $Re(s)>0$ and by taking the logarithm of two sides of the equation [5], we have:
\[\zeta(1-s) / \zeta(s)=2(2 \pi)^{-s} \cos(\pi / 2) \Gamma(s) \Rightarrow\]
\[\log [\zeta(1-s) / \zeta(s)]=\log [2]-s \log [2 \pi]+ \log [\cos(\pi s / 2) \Gamma(s)]+2k \pi i\]\\
but solving for $s$ and if $f(s)=\cos(\pi s / 2) \Gamma(s)$ and from Lemma 2 $[6, p. 6]$:\\
if $\lim\limits_{s \to s_{0}} \zeta(1-s) / \zeta(s)=1$ or $\log \left[\lim\limits_{s \rightarrow s_{0}} \zeta(1-s) / \zeta(s)\right]=0$ we get $s=\dfrac{\log [2]}{\log [2 \pi]}+\dfrac{\log [f(s)]+2 k \pi i}{\log [2 \pi]}$ with $Re\left(\dfrac{\log [f(s)]+2 k \pi i}{\log [2 \pi]}\right) \geq 0$, finally because we need real $s$ we will have $\boxed{Re(s) \geq \dfrac{\log [2]}{\log [2 \pi]}=0.3771}$. This is the lower bound, which gives us the first $\zeta(s)$ of Riemann's Zeta Function.\\\\\\\\\\
\textbf{2.2.} For the second equation, for real values if $Re(s)<1$ and by taking the logarithm of the two parts of the equation, we will have:\\
\[\zeta(s) / \zeta(1-s)=2(2 \pi)^{s-1} \sin (\pi s / 2) \Gamma(1-s) \Rightarrow\]
\[ \log [\zeta(s) / \zeta(1-s)]=\log [2]+(s-1) \log [2 \pi] + \log [\sin(\pi \cdot s / 2) \Gamma(1-s)]+2 k \pi i\]\\
but solving for $s$ and if $f(s)=\sin(\pi s / 2) \Gamma(1-s)$ and from Lemma 2[6, p.6] if $\lim\limits_{s \to s_{0}} \left[\zeta(s) / \zeta(1-s)\right]=1$ or $\log \left[\lim\limits_{s \to s_{0}} \zeta(s) / \zeta(1-s)\right]=0$ we get
$s=\dfrac{\log [\pi]}{\log [2 \pi]}-\dfrac{\log [f(s)]+2 k \pi i}{\log [2 \pi]}$ with $Re\left(\dfrac{\log [f(s)]+2 k \pi i}{\log [2 \pi]}\right)>=0$.\\
In the following, because we need real $s$ we will take. $\boxed{Re(s) \leq \dfrac{\log [\pi]}{\log [2 \pi]}=0.6228}$. This is the upper bound, which gives us the second of Riemann's Zeta Function of $\zeta(s)$. So, we see that the lower and the upper bound exist for $Re(s)$ and they are well defined.\\\\
\textbf{2.3.} Assuming that $s_{k}=\sigma_{\text {low }}+i t_{k}$ and $s_{s}^{\prime}=\sigma_{\text {upper }}+i t_{k}$ with $\sigma_{\text {low }}=\dfrac{\log [2]}{\log [2 \pi]}$ and $\sigma_{\text {upper }}=\dfrac{\log [\pi]}{\log [2 \pi]}$ we apply $\sigma_{1}=\sigma_{\text {low }}$ and $\sigma_{2}=\sigma_{\text {upper }}$, (Fig.4). If we evaluate the difference and 
\[Re \text{ } \Delta_{\kappa}^{\prime}=Re\text{ }s_{\kappa}^{\prime}-1 / 2=\sigma_{upper}-1 / 2=0.1228\] 
\[Re \text{ }\Delta s_{\kappa}=1 / 2-Re\text{ } s_{\kappa}=1 / 2-\sigma_{low}=0.1228.\]
This suggests for our absolute symmetry of 
\begin{figure}[h!]
\centering
\includegraphics[scale=0.27]{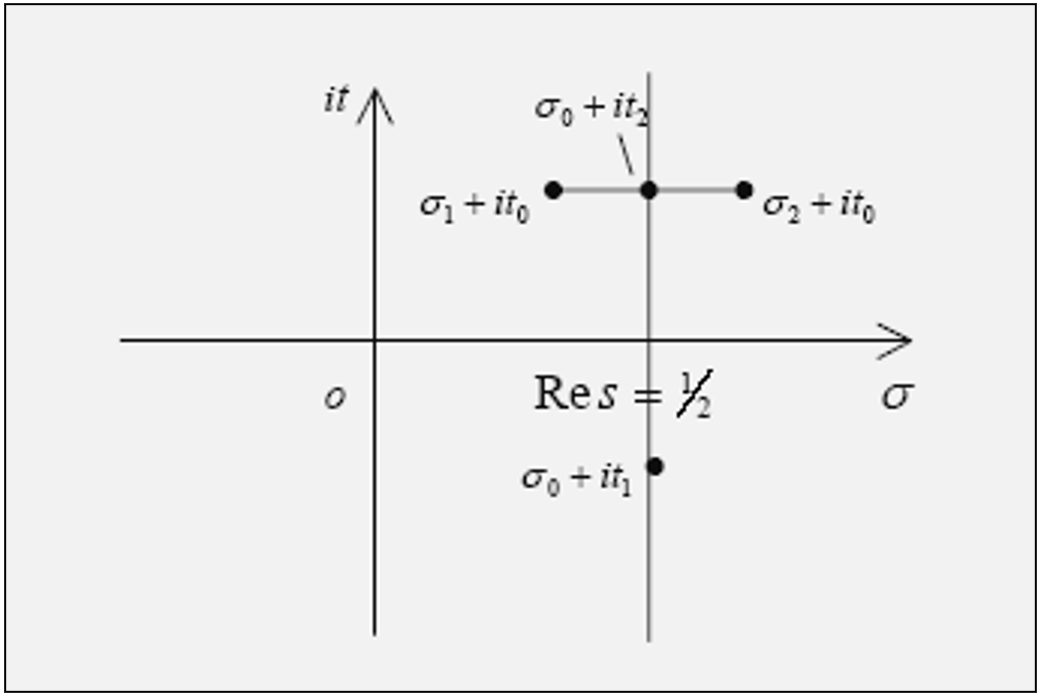}
\caption{Arrangement of low and upper of  real part of zeros from 2 functional equations of Riemann.}
\end{figure}\\\\
\textbf{2.4.}  The average value of the upper lower bound is $Re(s)=1/2$ because from (Figure 4): that is apply
\[Re(s)=1/2 \cdot \dfrac{\log[2]+\log[\pi]}{\log[2 \pi]}=1/2\]\\
\textbf{3. Theorem 3}\\
\textit{The Riemann Hypothesis states that all the non-trivial zeros of $\zeta(z)$ have real part equal to $1/2$.}\\\\
\justifying
\textbf{Proof:} In any case, we Assume that: The Constant Hypothesis $\zeta(z)=0$, $z \in \mathbb{C}$. In this case, we use the two equations of the Riemann zeta function, so if they apply what they represent the $\zeta(z)$ and $\zeta(1-\mathrm{z})$ to equality. Before developing the method, we make the three following assumptions:\\\\\\\\
\textbf{3.1. Analysis of specific parts of transcendental equations, which are detailed}\\\\
\textbf{3.1.1.} For $z \in \mathbb{C},\left\{a^{z}=0 \Rightarrow a=0 \wedge Re z>0\right\}$ which refers to the inherent function similar to two Riemann zeta functions as $(2 \pi)^{-z}=0$ or $(2 \pi)^{\tilde{-1}}=0$ and it seems that they do not have roots in $\mathrm{C}-\mathrm{Z}$, because $(2 \pi) \neq 0$.\\\\
\textbf{3.1.2.} \textbf{This forms $\boldsymbol{\Gamma(z)=0}$ or $\boldsymbol{\Gamma(1-z)=0}$ do not have roots in $\boldsymbol {\mathbb{C}-\mathbb{Z}}$.}\\\\
\textbf{3.1.3.} \textbf{Solution of $\boldsymbol{\sin(\pi/2z)=0}$ or $\boldsymbol{\cos(\pi/2 z)=0}$.} More specifically if $\boldsymbol{z=x+yi}$, then.\\\\
1. $\sin(\pi / 2 z)=0 \Rightarrow \cosh[(\pi \cdot y) / 2] \sin[(\pi \cdot x) / 2]+ I \cdot \cos[(\pi \cdot x) / 2] \sinh[(\pi \cdot y) / 2]=0$\\
If $k, m \in N$ \textit{(Integers)} then all the solutions can be found with a program by the language of Mathematica.\\
\[\begin{array}{rl}
i) &\{x=(2+4 \cdot k), y=2 \cdot i \cdot(1+2 \cdot m\}\\
ii) &\{x=4 \cdot k, y=2 \cdot i \cdot(1+2 \cdot m\}\\
iii)& \{x=(-1+4 \cdot k), y=i \cdot(-1+4 \cdot m\}\\
iv) &\{x=(-1+4 \cdot k), y=i \cdot(1+4 \cdot m\}\\
v) &\{x=(1+4 \cdot k), y=i \cdot(-1+4 \cdot m\}\\
vi)& \{x=(1+4 \cdot k), y=i \cdot(1+4 \cdot m\}\\
vii)& \{x=(2+4 \cdot k), y=i \cdot(4 \cdot m\}\\
viii) &\{x=(2+4 \cdot k), y=2 \cdot i \cdot(1+2 \cdot m\}\\
\end{array}\]
From the generalized solution, it seems that in the pairs $(x, y)$ will always arise integers which make impossible the case $0 < x <1$, therefore there are no roots of the equation $\sin(\pi z / 2)=0$ in $\mathbb{C}-\mathbb{Z}$.\\
Also, because as we see, all the roots given by the union of the sets $x \geq 1 \vee x=-1 \vee x=0$ all are Integers. 
So it cannot be true that for our case:

\[ \cos(\pi / 2 z)=0 \Rightarrow \cosh[(\pi y) / 2] \cos[(\pi x) / 2]-I \sin[(\pi \cdot x) / 2] \sinh[(\pi y) / 2]=0\]\\
If $k, m \in N$ \textit{(Integers)} and the solutions are:\\
\[\begin{array}{rl}
i) &\{x=4 \cdot k, y=i \cdot(-1+4 \cdot m)\}\\
ii) &\{x=4 \cdot k, y=i \cdot(1+4 \cdot m\}\\
iii) &\{x=(-1+4 \cdot k), y=i \cdot(4 \cdot m\}\\
iv) &\{x=(-1+4 \cdot k), y=2 \cdot i \cdot(1+2 \cdot m\}\\
v) &\{x=(1+4 \cdot k), y=4 \cdot i \cdot m\}\\
vi) &\{x=(1+4 \cdot k), y=i \cdot(1+4 \cdot m\}\\
vii) &\{x=(2+4 \cdot k), y=i \cdot(-1+4 \cdot m\}\\
viii) &\{x=(2+4 \cdot k), y=i \cdot(1+4 \cdot m\}\\
\end{array}\]
As we can see, again from the generalized solution, it seems that in the pairs $(x, y)$ will always arise integers, which make impossible the case is $0<x<1$, therefore there are not roots receivers of the equation $\cos(\pi z/2)=0$. Since all the roots are given by the union of the sets $x \geq 1 \vee x=-1 \vee x=0$ it follows that they are Integers.\\\\
\underline{\textit{Great Result}}\\
«With this three - cases analysis, we have proved that the real part of $x$, of the complex $z=x+y i$ cannot be in the interval $0<x<1$, and in particular in cases $\{3.1.1, 3.1.2, 3.1.3\}$, when they are zeroed. Therefore they cannot represent roots in the critical line».\\\\
\textbf{3.2.} Therefore, now we will analyze the two equations of the Riemann zeta function and we will try to find any common solutions.\\\\
\textbf{3.2.1.} For the first equation and for real values with $Re(z)>0$ we apply:\\
Where $f(z)=2(2 \pi)^{-2} \cos(\pi \cdot z / 2) \Gamma(z) \Rightarrow \zeta(1-z)=f(z) \cdot \zeta(z)$ but this means that the following two cases occur:\\\\
\textbf{3.2.1.1.} If $\boxed{\zeta(1-z)=\zeta(z)}$, where $z$ is complex number.\\
This assumption implies that $\zeta(z) \cdot(1-f(z))=0 \Rightarrow \zeta(z)=0$. In theorem 1.(I , II), (pages $1-8$), we showed that the functions $\zeta(z)$ and $\zeta(1-z)$ are 1-1 and therefore if $\zeta\left(x_{0}+y_{0} i \right)=\zeta\left(x_{0}^{\prime}+y_{0}^{\prime}i\right)$
then, but we also apply that $\zeta\left(x_{0}+y_{0} i\right)=\zeta\left(x_{0}-y_{0} i\right)=0$, because we apply it for complex roots. The form $\zeta(1-z)=\zeta(z)$ means that if $z=x_{0}+y_{0}$ i we apply $1-x_{0}=x_{0}$ namely $x_{0}=1 / 2$ because in this case it will be verified that $\zeta\left(1-x_{0}-y_{0}^{\prime} i\right)=\zeta\left(1 / 2-y_{0}^{\prime} i\right)=\zeta\left(1 / 2+y_{0}^{\prime} i\right)= \zeta\left(1-x_{0}+y_{0}^{\prime} i\right)=0$ and because $y_{0}=y_{0}^{\prime} \Rightarrow \zeta\left(x_{0}-y_{0} i\right)=\zeta\left(1 / 2-y_{0} i\right)=\zeta\left(1 / 2+y_{0} i\right)=\zeta\left(x_{0}+y_{0} i\right)=0$ and these are the two forms of the $\zeta$ equation, they have common roots and therefore it can be verified by the definition of any complex equation, when it is equal to zero. Therefore if $z=x_{0}+y_{0}$ i then $x_{0}=1 / 2$ which verifies the equation $\zeta(z)=0$.\\\\
\textbf{3.2.1.2.} If $\boxed{\zeta(1-z) \neq \zeta(z)}$, when $z$ is a complex number.\\
To verify this case must be:
\[\zeta(z)=\zeta(1-z) \cdot f(z) \Rightarrow \zeta(z)=0 \wedge f(z)=0.\]
But this case is not possible, because as we have shown in Section 3.1 (3.1.1, 3.1.2, 3.1.3), the individual functions of $f(z)$ cannot be zero when $z$ is a complex number.\\\\
\textbf{3.2.2. For the second equation and for any real values with $\boldsymbol{Re(z)<1}$ we apply}:\\
$\zeta(z) / \zeta(1-z)=2(2 \pi)^{z-1} \sin (\pi \cdot z / 2) \Gamma(1-z) \Rightarrow \zeta(z)=f(z) \cdot \zeta(1-z)$, where $f(z)=2(2 \pi)^{z-1} \sin(\pi \cdot z / 2) \Gamma(1-z)$ but this also means that two cases occur:\\\\
\textbf{3.2.2.1.} If $\boxed{\boldsymbol{\zeta(1-z)=\zeta(z)}}$, when $z$ is a complex number.\\
This case is equivalent to 3.2.1.1, and therefore if $z=x_{0}+y_{0}$ i then $x_{0}=1 / 2$ in order to verify the equation $\zeta(z)=0$ we follow exactly the same process algebraically.\\\\
\textbf{3.2.2.2.} If $\boxed{\boldsymbol{\zeta(1-z) \neq \zeta(z)}}$, where $z$ is a complex number. Similarly, the above case is equivalent to 3.2.1.2 and therefore it cannot be happening, as it has been proved.\\\\
\textbf{4. Forms of the Riemann $\boldsymbol{\zeta}$ Functional Equations}\\
The Riemann $\zeta$ function has three types of zeros: [7]\\
\[\textbf{Eq. Set } \left\langle 
\begin{array}{cllrr}
\zeta(1-s)&=& 2 \cdot(2 \cdot \pi)^{-s} \cdot \cos (\pi \cdot s / 2) \cdot \Gamma(s) \cdot \zeta(s), &\text{Res}>0 \text{ \& }\text{Res}<0 \\
\zeta(s)&=& 2 \cdot(2 \cdot \pi)^{s-1} \cdot \sin (\pi \cdot s / 2) \cdot \Gamma(1-s) \cdot \zeta(1-s), &\text{Res}<1
\tag{I, II}
\end{array}\right.\]\\\\
And\\
\[\boxed{\boldsymbol{\zeta(s)=\pi^{s-1/2} \cdot \dfrac{\Gamma\left(\dfrac{1-s}{2}\right)}{\Gamma\left(\dfrac{s}{2}\right)} \zeta(1-s)} \tag{III}}\]\\
Usually referred to as the trivial zeros, and non-trivial complex zeros. Therefore has been proved that any non-trivial zero \textbf{lies in the open strip} $\{s \in \mathbb{C}: 0<Re(s)<1\}$ that is called the \textbf{critical Riemann Strip}: And all complex zeros of the function $\zeta$ lie in the line $\{s \in \mathbb{C}: Re(s)=1/2\}$ which is called the \textbf{critical line}.\\\\\\
\textbf{Solving of the Riemann $\boldsymbol{\zeta}$ Functional-Equations}.\\\\
To find the imaginary part we must solve the functional equations (Eq. set), and for the cases which the real part of the roots lies on the critical line Cases 4.1, 4.2 [6].\\\\
\textbf{4.1. 1\st type roots of the Riemann zeta functions (1\st equation from the Eq. set)}\\
For the first category roots and by taking the logarithm of two sides of the equations, and thus we get [5]:\\
\[\zeta(1-z) / \zeta(z)=2(2 \pi)^{-x} \cos (\pi \cdot z / 2) \Gamma(z) \Rightarrow\]
\[\Rightarrow \log [\zeta(1-z) / \zeta(z)]=\log [2]-z \log [2 \pi] + \log [\cos(\pi \cdot z / 2) \Gamma(z)]+2 k \pi i \]
and the total form from the theory of Lagrange inversion theorem, [5] for the root is $p_{1}(z)=s$ which means that $f(s)=p_{1}^{-1}(s)=z$, but with an initial value for $s$ which is
\[\boxed{\boldsymbol{s_{in} \rightarrow \dfrac{\log [2]}{\log [2 \pi]}+\dfrac{2 \cdot k \cdot \pi \cdot i}{\log (2 \pi)}}},\]
by setting values for $k$ we can simply calculate the roots of the (1\st equation from the Eq. set, (Table.3)).
\begin{table}[h!]
\centering
\begin{tabular}{|lll|}
 \hline
  \rowcolor{Gray}
 1\st  0.377145562795  & -/+ &\\
  \rowcolor{Gray}
 2\nd  0.377145562795  & -/+ & 3.41871903296 I\\
 \rowcolor{Gray}
 3\rd  0.3771455627955 & -/+ & 6.83743806592 I\\
 \rowcolor{Gray}
 4\nth 0.3771455627955 & -/+ & 10.2561570988 I\\
 \rowcolor{Gray}
 5\nth 0.3771455627955 & -/+ & 13.6748761318 I\\
 \rowcolor{Gray}
 6\nth 0.3771455627955 & -/+ & 17.0935951646 I\\[1ex]  
 \hline
\end{tabular}
\caption{Solution of 1\st equation from the Eq. set, The firsts $6$ Roots.}
\end{table}\\
\textbf{A simple program in mathematica is:}\\
\[\boxed{\def\arraystretch{1.5}
\begin{array}{ll}
k \geq 4 ; \text{ } t=(\log[2] + 2 k \pi i) / \log [2 \pi]; \\
z= s + \sum\limits_{w=1}^{\infty} \left( (-1 / \log [2 \pi])^{w} / Gamma[w+1] \right)+ D\left[(s^{\prime}) \cdot  \log [\cos[\pi \cdot s/2]] + \right.\\
+ \log [Gamma[s]]-\log [Zetta[1-s]/Zeta[s]])^{w}, \{s, w-1\} ];\\
N[z/.s \rightarrow t, 10]
\end{array}}\]\\
With $k=0, \pm 1, \pm 2, \ldots$\\\\
But because the infinite sum approaching zero, theoretically $s$ gets the initial value\\
\[\boxed{\boldsymbol{z \rightarrow \dfrac{\log [2]}{\log [2 \pi]}+\dfrac{2 \cdot k \cdot \pi \cdot i}{\log (2 \pi)}}}\]
So we have in this case, in part, the consecutive intervals with $k=n$ and $k=n+1$ for any $n \geq 4$ and for the imaginary roots.\\\\
\textbf{4.2. 2\nd type of roots of the Riemann zeta functional equations (2\nd equation from the Eq. set)}\\\\
Same as in the first category roots by taking the logarithm of two sides of the equations, and thus we get:
\[\zeta(z) / \zeta(1-z)=2(2 \pi)^{z-1} \sin(\pi \cdot z / 2) \Gamma(1-z) \Rightarrow\]
\[\Rightarrow \log [\zeta(z) / \zeta(1-z)]=\log [2]+(z-1) \log [2 \pi] +\log [\sin(\pi \cdot z / 2) \Gamma(1-z)] +2 \cdot k \cdot \pi \cdot i\]\\
Now, we will have for total roots of $3$ groups fields (but I have interest for the first group), and therefore for our case we will get $p_{1}(z)=s$ which means that:
\[f(s)=p_{1}^{-1}(s)=z,\]
but with an initial value for $s$ that is 
\[\boxed{\boldsymbol{s_{in} \rightarrow \dfrac{\log [\pi]}{\log [2 \pi]}+\dfrac{2 \cdot k \cdot \pi \cdot i}{\log (2 \pi)}}}\]
the overall form from the \textbf{Lagrange inverse theory} succeeds after replacing the above initial value and therefore for the \textbf{first six roots} after calculating them we quote the following (Table.4).
\begin{table}[h!]
\centering
\begin{tabular}{|lll|}
 \hline
 \rowcolor{Gray}
 1\st  0.622854437204  & -/+ &  \\ 
 \rowcolor{Gray}
 2\nd  0.622854437204  & -/+ & 3.41871903296 I \\
 \rowcolor{Gray}
 3\rd  0.622854437204  & -/+ & 6.83743806592 I \\
 \rowcolor{Gray}
 4\nth 0.622854437204  & -/+ & 10.2561570988 I \\
 \rowcolor{Gray}
 5\nth 0.622854437204  & -/+ & 13.6748761318 I \\
 \rowcolor{Gray}
 6\nth 0.622854437204  & -/+ & 17.0935951648 I \\[1ex] 
 \hline
\end{tabular}
\caption{The firsts $6$ Roots of 2\nd equation from the Eq. set}
\end{table}\\
\textbf{A simple program in mathematica is:}\\
\[\boxed{\def\arraystretch{1.5}
\begin{array}{ll}
k \geq 4 ; \text{ } t=(\log[\pi] + 2 k \pi i) / \log [2 \pi]; \\
z= s + \sum\limits_{w=1}^{\infty} \left( (-1 / \log [2 \pi])^{w} / Gamma[w+1] \right)+ D\left[(s^{\prime}) \cdot  \log [\sin[\pi \cdot s/2]] + \right.\\
+ \log [Gamma[1-s]]-\log [Zetta[s]/Zeta[1-s]])^{w}, \{s, w-1\} ];\\
N[z/.s \rightarrow t, 10]
\end{array}}\]\\
Therefore, with $k=0, \pm 1, \pm 2, \ldots$\\\\
But because the infinite sum approaches zero, theoretically $s$ gets initial value\\
\[\boxed{\boldsymbol{z \rightarrow \dfrac{\log [\pi]}{\log [2 \pi]}+\dfrac{2 \cdot k \cdot \pi \cdot i}{\log (2 \pi)}}}\]
So we have in this case, in part, the consecutive intervals with $k=n$ and $k=n+1$ for any $n \geq 4$ and for the imaginary roots. And for the cases we have for $Im(s)$ the relationship
\[Im(z)=\dfrac{2 \cdot \pi \cdot k}{\log [2 \pi]}, k \in Z \]
For the 3\rd functional equation, one has been previously put on the (as in the other two in 4, page 10) and we take as imaginary part
\[\boxed{ Re(z)=\dfrac{1}{2} \wedge Im(z)=\dfrac{2 \cdot \pi \cdot k}{\log [\pi]}, k \in Z}\]
Following to complete the roots of the two sets of the functional equations Eq. Set, we solve the functions as
cosine or sin \textbf{according to its Generalized theorem of Lagrange [5].}\\\\
\textbf{4.3. Transcendental equations for zeros of the function (Explicit formula)}\\\\
The main new results presented in the next few sections are transcendental equations satisfied by individual zeros of some L-functions. For simplicity, we first consider the Riemann-function, which is the simplest Dirichlet L-function.\\\\
\textbf{*Asymptotic equation satisfied by the n-th zero on the critical line.} $[9, 10, 11]$\\\\
As above, let us define the function $x(z) \equiv \pi^{-z / 2} \Gamma(z / 2) \zeta(z)$ which satisfies the functional equation $x(z)=x(1-z)$. Now consider the Stirling's approximation $\Gamma(z)=\sqrt{2 \pi} z^{z-1 / 2} e^{-z}\left(1+O\left(z^{-1}\right)\right)$ where $z=x+iy$, which is valid for large $y$. Under this condition, we also have 
\[ z^{z}=\exp \left(i\left(y \log y+\dfrac{\pi x}{2}\right)+x \log y-\dfrac{\pi y}{2}+x+O\left(y^{-1}\right)\right). \] 
Therefore, using the polar representation $\zeta=|\zeta| e^{i \arg \zeta}$ and the above expansions, we can write as 
\[\boxed{\def\arraystretch{2.3}
\begin{array}{l}
A(x, y)=\sqrt{2 \pi} \pi^{-x / 2}\left(\dfrac{y}{2}\right)^{(x-1) / 2} e^{-\pi y / 4}|\zeta(x+i y)|\left(1+O\left(z^{-1}\right)\right), \\
\theta(x, y)=\dfrac{y}{2} \log \left(\dfrac{y}{2 \pi e}\right)+\dfrac{\pi}{4}(x-1)+\arg \zeta(x+i y)+O\left(y^{-1}\right).
\end{array}}\]
The final transactions we end up with
\[n=\dfrac{y}{2 \pi} \log \left(\dfrac{y}{2 \pi e}\right)-\dfrac{5}{8}+\lim _{\delta \rightarrow 0^{+}} \dfrac{1}{\pi} \arg \zeta\left(\dfrac{1}{2}+\delta+i y\right).\]
Establishing the convention that zeros are labeled by positive integers, $z_{n}=1 / 2+i \cdot y_{n}$ where $n=1,2,3,4, \ldots$, we must replace $ n \rightarrow n-2$. Therefore, the imaginary parts of these zeros satisfy the transcendental equation
\[\textbf{Eq. a  } \text{  } \text{ } \boxed{\dfrac{t_{n}}{2 \pi} \log \left(\dfrac{t_{n}}{2 \pi e}\right)+\lim _{\delta \rightarrow 0^{+}} \dfrac{1}{\pi} \arg \zeta\left(\dfrac{1}{2}+\delta+i t_{n}\right)=n-\dfrac{11}{8}.}\]
Let us recall the definition used in, namely:
\[S(y)=\lim _{\delta \rightarrow 0^{+}} \dfrac{1}{\pi} \arg \zeta\left(\dfrac{1}{2}+\delta+i y\right)=\lim _{\delta \rightarrow 0+} \dfrac{1}{\pi} \Im \left[\log \zeta\left(\dfrac{1}{2}+\delta+i y\right)\right].\]
These points are easy to find, since they do not depend on the uctuating $S(y)$.\\\\
We have:
\[\zeta\left(\dfrac{1}{2}-i y\right)=\zeta\left(\dfrac{1}{2}+i y\right) G(y), \quad G(y)=e^{2 i \theta(y)}\]
The Riemann-Siegel \textit{\#} function is defined by
\[\vartheta \equiv \arg \Gamma\left(\dfrac{1}{4}+\dfrac{i}{2} y\right)-y \log \sqrt{\pi}.\]
Since the real and imaginary parts are not both zero, at $y_{n}^{(+)}$ then $G=1$, whereas at $y_{n}^{(-)}$ then $G=-1$. Thus
\[\boxed{\def\arraystretch{2.3}
\begin{array}{l}
\Im\left[\zeta\left(\dfrac{1}{2}+i y_{n}^{(+)}\right)\right]=0 \text{ for } \vartheta\left(y_{n}^{(+)}\right)=(n-1) \pi, \\
\Re\left[\zeta\left(\dfrac{1}{2}+i y_{n}^{(-)}\right)\right]=0 \text{ for } \vartheta\left(y_{n}^{(-)}\right)=(n-1) \pi. \\
\end{array}}\]
they can be written in the form of Eq.b.
\[\textbf{Eq. b } \text{ } \text{ } \boxed{y_{n}^{(+)}=\dfrac{2 \pi(n-7 / 8)}{W\left[e^{-1}(n-7 / 8)\right]}, \quad y_{n}^{(-)}=\dfrac{2 \pi(n-3 / 8)}{W\left[e^{-1}(n-3 / 8)\right]}}\]\\
\\\\
where above $n=1, 2, \ldots$ and the Labert - Function $W$ denotes the principal branch $W(0)$. The $y_{n}^{(+)}$ are actually the Gram points. From the previous relation, we can see that these points (Fig. 5) are ordered in a regular manner $[8, 12]$
\[y_{1}^{(+)}<y_{1}^{(-)}<y_{2}^{(+)}<y_{2}^{(-)}<y_{3}^{(+)}<y_{3}^{(-)}<\cdots\]
\begin{figure}[h!]
\centering
\includegraphics[scale=1.1]{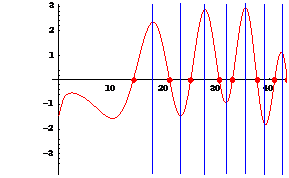}
\caption{Gram points}
\end{figure}\\\\
With this method, we are able to create intervals so as to approximate the correct values of imaginary part of non-trivial zeros. We did something similar in the cases 4.2. What it has left, is to find a method that approximates the values of imaginary part accurately.
\newpage
\begin{figure}[h!]
\centering
\begin{minipage}{0.7\textwidth}
\centering
\includegraphics[scale=1]{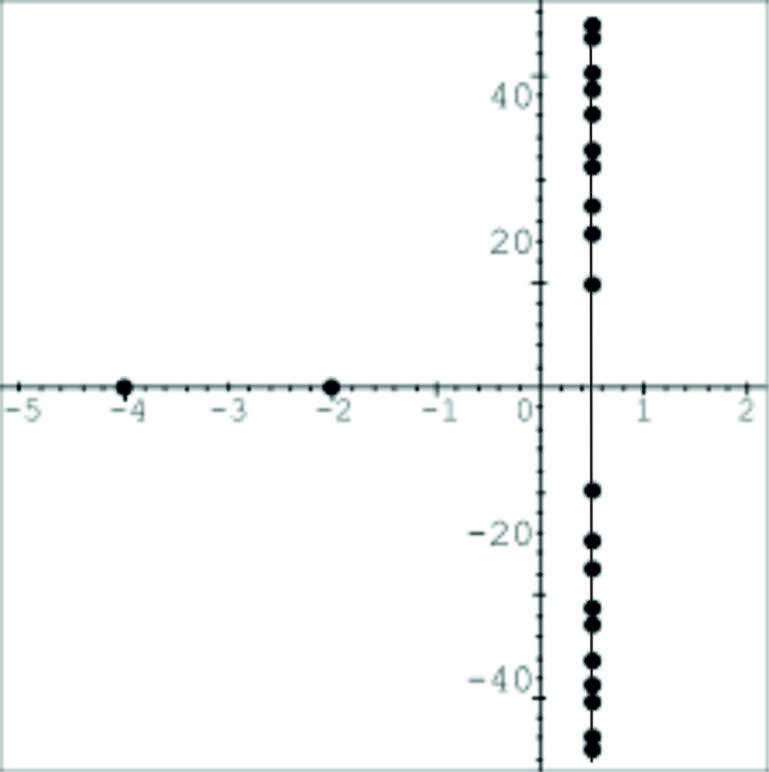}
\caption{The black dots represent the zeros of $\zeta(s)$ function including possible zeros which do not lie on the critical line.} 
\end{minipage}
\end{figure}
\vspace{3\baselineskip}
\begin{center}
\textbf{The zeros of the Riemann Zeta function.}\\
\textbf{“PROGRAMMING”}
\end{center}
\textbf{5.1. The M Function – Bisection Method}\\
Knowing the time of the successive steps $(k, k+1)$ of the relationship of imaginary parts $Im(z)=\dfrac{2 \cdot \pi \cdot k}{\log [2 \pi]}$, $Im(z)=\dfrac{2 \cdot \pi \cdot k}{\log [\pi]}$ (Eq I, II, III, page 10) with $k \in N$, we can calculate the roots by solving the equation $\zeta(1 / 2+i \cdot y)=0$ using the \textbf{Bisection Method}.\\\par
Bisection is the division of a given curve, figure, or interval into two equal parts (halves). A simple bisection procedure for iteratively converging on a solution, which is known to lie inside some interval $[a, b]$ proceeds by evaluating the function in question at the midpoint of the original interval $y=\dfrac{a+b}{2}$ and testing to see in which of the subintervals $[a,(a+b) / 2]$ or $[(a+b) / 2], b]$ the solution lies. The procedure is then repeated with the new interval as often as needed to locate the solution to the desired accuracy. Let $a_{n}, b_{n}$ be the endpoints at the n\nth iteration (with $a_{1}=a$ and $b_{1}=b$ ) and let $r_{n}$ be the n\nth approximate solution. Then, the number of iterations required to obtain an error smaller than $\varepsilon$ is
found by noting that... $b_{n}-a_{n}=\dfrac{b-a}{2^{n-1}}$ and that $r_{n}$ is defined by $r_{n}=\dfrac{1}{2}\left(a_{n}+b_{n}\right)$. In order for the error to be smaller than $\varepsilon$, then...
$\left|r_{n}-r\right| \leq \dfrac{1}{2}\left(b_{n}-a_{n}\right)=2^{-n}(b-a)<\varepsilon$. Taking the natural
logarithm of both sides gives $-n \ln 2<\ln \varepsilon-\ln (b-a)$. Therefore, we have for steps:
\[\boxed{n>\dfrac{\ln (b-a)-\ln \varepsilon}{\ln 2}.}\]
\textbf{5.2. M-function of the Bisection Method..}\\\\
We define the functions $M_{d}^{+}, M_{d}^{-}$ on an interval $(a, b)$ according to the scheme:\\\\
\textbf{I.} $M_{d}^{+}$, with $M_{d-1}$ the Nearest larger of $M_{d-m}$ where 
\[\dfrac{M_{d-1}+M_{d-m}}{2}, d \geq 2, d>m\geq d-2,\]
\[M_{1}=\dfrac{a+b}{2}, k=1\]
\textbf{II.} $\mathrm{M}_{\mathrm{d}}^{-}$, with $M_{d-1}$ the Nearest smaller of $M_{d-m}$ where
\[\dfrac{M_{d-1}+M_{d-m}}{2}, d \geq 2, d>m \geq d-2,\]
\[M_{1}=\dfrac{a+b}{2}, k=1\]
For calculating the roots on solving the equation $\zeta(1 / 2+i \cdot y)=0$ take the limit $r_{S}$ according to the scheme:
\[r_{s}=\lim _{n-\infty} M_{n}^{+}=\lim _{n-\infty} M_{n}^{-}, s \in N \text{ with } M_{d}^{+} \text{ and } M_{d}^{-},\]
belong in the interval $(a, b), d \in N$ and also $\boxed{\left(a=\dfrac{2 \cdot \pi \cdot k}{\log [2 \pi]}, \text{ } b=\dfrac{2 \cdot \pi \cdot(k+1)}{\log [2 \pi]}\right)}, k \in N, k \geq 3$ for the $x$ of $\zeta\left(1 / 2+i \cdot r_{s}\right)=0$.\\\\\\\\
\textbf{5.3 Program in Mathematica for the Bisection method of} $\zeta(1 / 2+i \cdot y)=0$.\\\\
Using the Intervals $\left(a=\dfrac{2 \cdot \pi \cdot k}{\log [2 \pi]}, \text{ } b=\dfrac{2 \cdot \pi \cdot(k+1)}{\log [2 \pi]}\right)$ and successive steps, we can compute all the roots of $\zeta(1 / 2+i \cdot y)=0$. We can of course use three types such intervals more specifically in generall than

\[\def\arraystretch{2.3}
\begin{array}{l}
1. a=\dfrac{2 \cdot \pi \cdot k}{\log [2 \pi]}, b=\dfrac{2 \cdot \pi \cdot(k+1)}{\log [2 \pi]} \text{ or}\\
2. a=\dfrac{2 \cdot \pi \cdot k}{\log [\pi]}, b=\dfrac{2 \cdot \pi \cdot(k+1)}{\log [\pi]} \text{ or}\\
3. a=y_{n}^{-}, b=y_{n}^{+}.
\end{array}\]

\[\boxed{y_{n}^{(-)}=\dfrac{2 \pi(n-7/8)}{W\left[e^{-1}(n-7/8)\right]}, \text{ } y_{n}^{(-)}=\dfrac{2 \pi(n-3/8)}{W\left[e^{-1}(n-3 / 8)\right]}}\]
We always prefer an interval that is shorter, in order to locate fewer non trivial zeros. The most important is to calculate all the roots in each successive interval and therefore only then we will have the program for data {example: Integer $k = 4$, and $\left(a=\dfrac{8 \cdot \pi}{\log [2 \pi]}, \text{ } b=\dfrac{10 \cdot \pi}{\log [2 \pi]}\right)$ and Error approximate $\text{tol }=10^{-6}$ and Trials $n=22$}..\\
\begin{center}
\textbf{“Programm for Bisection method”}
\end{center}
A program relevant by dividing intervals..\\\\
This program gives very good values as an approach to the roots we ask if we know the interval. Selecting the interval for the case $\boxed{\left(a=\dfrac{2 \cdot \pi \cdot k}{\log [2 \pi]}, \text{ } b=\dfrac{2 \cdot \pi \cdot(k+1)}{\log [2 \pi]}\right)}, k \in N, k \geq 3$ the results are given below in (Table.5)..
\newpage\noindent
Clear $[" * "]$;\\[4pt]
$f [ x_{-} ]:=$ Zeta $\left[ 1 / 2+ x ^* I \right]$;\\[4pt]
$k =\operatorname{Input}\left[\right.$"Epilogh $\left.k "\right]$;\\[4pt]
$a =\left(2^* k \right)^* \pi / \log [2 \pi]$;\\[4pt]
$b =2 *( k +1) * \pi / \log [2 \pi]$;\\[4pt]
tol=Input ["Enter tolerance"];\\[4pt]
$n =\operatorname{Input}[$ "Enter total iteration"];\\[4pt]
$g = N [(-\log [\operatorname{tol}] / \log [10]+4)]$;\\[4pt]
If$[ Arg [ f [ c ]]>0,\{$ Print["No solution exists"] $\}]$;\\[4pt]
Print $[" n$ a b c ........ $f(c) . . . . "]$;\\[4pt]
$Do [\{ c = N [( a + b ) / 2, g ]$;\\[4pt]
If $[$ Arg $[f[c]]<0, a=c, b=c]$, Print[PaddedForm[i,10],PaddedForm[N[a], 7,7$\}]$,\\[4pt]
PaddedForm $[ N [ b ],\{7,7\}],$ PaddedForm$[N [ c ],\{7,7\}]$,\\[4pt]
PaddedForm $[ N [ f [ c ]],\{7,7\}]]$ If $[ Abs [ a - b ]<$ tol \& \& Abs $[ N [ f [ c ]]]<$ tol*1000,\\[4pt]
\{Print["The solution is:",N[c,g]] Exit[]\}]\},\{i,1,n\}];\\[4pt]
Print ["The maximum iteration failed,No solution exists"];\\\\
\textbf{Calculation 1\st Zetazero} $z_1=\dfrac{1}{2}+14.1347216..$
\begin{table}[h!]
\centering
\includegraphics[width =7cm, height=7cm]{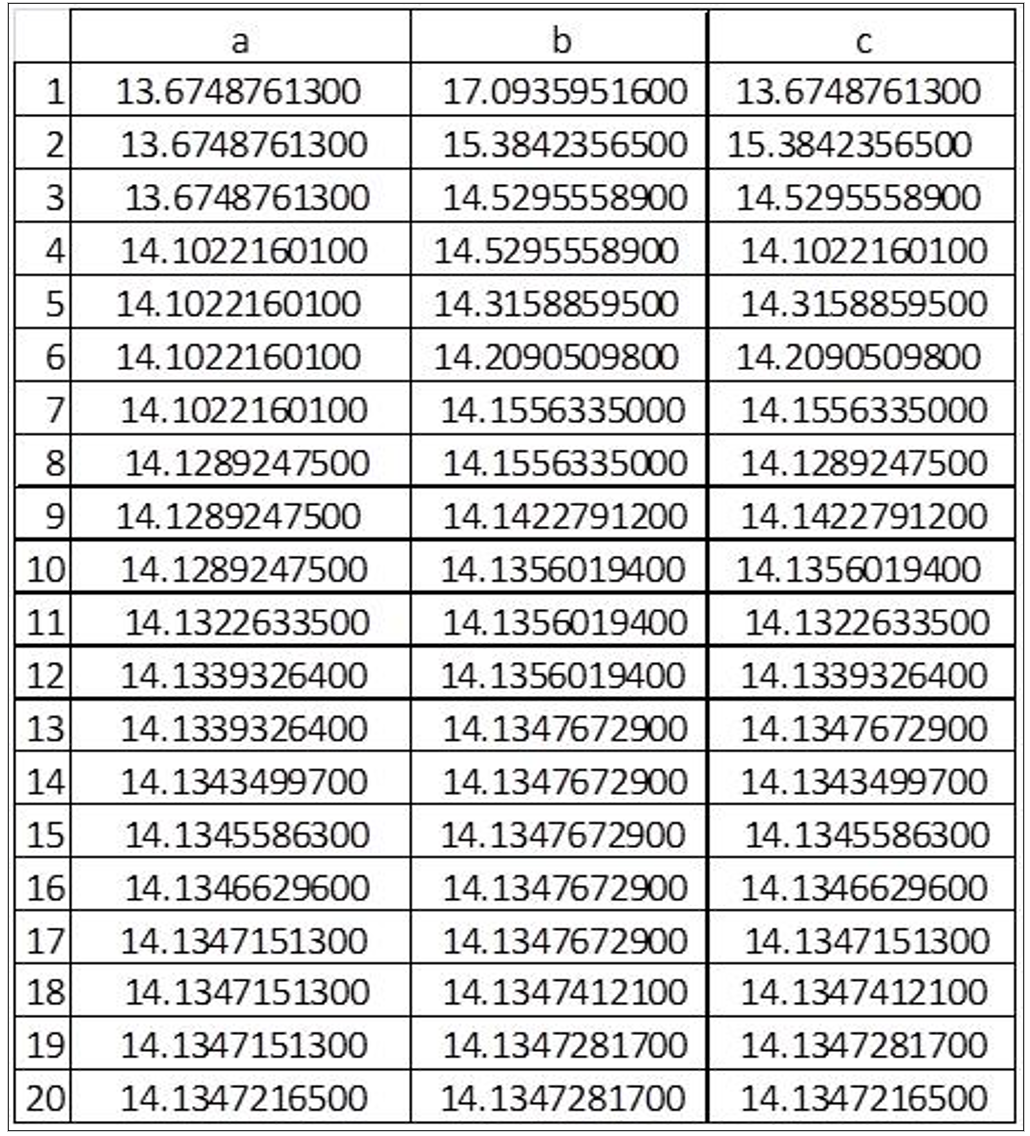}
\caption{The first non-trivial Root, Imaginary Part.} 
\end{table}\\
With final value $14.1347216500$, error near of $10^{-7}$. This value is the approximate root of the by the nearest error $<10^{-7}$.\\
In the event that we have two or more uncommon roots within the interval, we divide similar successive intervals in the order of finding of the first root either above or below. In such a case we have the $k = 13$, and at the interval $(47.8620664 , 51.280785)$ the two roots are $48.0051088$ and $49.7738324$.\\\\
\textbf{5.4 Explicit formula and the Zeros of} $\zeta(1 / 2+i \cdot y)=0$.\\\\
Consider its leading order approximation, or equivalently its average since $<\arg \zeta(1 / 2+i \cdot y)> \geq 0$ [11, 12].\\
Then we have the transcendental equation $\dfrac{t_{n}}{2 \pi} \log \left(\dfrac{t_{n}}{2 \pi e}\right)=n-\dfrac{11}{8}$. Through the transformation $t_{n}=2 \pi\left(n-\dfrac{11}{8}\right) x_{n}^{-1}$ this equation can be written a $x_{n} e^{x_{n}}=e^{-1}\left(n-\dfrac{11}{8}\right)$. Comparing the previous results, we obtain
\[\boxed{ t_{n}=\dfrac{2 \pi\left(n-\dfrac{11}{8}\right)}{W\left[e^{-1}\left(n-\dfrac{11}{8}\right)\right]}} \text{ where } n=1,2,3 \ldots\]\\
\textbf{5.5 Programm by Newton's method, which finds the Zeros of} $\zeta(1 / 2+i \cdot y)=0$\\\\
Using Newton's method we can reach the roots of the equation $\zeta(1 / 2+i \cdot y)=0$ at a very good initial value from $y \rightarrow t_{n}$ by the explicit formula.\\\\
It follows a mathematica program for the first 50 roots by the Newton's method (Table.6). This method determines
and detects the roots at the same time in order to verify the relation $\zeta(z)=0$ and always according to the relation... 
\[\boxed{t_{n}=\dfrac{2 \pi\left(n-\dfrac{11}{8}\right)}{W\left[e^{-1}\left(n-\dfrac{11}{8}\right)\right]}} \text{ where } n = 1,2,3 \ldots \text{ and where W is the W-function. }\]
Following is the program in mathematica and the Table.6 of complex roots…
\begin{table}[h!]
\centering
\includegraphics[width=9cm, height=6.6cm]{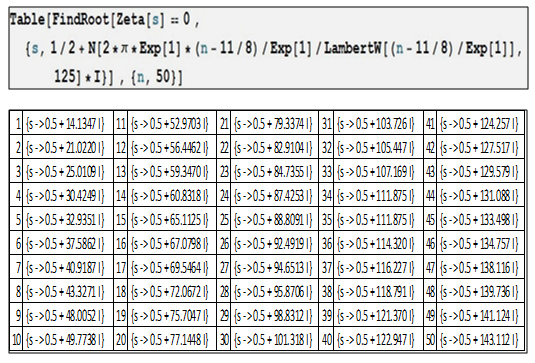}
\caption{The firsts 50 complex Roots(with W-function)} 
\end{table}\\
A very fast method that ends in the root very quickly and very close to the roots of $\zeta(s)=0$ as shown.\\\\
\textbf{5.6 Directly (from Explicit form) with the solution of this equation [8]}\\\\
We know that the equation
\[\boxed{\dfrac{t_{n}}{2 \pi} \log \left(\dfrac{t_{n}}{2 \pi e}\right)+\lim _{\delta \rightarrow 0^{+}} \dfrac{1}{\pi} \arg \zeta\left(\dfrac{1}{2}+\delta+i t_{n}\right)=n-\dfrac{11}{8}.}\]
Using the initial value the relation
\[\boxed{t_{n}=\dfrac{2 \pi\left(n-\dfrac{11}{8}\right)}{W\left[e^{-1}\left(n-\dfrac{11}{8}\right)\right]}} \text{ the explicit formula..} \] 
As already discussed, the function $\arg(\zeta (1/2+\delta+i \cdot y))$ oscillates around zero. At a zero it can be well-defined by the limit, which is generally not zero. For example, for the first Riemann zero 
$y_1=14.1347$ the limit $\delta \rightarrow 0^{+}$ has value as..
\[\lim \arg \zeta(1 / 2+\delta+i \cdot y)=0.157873.\]
The $\arg \zeta$ term plays an important role and indeed improves the estimate of the n\nth zero. We can calculate by Newton's method, and we locate the first $30$ imaginary part roots, of the equation $\zeta(1/2+i \cdot y)=0$, on the bottom of the (Table.7), and where $x=Im(1/2+yi)$. The process we use here does not cease to be approximate, although it gives us very good results quickly. The real solution to the problem, that is, the correct solution of the before transcendental equation, is dealt with most effectively in Chapter 5.7.
\begin{table}[h!]
\centering
\includegraphics[width=9.5cm, height=6.4cm]{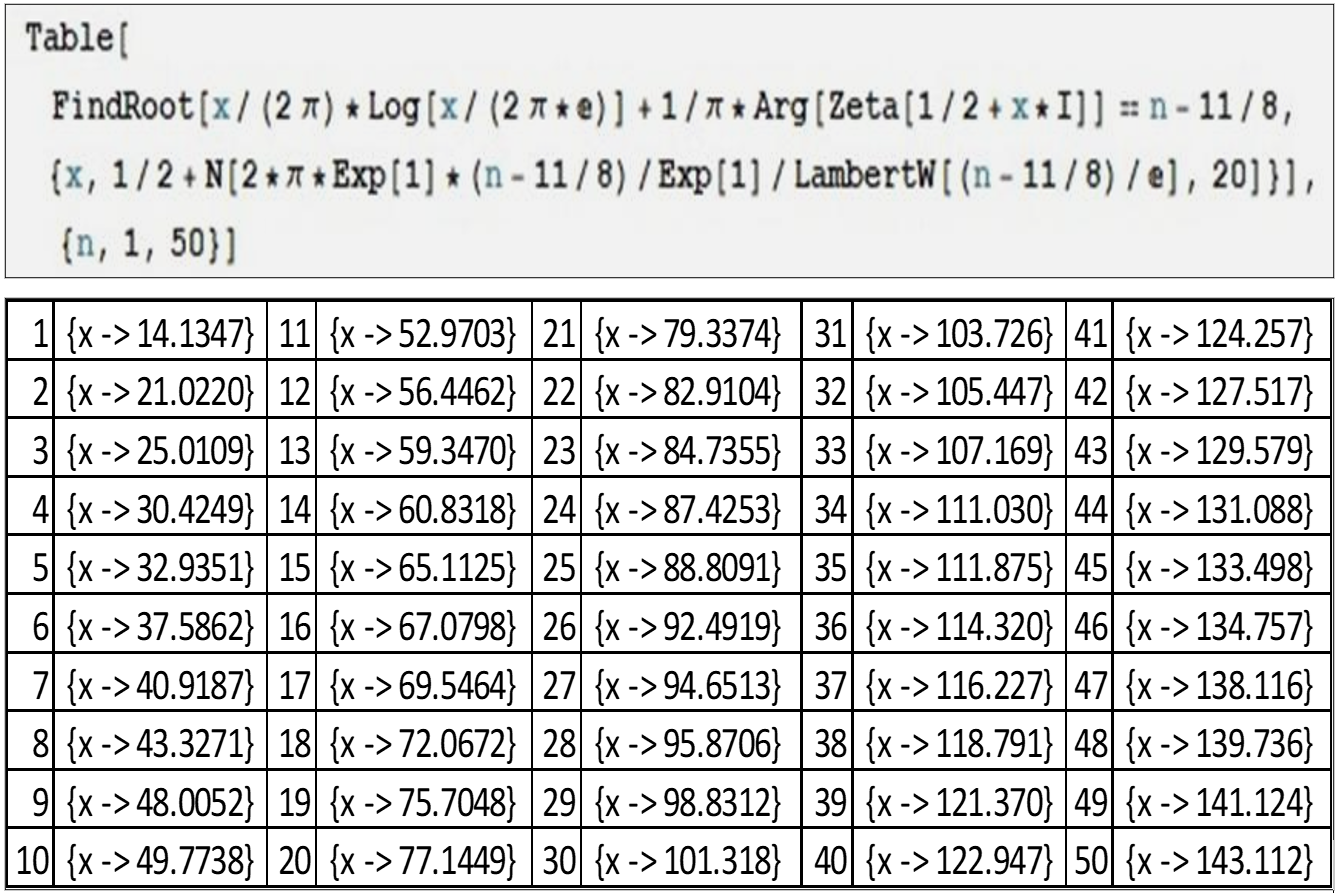}
\caption{The firsts 50 roots of ZetaZero.(Arg() Method)} 
\end{table}\\
\textbf{5.7 From Explicit form with the solution of method of Periodic Radicals..}\\\\
Using the relation Explicit form..
\[\boxed{\dfrac{t_{n}}{2 \pi} \log \left(\dfrac{t_{n}}{2 \pi e}\right)+\lim _{\delta \rightarrow 0^{+}} \dfrac{1}{\pi} \arg \zeta\left(\dfrac{1}{2}+\delta+i t_{n}\right)=n-\dfrac{11}{8}.}\]
We divide all the factors of the equation by $1/e$ and we take the more specialized form, according to the replace $x=t_n/(2 \pi e)$ therefore we have the form
\[\boxed{x \cdot \log(x)+1 /(\pi \cdot e) \cdot \arg(z(1/2+I \cdot(2 \cdot \pi \cdot e \cdot x))=1 / e \cdot(n-11/8)}\]
There are 2 ways to solve this equation or find the inverse of the expression $x\cdot\log(x)$ or of the expression $\arg(z(1 / 2+I \cdot(2 \cdot \pi \cdot e \cdot x)))$. Respectively the inversions are for each performance $\boxed{f_{1}(u)=e^{\text {Productlog(u)}}}$ and $\boxed{f_{2}(u)=\operatorname{Zeta}^{-1}\left(\arg^{-1}(u)\right)}$. The preferred procedure for finding a solution is the first because it is simpler, because the second it's more tricky program. Following is a program in mathematica.\\
\begin{figure}[h!]
\includegraphics[scale=0.355]{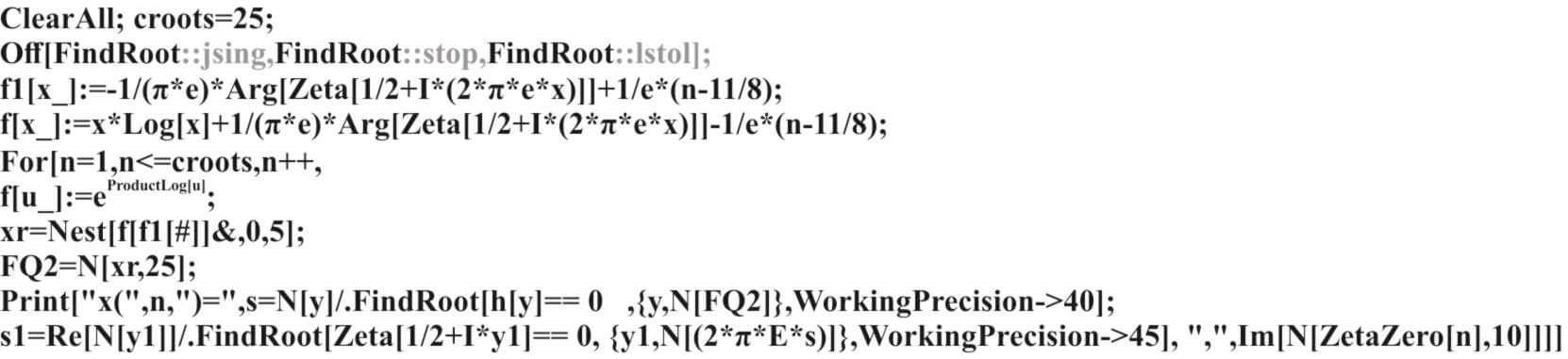}
\end{figure}\\
This program is based on the reversal of $x \log x$ term, the reverse of which is associated with W - Function, a process that gives us a multiplicity to the Improved programming, which we take advantage in order to find a plurality of roots. Those with the creation of this small Program we find it quickly and simply, any set of complex roots we need.
\begin{table}[h!]
\centering
\hspace*{6.5ex}
\includegraphics[width=8.8cm, height=8.5cm]{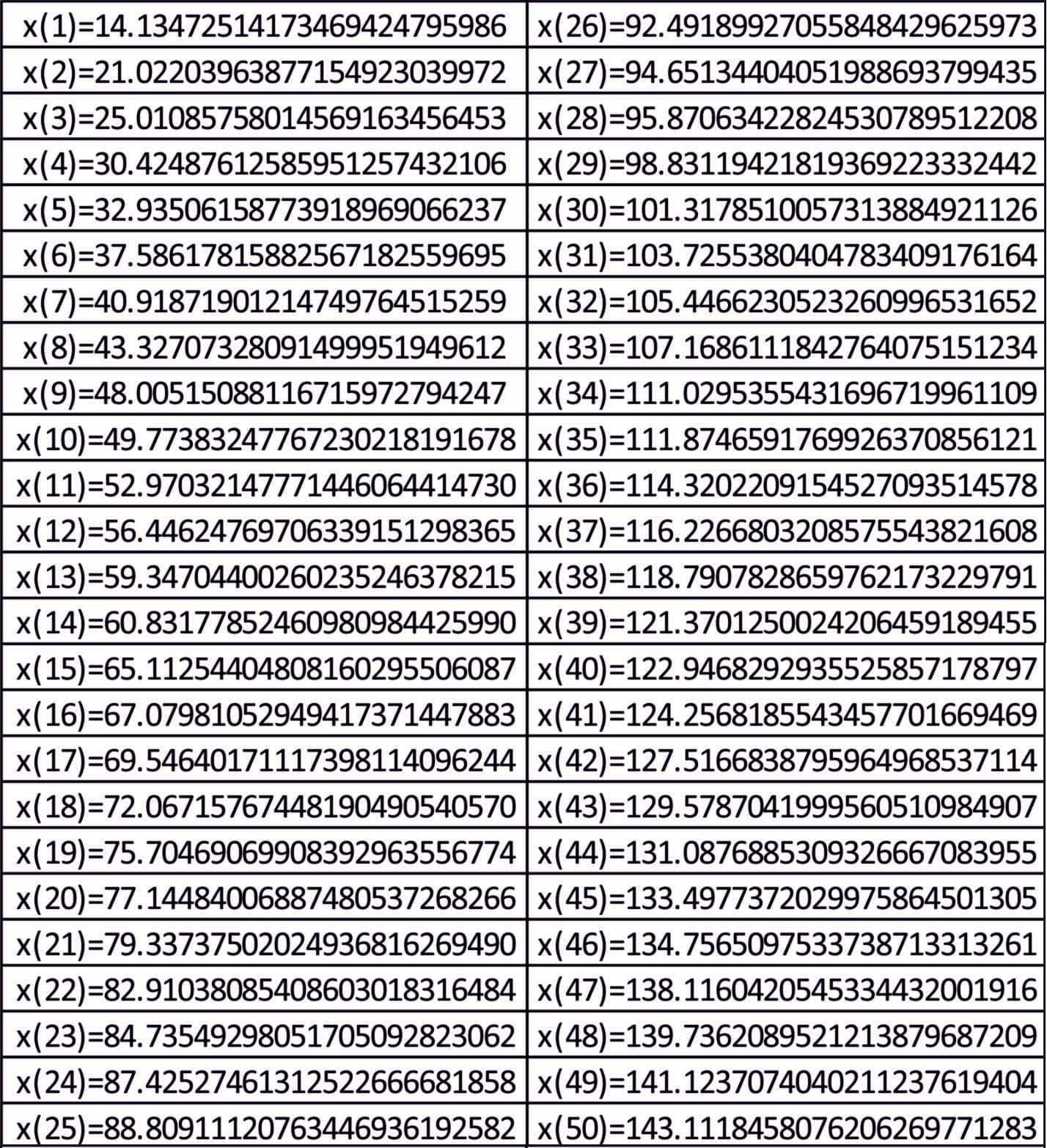}
\hspace*{0.8cm}\caption{The firsts $50$ roots of ZetaZero.($xlog(x)$ Method)} 
\end{table}
As we observe with this process we achieve a better approach and the roots are not removed from their correct values Perhaps it is the only case of finding complex roots in a large area and number that we could achieve computationally. Unlike case \textit{\#} 5.6 which is widely used based on the approximate type, 
\[\boxed{t_{n}=\dfrac{2 \pi\left(n-\dfrac{11}{8}\right)}{W\left[e^{-1}\left(n-\dfrac{11}{8}\right)\right]}}\] 
from explicit formula we may have discrepancies. That is, what is clearly seen in the calculation of complex roots of $\zeta(z)=0$ is the fact that they must treat the equation as a whole and not in parts. Finally, as we saw in this chapter, we fully tackled the transcendental equation as a whole and formulated without gaps-free algorithm with a simple process that works perfectly and efficiently.\\\\
\newpage
\noindent\textbf{6. }$\boldsymbol{100 \%}$\textbf{ of the zeros of $\boldsymbol{\zeta(z)=0}$ are on the critical line}\\\\
The elementary in our proof will be broken down into 5 parts.\\\\
\textbf{$\boldsymbol{6.1.1}$ Using direct method of inverse of $\boldsymbol{\zeta ()}$ on functional equations.}\\\\
From the proofs of Theorems $1,2,3$ we find it necessary to find the inverse of the function $\phi(z)=\zeta(z)-\zeta(1-z)$ by the lagrange method. First of all we have to mention why we got this equation now and we want to examine its inverse.  we saw in Theorem 3 since we proved from Th.1 that it is $1-1$ and therefore the function $\phi(z)$ is also inverted. We saw that the condition for the functional equations to hold is necessary in the form $\ldots, \phi(z)=0 .$ But from the fact that $\phi(z)$ is $1-1$ since the partial functions that make it up are also $1-1$. So $\zeta(z)=\zeta(1-z)=>z=1-z \Rightarrow \operatorname{Re}(z)=1 / 2$. Also a very simple example is the fact that the series
$1^{-z}+2^{-z}+3^{-z}+\ldots +n^{-z}=1^{z-1}+2^{z-1}+3^{z-1}+\ldots +n^{z-1}$ it is very simple to prove that take as solution the unique real $z=1 / 2$ comparing the terms one by one. \textbf{Therefore, according to Riemann's functional equations, we have one critical critical line on which the non-significant zeros will lie}.\\\\
\textbf{6.1.2 Solution of the functional Equation $\boldsymbol{\zeta(\mathrm{z})}$, $\boldsymbol{z}$ Complex number}\\\\
For real values with $\operatorname{Re}(z)>0$ and with form
$$
\begin{aligned}
\pi^{-\frac{x}{2}} \Gamma\left(\frac{z}{2}\right) \zeta(z) &=\pi^{-\frac{1-\nu}{2}} \Gamma\left(\frac{1-z}{2}\right) \zeta(1-z) \Rightarrow \\
\zeta(z) &=\pi^{s-\frac{1}{2}} \frac{\Gamma\left(\frac{1-g}{2}\right)}{\Gamma\left(\frac{g}{2}\right)} \zeta(1-z)
\end{aligned}
$$
if we take the logarithm but or $\mathrm{z}$ (we are making a change in the variable s), of the two parts of equation
[\#4. Eq-set(III)] then we have
$$
\zeta(z) / \zeta(1-z)=\pi^{z-\frac{1}{2}} 2(2 \pi)^{2-1} \Gamma\left(\frac{1-z}{2}\right) / \Gamma\left(\frac{z}{2}\right)=>
$$
$$
\log [\zeta(z) / \zeta(1-z)]=\left(\frac{1}{2}-z\right) \log [\pi]+\log \left[\Gamma\left(\frac{1-z}{2}\right) / \Gamma\left(\frac{z}{2}\right)\right]+2 k \pi i
$$
The resolution analysis will be performed by the Lagrange inversion theorem, which distinguished the best approach for transcendental equations. Using the correlation theory and after relations with 1 from 3 groups fields and therefore for our case we will get:\\\\
$p_{1}(z)=s$ which means that $f(s)=p_{1}^{-1}(s)=z$, but with an initial value\\\\
\[\mathrm{s}_{\mathrm{in}}=\frac{1}{2}+\frac{2 \cdot k \cdot \pi \cdot i}{\log (\pi)}\]
and \textbf{total form} from theory Lagrange for the root is:\\\\
\textbf{Program in mathematica }\\
\textbf{\bigg\{}Needs["\textit{NumericalCalculus}"]\\ $\mathrm{k}:=1 ; \mathrm{q}:=25 ;$\\\\
$z=N[\mathrm{s}_{\mathrm{in}}]+$\\\\
$\sum\limits_{w=1}^{q}\left(\left(-1/\log [\pi]\right)^{w}/\operatorname{Gamma} [w+1]\right)*$\\\\
$ND[(s')*\log [\operatorname{Gamma} [1-s]/2]/\operatorname{Gamma}[s/2]]-\log [\operatorname{Zeta}[s]/\operatorname{Zeta}[1-s]])^{w},\{s,w-1\},$\\\\
$t,\operatorname{WorkingPrecision}\rightarrow 250, \operatorname{Scale}\rightarrow .00001]$\\\\
$\mathrm{FQ}=\mathrm{N}[\mathrm{z} / .\mathrm{Sin}\rightarrow\mathrm{t}, 25]; \mathrm{N}\left[\pi \mathrm{FQ}-1 / 2^{*}\right.$ Gamma $[(1-\mathrm{FQ}) / 2] /$ Gamma $[\mathrm{FQ} / 2]$ - Zeta $[\mathrm{FQ}] /$ Zeta $\left.[1-\mathrm{FQ}]\right]$ \textbf{\bigg\}}\\\\
The results for $k=1$ and $q=25$ limit indicator of series is $0.5+5.48879293259 . .$ and with approximation $2.22044* 10^{-15}+2.664535* 10^{-15}$I of equation $\zeta$. Note : The program is in language mathematica
and where ND is $\frac{\mathrm{d}^{\mathrm{w}-1}}{\mathrm{~d} s^{\mathrm{w}-1}}$ and $s^{\prime}$ (is first derivative of $\left.\mathrm{s}\right)=1$, also $k \in N^{0}$.\\\\
Because of its limit summation tens to zero for a large $\mathrm{n}$, we can roots of the functional equation to write it simply as:
\[ z \rightarrow \frac{1}{2}+\frac{2 \cdot k \cdot\pi \cdot i}{\log (\pi)}, k \in Z\]\\
We see, therefore, that the roots of the functional equation of the $\zeta$, are on the critical line, with $\operatorname{Re}(z)=1 / 2$, a very serious result for its complex roots. Below we see a list with $k=0,1,2, \ldots, 10$ first values
which are the roots with a very good approach of the functional equation of the $\zeta()$ and if z a zeros we
have solutions:\\\\
\[\begin{aligned}
\{z==0.5, z= & =0.5+5.4887929325942295 i, z==0.5+10.977585865188459 i, z= \\[3pt]
& =0.5+16.46637879778269 i, z==0.5+21.955171730376918 i, z= \\[3pt]
& =0.5+27.443964662971148 i, z==0.5+32.93275759556538 i, z= \\[3pt]
& =0.5+38.42155052815961 i, z==0.5+43.910343460753836 i, z= \\[3pt]
& =0.5+49.399136393348066 i, z==0.5+54.887929325942295 i\}
\end{aligned}\]\\
If we want to see the roots of $\zeta (z) = 0$ ie analytically for the first 12 values:\\
\[\begin{aligned}
\{z==0.5+14.&134725141734695 i, z==0.5+21.022039638771556 i, z= \\[3pt]
&=0.5+25.01085758014569 i, z==0.5+30.424876125859512 i, z= \\[3pt]
&=0.5+32.93506158773919 i, z==0.5+37.586178158825675 i, z= \\[3pt]
&=0.5+40.9187190121475 i, z==0.5+43.327073280915 t, z= \\[3pt]
&=0.5+48.00515088116716 i, z==0.5+49.7738324776723 i, z= \\[3pt]
&=0.5+52.970321477714464 i, z==0.5+56.44624769706339 i\}
\end{aligned}\]\\
we can compare them to the previous values that we found of the functional equation of the $\zeta$. Below is given the graphical representation of  $\zeta$- function in C for $x\in (0,1)$ and $y\in (-31,31)$,
\newpage
\begin{figure}[h!]
    \centering
    \includegraphics[width=7.7cm, height=6.6cm]{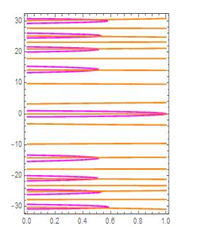}
    \caption{Plot of $\zeta (z)$ function}
\end{figure}
\vspace{\baselineskip}
\noindent we notice in the diagram Fig.7 that for the blue lines apply. $\operatorname{Im}\left(\frac{\operatorname{Zeta}(1 / 2+I \cdot y)}{\operatorname{Zeta}(1 / 2-I \cdot y)}\right)=0$ and for the Red
lines $\operatorname{Re}\left(\frac{\operatorname{Zeta}(1 / 2+I \cdot y)}{\operatorname{Zeta}(1 / 2-I \cdot y)}\right)=0$.\\\\\\
\textbf{6.1.3 Solutions of the Equation $\boldsymbol{\zeta (s) / \zeta(1-s)=1, s}$ is Complex number.}\\\\
A functional equation that hides all
secrets of the $\zeta(s)=0$. Based on the analytical extension of the zeta function we conclude that $\zeta (1-s)$ is analytical at $\mathrm{C}-\{0\}$ and zero will be a simple pole. So we will have:
$$
\begin{gathered}
\lim _{s \rightarrow 0} \zeta(1-s)=\infty \\
\lim _{s \rightarrow 0}-\frac{H_{\epsilon}(s)}{2 i \sin \frac{\pi}{2} s}(2 \pi)^{s}=\infty
\end{gathered}
$$
ie equality applies to $s_{0}$, where
$$
H(s)=\lim _{\epsilon \rightarrow 0} H_{\epsilon}(s)=-2 i \sin (\pi s) \Gamma(s) \zeta(s).
$$
As for the regard 4 categories roots of the equation $\mathrm{F}(\mathrm{s})=\zeta(\mathrm{s})-\zeta(1-\mathrm{s})=0$ and $s=x_{0}+y_{0} i$ we have a brief
summary of the form:
$$
\mathrm{s}=\left[\def\arraystretch{1.7}\begin{array}{ll}
\mathrm{x}_{0}=\frac{1}{2}, \pm \mathrm{y}_{0}, &\mathrm{if}\left(\mathrm{y}_{0} \neq 0\right) \\
\left\langle\begin{array}{l}
\mathrm{x}_{0}, \pm \mathrm{y}_{0} \mathrm{i} \\
1-\mathrm{x}_{0}, \pm \mathrm{y}_{0} \mathrm{i} \end{array}\right.
&\mathrm{if}\left(\mathrm{x}_{0} \neq \frac{1}{2}, \mathrm{y}_{0} \neq 0\right)
 \\
\left\langle\begin{array}{l}
\mathrm{x}_{0},\\
1-x_{0},
\end{array}\right.\quad &\text {if }\left(\mathrm{x}_{0} \neq \frac{1}{2}, \mathrm{y}_{0}=0\right) \\
\mathrm{x}_{0}=1-\mathrm{x}_{0}=\frac{1}{2}, &\mathrm{if}\left(\mathrm{y}_{0}=0\right)
\end{array}\right]
$$
The equation $\mathrm{F}(\mathrm{s})$ is related to the Hypothesis Riemann and helps in the final proof of Theorem 3(pages 8-10) seen before. It is a general form that also contains the non-trivial roots of $\zeta(s)=0$ and $\zeta(1-s)=0$.\\\\
Below we see a list with 15 first values, which are the roots with a very good approach of the equation of the $F(s)=0$ with method Newton in language Mathematica.\\\\
Reduce[Zeta[s]==Zeta[1-s]\&\& Abs[s]$<$15,s]\\\\
$s = \{ -12.16278-2.58046\text{ I}\}||s = \{ -12.1627878+2.5804649\text{ I}\}||s = \{ -7.9909145-4.5105941\text{ I}\}||s = \{ -7.99091453 + 4.51059414\text{ I}\}||s = \{ 0.5 - 14.13472514 \text{ I}\}||s = \{ 0.5 - 9.666908  \text{ I}\}||s = \{ 0.5 - 3.4362182   \text{ I}\}||s = \{ 1/2  \text{ I}\}||s = \{ 0.5 + 3.4362182    \text{ I}\}||s = \{ 0.5 + 9.66690805     \text{ I}\}||s = \{ 0.5 + 14.1347251      \text{ I}\}|| s = \{ 8.990914533 - 4.510594140      \text{ I}\}|| s = \{ 8.990914533 + 4.51059414069       \text{ I}\}||s = \{13.1627870 - 2.580464       \text{ I}\}||s = \{13.16278786 + 2.58046497       \text{ I}\}$\\\\
Further, is given the graphical representation of $F(s)=\zeta(s)-\zeta(1-s)$ functional equation.\\
\begin{figure}[h!]
    \centering
    \includegraphics[width=6.9cm, height=6.9cm]{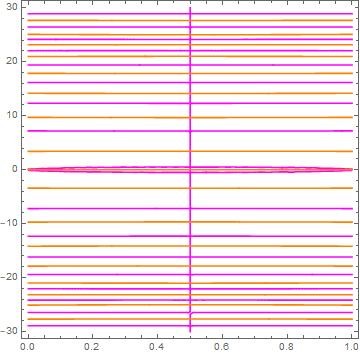}
    \caption{Plot of $\zeta(z)-\zeta(1-z)$ function}
\end{figure}\\\\
\textbf{6.2. With general inverse random function theory.}\\\\
According to the theory of inversion of functions that help us to solve them as equations, it completely guides us to a single solution i.e. $\operatorname{Re}(s)=1 / 2$ for $\zeta(s)=0$. In this article, we seek to develop a formula for an inverse Riemann zeta function $\mathrm{w}=\zeta(s)$ we have an inverse function $\mathrm{s}=\zeta^{-1}(\mathrm{w})$ and whics implies that $\zeta^{-1}(\zeta(s))=s \wedge \zeta^{-1}(\zeta(w))=w$ for real and complex Domain w and s. The presented method can also recursively compute these multiple solutions on other branches, but as we will find, the computational requirements become very high and start exceeding the limitations of the test computer, so we will primarily focus on the principal solution. We develop a recursive formula for an inverse Riemann zeta function as:
\[\left.s=\zeta^{-1}(w)=\lim _{m \rightarrow \infty} m\left[-\left.\frac{1}{(2 m-1)} \frac{d^{2 m}}{d s^{2 m}} \log \left[(\zeta(s)-w)(s-1)\right]\right|_{s\rightarrow 0} -\sum_{k=1}^{n} \frac{1}{s_{k}^{2 m}}\right)\right]^{-\frac{1}{2 m}}\]\\
Where $S_{n}$ is the multi-value for which $w=\zeta(s)$.
\newpage
\noindent\textbf{6.2.1. Theorem 4.} To prove that for the formula [26]\\\\
\[-\frac{1}{(m-1) !} \frac{d^{n}}{d z^{m}} \log |f(z)|_{z \rightarrow 0}=Z(m)-P(m)\]\\
that valid for a positive integer variable $m>1$. This formula generates $\mathrm{Z}(\mathrm{m})-\mathrm{P}(\mathrm{m})$ over all zeros and poles in the whole complex plane, as opposed to being enclosed in some contour.\\\\\\
\textbf{Proof.} If we model an analytic function $\mathrm{f}(\mathrm{z})$ having simple zeros and poles by admitting a factorization of the form
\[f(z)=g(z) \cdot \prod_{n=1}^{N_{2}}\left(1-\frac{z}{z_{n}}\right) \cdot \prod_{n=1}^{N_{p}}\left(1-\frac{z}{p_{n}}\right)^{-1}\]\\
Where $\mathrm{g}(\mathrm{z})$ is a component not any zeros or poles, then
\[\log [f(z)]=\log [g(z)]+\sum_{n=1}^{N_{z}} \log \left(1-\frac{z}{z_{n}}\right)-\sum_{n=1}^{N_{p}} \log \left(1-\frac{z}{p_{n}}\right)\]
From theorem Residue if i take
\[\frac{1}{2 \mathrm{i} \pi} \int \mathbf{P}_{\Omega} \frac{f^{\prime}(\mathrm{z})}{f(z)} \frac{1}{z^{s}} \mathrm{dz}=\mathrm{Z}(\mathrm{s})-\mathrm{P}(\mathrm{s})\]\\
And with use Taylor Theorem and $|\mathrm{z}|<1$ we take in final
\[\log [\mathrm{f}(\mathrm{z})]=\log [\mathrm{g}(\mathrm{z})]-\sum_{\mathrm{k}=1}^{\infty} \mathrm{Z}(\mathrm{k}) \frac{\mathrm{z}^{\mathrm{k}}}{\mathrm{k}}+\sum_{\mathrm{k}=1}^{\infty} \mathrm{P}(\mathrm{k}) \frac{\mathrm{z}^{\mathrm{k}}}{\mathrm{k}}\]
From this form, we can now extract $\mathrm{Z}(\mathrm{m})-\mathrm{P}(\mathrm{m})$ by the nth order differentiation as
\[-\frac{1}{(\mathrm{~m}-1) !} \frac{\mathrm{d}^{\mathrm{m}}}{\mathrm{d} z^{\mathrm{m}}} \log |\mathrm{f}(\mathrm{z})|_{\mathrm{z} \rightarrow 0}=\mathrm{Z}(\mathrm{m})-\mathrm{P}(\mathrm{m})\]\\\\
\textbf{6.2.2. Theorem 5.} If $Z_{n}$ is a set of positive real numbers ordered such that $0<z_{1}<z_{2}<\ldots<z_{n}$, and so
on, then the recurrence relation for the nth+1 term is
\[\mathrm{Z}_{\mathrm{n}+1}=\lim _{\mathrm{s} \rightarrow \infty}\left(\mathrm{Z}(\mathrm{s})-\sum_{\mathrm{k}=1}^{\mathrm{n}} \frac{1}{\mathrm{Z}_{\mathrm{k}}^{\mathrm{s}}}\right)^{-1 / \mathrm{s}}\]\\
\textbf{Proof.}\\\\
The generalized zeta series over zeros $\mathrm{Z}_{\mathrm{n}}$ of a function
\[\mathrm{Z}(\mathrm{s})=\sum_{n=1}^{\mathrm{N}_{z}} \frac{1}{\mathrm{z}_{\mathrm{n}}^{\mathrm{s}}}=\frac{1}{\mathrm{z}_{1}^{\mathrm{s}}}+\frac{1}{\mathrm{z}_{2}^{s}}+\frac{1}{\mathrm{z}_{3}^{\mathrm{s}}}+\ldots\]
and let us also assume that the zeros are positive, real, and ordered from smallest to largest such that $0<\mathrm{z}_{1}<\mathrm{z}_{2}<\ldots<\mathrm{z}_{\mathrm{n}}$ then the asymptotic relation holds
\[\frac{1}{z_{n}^{s}} > \frac{1}{z_{n+1}^{s}} \Rightarrow O\left(z_{n}^{-s}\right) \text{ } > \text{ }  O\left(z_{n+1}^{-s}\right),\text{ } s \rightarrow \infty\]
Final we take
\[\begin{aligned}
&\mathrm{z}_{1}=\lim _{s \rightarrow \infty}(\mathrm{Z}(\mathrm{s}))^{-1 / \mathrm{s}}, \text{ } \mathrm{z}_{2}=\lim _{s \rightarrow \infty}\left(\mathrm{Z}(\mathrm{s})-\frac{1}{\mathrm{z}_{1}^{s}}\right)^{-1 / \mathrm{s}}, \ldots, \\
&\mathrm{z}_{\mathrm{n}+1}=\lim _{\mathrm{s} \rightarrow \infty}\left(\mathrm{Z}(\mathrm{s})-\sum_{\mathrm{k}=1}^{\mathrm{n}} \frac{1}{\mathrm{z}_{\mathrm{k}}^{\mathrm{s}}}\right)^{-1 / \mathrm{s}},
\end{aligned}\]\\
thus all zeros up to the nth order must be known in order to generate the nth+1 zero. Exist the Voros's closed-form formula for depends on (RH), and the formula which we can express the zeta and beta terms in terms of a Hurwitz zeta function, and then substituting the the Voros's closed-form formula we obtain another formula for non-trivial zeros a for in terms of the Zeta .If the limit converges to $1 / 2$, it would imply (RH) is correct. \textbf{The formula of Voros's} is especially for real part (of $\rho_{\mathrm{n} . n t}$ the non-trivial zeros) of $\zeta(s)=0$ is
\[\begin{aligned}
\Re\left(\rho_{1, n t}\right)=& \sigma_{1}=\lim _{m \rightarrow \infty}\left[\left(\frac{1}{2} Z_{n t}^{2}(m)-\frac{1}{2} Z_{n t}(2 m)\right)^{-\frac{1}{s}}+\right.\\
&\left.-\left(\frac{(-1)^{m}}{2}\left(2^{2 m}-\frac{1}{(2 m-1) !} \log (|\zeta|)^{(2 m)}\left(\frac{1}{2}\right)-\frac{1}{2^{2 m}} \zeta\left(2 m, \frac{5}{4}\right)\right)\right)^{-\frac{1}{s}}\right]^{\frac{1}{2}}=\frac{1}{2}.\\
&&\text{V.roots}
\end{aligned}\]
\vspace{-\baselineskip}
\begin{figure}[h!]
\begin{minipage}{0.95\textwidth}
\begin{center}
\begin{tabular}{|c|c|}
\hline $\mathbf{m}$ & $\mathbf{R}(\boldsymbol{\rho} \mathbf{1}, \mathbf{n t})$ \\
\hline 15 & $0.4730925331369$ \\
\hline 20 & $0.4898729067547$ \\
\hline 25 & $0.4993065936936$ \\
\hline 50 & $0.5000000028549$ \\
\hline 100 & $0.4999999999999$ \\
\hline 150 & $0.5000000000000$ \\
\hline 200 & $0.4999999999999$\\
\hline
\end{tabular}
\end{center}
\end{minipage}\hfill
\begin{minipage}{0.05\textwidth}
\vspace{6.5\baselineskip}
\hspace{-38ex}
Table 9: Table Voro's. Approximation of\\
\hspace{-28.5ex}
 $Re(z)$ of $z ()$ function
\end{minipage}
\end{figure}
\phantom{ }\\
The computation of the real part of the first non-trivial zero for $m=200$ places (13 digits) using this relationship Vroots and with 200 iterations we achieve a very satisfactory result. It is of course obvious that it approximates by the direct inverse of $\zeta()$ , the real part of the non-trivial roots to  $1/2$.\\\\\\
\textbf{6.3. Statistical evidence}\\\\
Here, we formulate two of the functional equations from E-q.Set (Theorem 2, page 7)
\[\begin{aligned}
&\zeta(1-s) / \zeta(s)=2(2 \pi)^{-s} \operatorname{Cos}(\pi s / 2) \Gamma(s), \text{ }  s>0 \\
&\zeta(s) / \zeta(1-s)=2(2 \pi)^{s-1} \operatorname{Sin}(\pi s / 2) \Gamma(1-s), \text{ } \operatorname{Re}(s)<0 and \operatorname{Re}(s)<1
\end{aligned}\]
We look at each one equation individually in order to identify the set of values that we want each time. For the first equation and for real values with $\operatorname{Re}(s)>0$ and by taking the logarithm of two sides of the equation [5], we have...\\\\
real $s$ we will have $\boxed{\operatorname{Re}(s) \geq \dfrac{\log [2]}{\log [2 \pi]}=0.3771}$. This is the lower bound, which gives us the first $\zeta(s)$ of Riemann's Zeta Function.\\\\\\
For the second equation, for real values if $\operatorname{Re}(s)<1$ and by taking the logarithm of the two parts of the
equation, we will have:\\\\
In the following, because we need real $s$ we will take $\boxed{  \operatorname{Re}(s) \leq \frac{\log [\pi]}{\log [2 \pi]}=0.6228}$\\\\\\
average value of the upper lower bound is $\operatorname{Re}(\mathrm{s})=1 / 2$ because from that is apply \[\boxed{\operatorname{Re}(s)=1 / 2 \cdot \dfrac{\log [2]+\log [\pi]}{\log [2 \pi]}=1 / 2}\]\\\\
\textbf{6.4. The critical line also follows from the definition of $\boldsymbol{\lambda (N)}$}\\\\
Will define the function $\lambda(N)$ as follows:
\[\lambda\left(\mathrm{N}\right)=1+2^{-\mathrm{S}}+3^{-\mathrm{S}}+\ldots+\mathrm{N}^{-\mathrm{S}}\]\\
The next step is to perform a mathematical analysis on the functions of $\lambda$ after allocating $N$ value by equalizing it with $P_{n+1}$, because the values of $\lambda$ function in mainly depends on the prime numbers. In this way, the equations that represent the imaginary part of zeros which we are going to have will be as a function of the prime numbers. But the equations include all the natural numbers, so we'll generalize the solutions which depend on the prime numbers after obtaining all natural numbers with the same values which approach to those prime numbers.\\\\
Now we get the equation...\\\\
\[\dfrac{s(1-s) p_{n}^{-s}}{p_{n+1}^{1-s}}=\left(1-\dfrac{1}{2}\right)\left(1-\dfrac{1}{3}\right) \cdots\left(1-\dfrac{1}{p_{n-1}}\right) \dfrac{1}{p_{n}}\]\\\\
With replacing $1-s$ by $s$ and vice versa in the same equation, we get the\\
\[\dfrac{s(1-s) p_{n}^{s-1}}{p_{n+1}^{s}}=\left(1-\dfrac{1}{2}\right)\left(1-\dfrac{1}{3}\right) \cdots\left(1-\dfrac{1}{p_{n-1}}\right) \dfrac{1}{p_{n}}\]\\\\
By dividing equations the previous two, we get rid of the two factors $s$ and $(1-s)$ which cause the
emergence of the classical zeros of eta function\\\\
\[\begin{aligned}
&\dfrac{p_{n+1}^{s}}{p_{n}^{s-1}} \dfrac{p_{n}^{-s}}{p_{n+1}^{1-s}}=1 \Rightarrow \dfrac{p_{n+1}^{2 s-1}}{p_{n}^{2 s-1}}-1=0\\\\
&\Rightarrow \dfrac{p_{n+1}^{s-1 / 2}}{p_{n+1}^{s-1 / 2}}\left\{\left(\dfrac{p_{n+1}}{p_{n}}\right)^{s-\frac{1}{2}}-\left(\dfrac{p_{n+1}}{p_{n}}\right)^{-\left(s-\frac{1}{2}\right)}\right\}=0
\end{aligned}\]
\newpage
\noindent And therefore\\
\[\begin{aligned}
&\frac{1}{2 i}\left\{e^{i(s-1 / 2) \ln \left(\frac{P_{n+1}}{P_{n}}\right)^{-i}}-e^{-i(s-1 / 2) \ln \left(\frac{P_{n+1}}{P_{n}}\right)^{-i}}\right\}=0\\
&\Rightarrow \sin \left\{(s-1 / 2) \ln \left(\frac{P_{n+1}}{P_{n}}\right)^{-i}\right\}=0\\
&\Rightarrow-i(s-1 / 2) \ln \left(\frac{P_{n+1}}{P_{n}}\right)=\pm k \pi, k=1, 2, 3, \ldots\\
&\Rightarrow s=\frac{1}{2} \pm i \frac{k \pi}{\ln \left(\frac{P_{n+1}}{P_{n}}\right)}, k=1, 2, 3, \ldots\\
&&\text{(6.4.I)}
\end{aligned}\]\\
this last equation (6.4.I) proves the truth of the Hypothesis Riemann.\\\\\\
\textbf{6.4.1. A very important inequality arises from the 3rd functional equation for Number Theory}\\\\
From the 3rd Solve of the functional equations equation of $\zeta (s)$(page 10, III) we have for real values with $\mathrm{Re}(\mathrm{s})>0$ where $\mathrm{s}$ is Complex number and with form
$$
\begin{aligned}
\pi^{-\frac{x}{2}} \Gamma\left(\frac{s}{2}\right) \zeta(s) &=\pi^{-\frac{1-x}{2}} \Gamma\left(\frac{1-s}{2}\right) \zeta(1-s) \Rightarrow \\
\zeta(s) &=\pi^{s-\frac{1}{2}} \frac{\Gamma\left(\frac{1-s}{2}\right)}{\Gamma\left(\frac{g}{2}\right)} \zeta(1-s)
\end{aligned}
$$
if we take the logarithm but or $z$ (we are making a change in the variable s), of the two parts of equation
then we have
$$
\begin{aligned}
&\zeta(z) / \zeta(1-z)=\pi^{z-\frac{1}{2}} 2(2 \pi)^{z-1} \Gamma\left(\frac{1-z}{2}\right) / \Gamma\left(\frac{z}{2}\right) \Rightarrow \\
&\log (\zeta(z) / \zeta(1-z))=\left(\frac{1}{2}-z\right) \log (\pi)+\log \left(\Gamma\frac{1-z}{2}\right) / \Gamma\left(\frac{z}{2}\right)+2 k \pi i
\end{aligned}
$$\\
The resolution analysis will be performed by the Lagrange inversion theorem, which distinguished the best approach for transcendental equations. Using the correlation theory and after relations with 1 from
3 groups fields and therefore for our case we will get\\\\
$\mathrm{p}_{1}(\mathrm{z})=\mathrm{s}$ which means that $\mathrm{f}(\mathrm{s})=\mathrm{p}_{1}^{-1}(\mathrm{~s})=\mathrm{z}$ but with an initial value
\[\mathrm{s}_{\mathrm{in}}=\frac{1}{2}+\frac{2 \cdot \mathrm{k} \cdot \pi \cdot \mathrm{i}}{\log (\pi)}\]
And total \textbf{from Lagrange inversion theorem} for the root $\mathrm{Z}$ is if $\mathrm{q} \rightarrow \infty$
\[
\begin{aligned}
z=s_{i n}+& \sum_{w=1}^{q}\left((-1 / \log [\pi])^{w} / G a m m a[w+1]\right) \\
& \ast\left[\frac{d^{w-1}}{d s^{w-1}}(\log [(\operatorname{Gamma}[(1-s) / 2]) *(\operatorname{Gamma}[x / 2])]-\log [(\operatorname{Zeta}[x])(\operatorname{Zeta}[1-x])])^{w}\right.
\end{aligned}\]
Because of its limit summation tens to zero for a large $n$, we can roots of the functional equation to write
it simply as:
\[\mathrm{z} \rightarrow \frac{1}{2}+\frac{2 \cdot \mathrm{k} \cdot \pi \cdot \mathrm{i}}{\log (\pi)}, \mathrm{k} \in \mathrm{Z}\]
We see, therefore, that the roots of the functional equation of the $\zeta$, are on the critical line, with
$\operatorname{Re}(s)=1 / 2$, a very serious result for its complex roots. Below we see a list with $k=0,1,2, . .10$ first values,
which are the roots with a very good approach of the functional equation of the $\zeta$ and if $\mathrm{z}$ we have
solutions:
$$
\begin{aligned}
\{z==0.5, z=&=0.5+5.4887929325942295 t, z==0.5+10.977585865188459 i, z=\\
&=0.5+16.46637879778269 i, z==0.5+21.955171730376918 i, z=\\
&=0.5+27.443964662971148 i, z==0.5+32.93275759556538 i, z=\\
&=0.5+38.42155052815961 i, z==0.5+43.910343460753836 i, z=\\
&=0.5+49.399136393348066 i, z==0.5+54.887929325942295 t\}
\end{aligned}
$$
we can compare them to the previous values that we found of the functional equation of the $\zeta(z)$ equation. But from the second proof \textbf{the non-trivial zeros $\mathbf{z}$} has \textbf{a lower value} than the first one. This
means that that is: $z<s ($from $6.4 .1$) that is
\[\frac{2 \kappa \pi i}{\log (\pi)}<\frac{\kappa \pi i}{\log \left(\frac{p_{n+1}}{p_{n}}\right)} \Rightarrow \log \left(\frac{p_{n+1}}{p_{n}}\right)<\frac{1}{2} \log (\pi) \Rightarrow \frac{p_{n+1}}{p_{n}}<\sqrt{\pi}\]
Sure it is obvious
\[\boxed{\boldsymbol{\frac{p_{n+1}}{p_{n}}>1\Rightarrow 1<\frac{p_{n+1}}{p_{n}}<\sqrt{\pi}}\tag{6.4.2.II}} \]
The ratio of the largest to the smallest consecutive prime has an upper bound of the real number $\sqrt{\pi}$.
A relation that enters into many prime number conjectures.\\\\\\
\textbf{6.5. From the inverse of the function $\boldsymbol{\zeta ()}$ }\\\\
First we need to formulate the Riemann functional equation. It is necessary to define the function $\xi(\mathrm{s})$
\[\xi(s):=\frac{1}{2} s(s-1) \pi^{-s / 2} \Gamma\left(\frac{1}{2} s\right) \zeta(s)=(s-1) \pi^{-s / 2} \Gamma\left(\frac{1}{2} (s+1)\right) \zeta(s)\]
Where we have used the identity $\frac{1}{2} s \Gamma\left(\frac{1}{2} s\right)=\Gamma\left(\frac{1}{2} s+1\right)$.\\\\
We use the following theorem,\\\\
\textbf{Theorem 6.} The function $\xi$ has an analytic continuation to C. For continuation we have
$\xi(1-s)=\xi(s)$ for $s \in C$.\\\\
\textbf{6.5.2. Corollary 2.} The function $\zeta$ an analytic continuation to $\mathrm{C} \backslash\{1\}$ with a simple pole with Residue 1 at $\mathrm{s}=1$.\\\\
For this continuation we have
\[\zeta(1-s)=2^{1-s} \pi^{-s} \cos \left(\frac{1}{2} \pi s\right) \Gamma(s) \zeta(s) \text{ for } s \in C \backslash\{0,1\}\]\\
\textbf{Proof.} We define the analytic continuation of $\zeta$ by\\\\
\[\zeta(\mathrm{s})=\frac{\xi(\mathrm{s}) \pi^{s / 2} \frac{1}{\Gamma\left(\frac{1}{2} \mathrm{~s}+1\right)}}{\mathrm{s}-1}\]\\
$1/ \Gamma$ is analytic on $\mathrm{C}$ and the other function in the numerator are the analytic on C.Hence $\zeta$ is analytic on $\mathrm{C} \backslash\{1\}$. The  analytical continuation stated here coincides with that stated by the well-known theory $\{\mathrm{s}$ in $\mathrm{C}: \operatorname{Re}\text{ } \mathrm{s}>0\} \backslash\{1\}$ since analytic continuations to connected sets are uniquely determined. Hence $\zeta(\mathrm{s})$ has a simple pole with residue 1 at $\mathrm{s}=1$. We derive the functional equation. By Theorem $6.5 .1$ we have, for $\mathrm{s}$ in $\mathrm{C} \backslash\{1\}$,
\[\begin{aligned}
\zeta(1-s) &=\frac{\xi(1-s)}{\frac{1}{2}(1-s)(-s) \pi^{-(1-s) / 2} \Gamma\left(\frac{1}{2}(1-s)\right)}=\frac{\xi(s)}{\frac{1}{2} s(s-1) \pi^{-(1-s) / 2} \Gamma\left(\frac{1}{2}(1-s)\right)} \\
&=\frac{\frac{1}{2} s(s-1) \pi^{-s / 2} \Gamma\left(\frac{1}{2} s\right)}{\frac{1}{2} s(s-1) \pi^{-(1-s) / 2} \Gamma\left(\frac{1}{2}(1-s)\right)} \cdot \zeta(s)=F(s) \zeta(s)
\end{aligned}\]\\
Now we have\\
\[\begin{aligned}
&F(s)=\pi^{(1 / 2)-s} \cdot \frac{\Gamma\left(\frac{1}{2} s\right) \Gamma\left(\frac{1}{2} s+\frac{1}{2}\right)}{\Gamma\left(\frac{1}{2}-\frac{1}{2} s\right) \Gamma\left(\frac{1}{2}+\frac{1}{2} s\right)} \\
&=\pi^{(1 / 2)-s} \frac{2^{1-s} \sqrt{\pi} \Gamma(s)}{\pi / \sin \left(\pi\left(\frac{1}{2}-\frac{1}{2} s\right)\right)} \\
&=\pi^{-s} 2^{1-s} \cos \left(\frac{1}{2} \pi s\right) \Gamma(s).
\end{aligned}\]\\
This implies this Corollary.\\\\\\
\textbf{6.5.3. Direct solution of the equation $\boldsymbol{\zeta (q \ast z)=0}$, $\boldsymbol{q\in{R}}$}\\\\ Keeping the basic logic of the corollary 6.5.2., declare the definition $\mathrm{s} \mapsto \mathrm{q} \cdot \mathrm{z}, \mathrm{q} \in \mathrm{Z}$ who is analytic on $\mathrm{C}$, for every fixed $\mathrm{q} \in \mathrm{R}, \mathrm{q} \neq 0$. Then is obvious. According Corollary 6.5.2. we have \\
\[\zeta(1-\mathrm{q} \cdot \mathrm{z})=2^{1-\mathrm{q}\cdot \mathrm{z}} \pi^{-\mathrm{q}\cdot \mathrm{z}} \cos \left(\frac{1}{2} \pi \cdot \mathrm{q} \cdot \mathrm{z}\right) \Gamma(\mathrm{q} \cdot \mathrm{z}) \zeta(\mathrm{q} \cdot \mathrm{z}) \text{ for } \mathrm{q} \cdot \mathrm{z} \in \mathrm{C} \backslash\{0,1\}.\]\\
But according to theorems 1 \& 3 the function $\zeta()$ is 1-1 in the critical strip, which means that:\\
\[\zeta(1-\mathrm{q} \cdot \mathrm{z})=\zeta(\mathrm{q} \cdot \mathrm{z}) \text{ and therefore}\]
\[\mathrm{z}=\frac{1}{2 \mathrm{q}}, \mathrm{q} \neq 0, \mathrm{q} \in \mathrm{R} \text{ (a)}\]\\
This relationship always determines the position of the critical line, and of course atom is directly related to the $\zeta()$ function as we proved in Theorem 3. Let us now define a function ZetaZero()  which gives us \textbf{the zeroes} of $\zeta(\mathrm{s})=0$ as defined in sections $\{5.4,5.5,5.6,5.7\}$.\\\\
We therefore have the relation $\mathrm{s}=\operatorname{ZetaZero}(\pm \mathrm{k}), \mathrm{k} \in \mathrm{N}, \mathrm{k} \geq 1$. Therefore\\
\[\mathrm{s}=\operatorname{Re}( \text{ZetaZero}(\mathrm{k})) \pm \operatorname{Im}(\text{ZetaZero})(\mathrm{k}) \cdot \mathrm{i(b)},\text{ }\mathrm{i}\text{ is imaginary unit, }\mathrm{k} \in \mathrm{N}^{+}\]\\
As we know that for the simple function $\zeta(s)=0$ the\\
\[\operatorname{Re}(\mathrm{s})=\operatorname{Re}(\operatorname{ZetaZero}(\mathrm{k}))=\frac{1}{2}, \forall \mathrm{k} \in \mathrm{N}^{+}\]\\
From relations (a,b) we have the relation for \textbf{non trivial zeroes},
\[\boxed{\begin{aligned}
&\mathrm{z}_{\text {zeroes }}=\frac{1}{\mathrm{q}} \cdot \operatorname{Re}(\operatorname{ZetaZero}(\mathrm{k})) \pm \frac{1}{\mathrm{q}} \cdot \operatorname{Im}(\operatorname{Zeta} \operatorname{Zero}(\mathrm{k})) \cdot \mathrm{i} \Leftrightarrow \\
&\mathrm{z}_{\text {zeroes }}=\frac{1}{2 \cdot \mathrm{q}} \pm \frac{1}{\mathrm{q}} \cdot \operatorname{Im}(\text {ZetaZero }(\mathrm{k})) \cdot \mathrm{i}, \\
&\mathrm{i} \text { is imaginary unit}, \mathrm{k} \in \mathrm{N}^{+}, \mathrm{q} \in \mathrm{R}, \text { where } \mathrm{q} \neq 0.\quad\quad\quad \text{(c)}\\
\end{aligned}}\]
from the last relation (c) it follows that the Critical line will be given by the relation
\[\boxed{\begin{aligned}
&\mathrm{z}_{\mathrm{c} . \text {line }}=\frac{1}{\mathrm{q}} \cdot \operatorname{Re}(\operatorname{ZetaZero}(\mathrm{k})) \Leftrightarrow\\
&\mathrm{z}_{\mathrm{c} . \mathrm{line}}=\frac{1}{2 \mathrm{q}}, \forall \mathrm{k} \in \mathrm{N}^{+}, \mathrm{q} \in \mathrm{R}, \text { where } \mathrm{q} \neq 0.\quad \text{(d)}\\
\end{aligned}}\]
Finaly for the case of trivial zeroes we have the relation,
\[\boxed{\begin{aligned}
&z_{\text {zeroes }}=-\frac{2 \cdot \mathrm{k}}{\mathrm{q}}, \mathrm{k} \in \mathrm{N}^{+}, \mathrm{q} \in \mathrm{R}, \text { where } \mathrm{q} \neq 0. \quad \text{(e)}
\end{aligned}}\]
With this analysis we have essentially achieved the general solution of the transcendental equation
$\zeta(\mathrm{q} \cdot \mathrm{z})=0, \mathrm{q} \in \mathrm{R}^{+}$, and we observe that for positive values of $\mathrm{q}$, the critical lines, infinite in number, run
through the interval $(0,1)$ for 2 relations $(1 \& 2)$, i.e. $q>1/2$.\\
\[\boxed{\begin{aligned}
&1.\text{ }\text{ }0<\mathrm{z}_{\mathrm{c}\text{.line}}<\frac{1}{2} \text { then } 0<\frac{1}{2 \mathrm{q}}<\frac{1}{2} \text { therefore } \mathrm{q}>1 \\
&2.\text{ }\text{ }\frac{1}{2}<\mathrm{z}_{\mathrm{c}\text{.line}}<1 \text { then } \frac{1}{2}<\frac{1}{2 \mathrm{q}}<1 \text { therefore } \frac{1}{2}<\mathrm{q}<1 \\
&3.\text{ }\text{ }\mathrm{z}_{\mathrm{c}\text{.line}}>1 \text { then } \frac{1}{2 \mathrm{q}}>1 \text { therefore } 0<\mathrm{q}<\frac{1}{2}
\end{aligned}}\]\\
Similarly for the interval $(-1,0)$ we will have only 2 relations  $(4 \& 5)$, i.e. $q<-1/2$. If we want and $q=1/2$ then $q\leq -1/2$.\\
\[\boxed{\begin{aligned}
&4.\text{ }\text{ }-\frac{1}{2}<\mathrm{z}_{\text{c.line}}<0 \text { then }-\frac{1}{2}<\frac{1}{2 q}<0 \text { therefore } q<-1\\
&5.\text{ }\text{ }-1<\mathrm{z}_{\text {c.line }}<-\frac{1}{2} \text { then }-1<\frac{1}{2 \mathrm{q}}<-\frac{1}{2} \text { therefore }-1<\mathrm{q}<-\frac{1}{2}\\
&6.\text{ }\text{ }\mathrm{z}_{\mathrm{c}\text{.line}}<-1 \text { then } \frac{1}{2 \mathrm{q}}<-1 \text { therefore }-\frac{1}{2}<\mathrm{q}<0
\end{aligned}}\]\\\\\\\\
\textbf{6.5.4. Two Characteristic examples}\\\\
For the Riemann hypothesis the most important thing we are interested in is the interval $(0,1)$. The transcendental equation $z\left(q \ast z\right)=0$ as we have shown has infinite critical lines in this interval, which depend exclusively on ato $\mathrm{q}$.\\\\
i) I take for begin the mean value of interval $(0,1 / 2)$ and we have:\\\\
1. $z_{\text {c.line }}=\frac{1}{4}=\frac{1}{2 q} \Leftrightarrow q=2$\\\\
2. $\left\{\begin{array}{l}
\text {non trivials zeroes } \\ z_{\text {zeroes }}=\frac{1}{\mathrm{q}}(\text { ZetaZero }(\pm \mathrm{m})), \mathrm{m}=1 \div 5\\
z_{\text{zeroes}}=\{0.25+7.06736I,0.25+10.511I,0.25+12.5054I,\\
0.25+15.2124I,0.25+16.4675I\}\end{array}\right\}$\\\\\\
3. $\left\{\begin{array}{l}\text {trivials zeroes } \\ z_{\text {zeroes }}=-2 \times \frac{m}{q}, m=1 \div 5 \\ z_{\text {zeroes }}=\{-1 .,-2 .,-3 .,-4 .,-5 .\}\end{array}\right\}$\\\\\\
ii) Also the mean value of interval $(1 / 2,1)$ is $3 / 4$. Therefore we get the relations for this case:\\\\
1. $z_{\text {c.line }}=\frac{3}{4}=\frac{1}{2 q} \Leftrightarrow q=\frac{2}{3}$\\\\
$2 .\left\{\begin{array}{l}\text {non trivials zeroes } \\ z_{\text {zeroes }}=\frac{1}{\mathrm{q}}(\text { ZetaZero }(\pm m)), \mathrm{m}=1 \div 5 \\ z_{\text {zeroes }}=\{0.75+21.2021 \mathrm{I}, 0.75+31.533 \mathrm{I}, 0.75+37.5163 \mathrm{I} \\ 0.75+45.6373 \mathrm{I}, 0.75+49.4016 \mathrm{I}\}\end{array}\right\}$\\\\\\
3. $\left\{\begin{array}{l}\text {trivials zeroes } \\ z_{\text {zeroes }}=-2 \times \frac{\mathrm{m}}{\mathrm{q}}, \mathrm{m}=1 \div 5 \\ z_{\text {zeroes }}=\{-3 .,-6 .,-9 .,-12 .,-15 .\}\end{array}\right\}$\\\\\\\\
\textbf{6.6. Deniers of the Riemann Hypothesis.}\\\\
The most characteristic case that can be considered in relation to $\zeta(\mathrm{z})=0$ is the Davenport-Heilbronn function (introduced by Titchmarsh), is a linear combination of the two $L$ -functions with a complex character mod 5, with a functional equation of $L$ -function type but for which the analogue of the Riemann hypothesis fails. In this lecture, we study the Moebius inversion for functions of this type and show how its behavior is related to the distribution of zeros in the half-plane of absolute convergence. In general,
however, these cases do not refer to $\zeta(\mathrm{z})=0$ but to the $\mathrm{L}$ -function. The most important element in the negators of the Riemann hypothesis to consider is whether it is consistent with the Riemann equation $\zeta(\mathrm{z})=0$. Therefore, we will not consider cases that do not fall within our hypothesis, as defined by Riemann.\\\\\\
\textbf{6.6.1. Introduction of the denial conjecture through control one case of $\boldsymbol{\zeta(s,6)=0}$.}\\\\
However, there is one interesting case announced by the Enrico Bombieri .In a video by Enrico Bombieri[23] it is clear that there is a root that belongs to the
\[\zeta\left(\mathrm{s},\left(\mathrm{m}^{2}+5 \cdot \mathrm{n}^{2}\right)\right)=0\]
The root is $\mathrm{s}=0.93296997+15.668249531^{*} \mathrm{I}$ and is the prime for $\mathrm{m}=\mathrm{n}=1$. Thus the equation becomes
$\zeta(\mathrm{s}, 6)=0$ as he mentions. We have to prove that it is not the root of the equation.\\\\
True that $\zeta(\mathrm{s}, 6)=-0.101884-0.00598394I$ approximately. In contrast, the $\zeta(6, z)$ case gives a value to zero.
\newpage
\noindent In more detail approximately,\\
\[\zeta(6, \mathrm{z})=2.93166 \cdot 10^{-8}+2.10446 \cdot 10^{-7} \cdot \mathrm{I}\]\\
But if we do a simulation with the equation $\zeta\left(q^{*} z\right)=0$ then we need $q=1 /(2(\operatorname{Re}(z)))=1 /\left(2^{*} 0.93296997\right)$ i.e $\mathrm{q}=0.535923$ so we will have another root $\mathrm{z}=0.93297+26.3745 \mathrm{I}$ approximately, for which\\\\
\[\zeta(q \cdot z)=-2.46975 \cdot 10^{-15}+5.81423 \cdot 10^{-15} \cdot \mathrm{I}\]\\
a very good approximation and correct this time. To verify the fact that this root is a false root we will use 2 methods. One is with the \textbf{Newton method} locally and the other is to solve also locally the transcendental equation $\zeta(6, z)$ i.e. near the real root of the complex $\operatorname{Re}(\mathrm{s})=0.93296997$. If in this zone we do not find such a value this root will no longer exist.\\\\\\
\textbf{6.6.2. Approximate proof of the Newton method for checking $\boldsymbol{\zeta(6, \mathbf{z})}$.}\\\\
So assuming that the root is actually $\mathrm{s} 1=0.93296997-15.668249531^{*} \mathrm{I}$, we set it as the initial root in \textbf{FindRoot[Zeta $\boldsymbol{[6, \mathbf{z}],\{\mathbf{z}, \mathrm{s} 1\}}$, WorkingPrecision $\boldsymbol{\rightarrow\mathbf{5 0}]}$} with a dynamic iteration 50 times. We immediately see that it locally diverges to infinity and gives us a value
$$
\mathrm{z} \rightarrow 3.5898 \cdot 10^{7}-1.29617 \cdot 10^{9} \mathrm{I}
$$
This means that it failed near there to approach any root that would zero the equation $\zeta(6, z)$. In contrast, by the same method we observe that the local root $\mathrm{s} 1=0.93297+26.3745 \mathrm{I}$ for the equation $\zeta\left(\mathrm{q}^{*} \mathrm{z}\right)$ very quickly locates root i.e. using the command \textbf{FindRoot[Zeta $\boldsymbol{\left[q^{*} \mathbf{z}\right],\{\mathbf{z}, \mathrm{s} 1\} }$, WorkingPrecision $\boldsymbol{\left.->\mathbf{5 0}\right]}$} with $\mathrm{q}=0.535923$ then
$z=0.93296997+26.374547818I$ with approximation
$$
\zeta(\mathrm{q} \cdot \mathrm{z}) \approx 2.46975 \cdot 10^{-15}-1.29617 \cdot 10^{-15} .
$$
Therefore the root s1=0.93296997-I*15.668249531, \textbf{failed the first test} to confirm that it is a root.\\\\\\
\textbf{6.6.3. Direct proof of $\boldsymbol{\zeta(6, \mathbf{z})=0}$.}\\\\
In mathematics, \textbf{the polygamma function of order $\boldsymbol{m}$} is a meromorphic function on the complex numbers C, defined as the $(m+1)$ th derivative of the logarithm of the gamma function.\\\\
It satisfies the Reflection relation:
\[\boxed{(-1)^{\mathrm{m}} \psi(\mathrm{m}, 1-\mathrm{z})-\psi(\mathrm{m}, \mathrm{z})=\pi \cdot \frac{\partial^{\mathrm{m}}}{\partial \mathrm{z}^{\mathrm{m}}} \cot (\pi \cdot \mathrm{z})}\tag{6.6.3.1}\]
The polygamma function has the series representation
\[\boxed{\psi(\mathrm{m}, \mathrm{z})=(-1)^{\mathrm{m}+1} \mathrm{~m} ! \sum_{\mathrm{k}=0}^{\infty} \frac{1}{(\mathrm{z}+\mathrm{k})^{\mathrm{m}+1}}}\tag{6.6.3.2}\]
which holds for integer values of $m>0$ and any complex $z$ not equal to a negative integer. This representation can be written more compactly in terms of the \underline{Hurwitz zeta function} as\\
\[\boxed{\psi(\mathrm{m}, \mathrm{z})=(-1)^{\mathrm{m}+1} \mathrm{~m} ! \zeta(\mathrm{m}+1, \mathrm{z})}\tag{6.6.3.3}\]\\
Alternately, the Hurwitz zeta can be understood to generalize the polygamma to arbitrary, non-integer order. If we generally define $\zeta(\mathrm{m}+1, \mathrm{z})=\mathrm{f}(\mathrm{z})$ we can with improved programming solve this equation based on the 2 previous equations (6.6.3.1) \& (6.6.3.3). In this case we ask to find the roots when $\zeta(\mathrm{m}+1, \mathrm{z})=0$ therefore $\mathrm{f}(\mathrm{z})=0$ with method Periodic Radicals.\\\\\\
In this case we are interested in the fact that $\mathrm{\zeta}(6, \mathrm{z})=0$ and therefore $\psi(5, \mathrm{z})=0$ from relation (6.6.3.3). From (6.6.3.1) we have\\
\[\boxed{(-1)^{5} \psi(5,1-\mathrm{z})-\psi(5, \mathrm{z})=\pi \cdot \frac{\partial^{5}}{\partial \mathrm{z}^{5}} \cot (\pi \cdot \mathrm{z})}\tag{6.6.3.4}\]\\
But because $\Psi(5, \mathrm{z})=0$ the resulting final relation will be\\
\[\boxed{-\psi(5,1-z)=\pi \cdot \frac{\partial^{5}}{\partial z^{5}} \cot (\pi \cdot z)}\tag{6.6.3.5}\]\\
Performing the process of finding roots according to the Generalized Root Existence Theorem (of an equation) [24] we find the complementary function after solving the form\\\\
\[\boxed{\psi(5,1-z)+\pi^{6} \cdot\left(-16 \cdot(\cot (\pi z))^{4} \cdot(\csc (\pi z))^{2}-88 \cdot(\cot (\pi z))^{2} \cdot(\csc (\pi z))^{4}-16 \cdot(\csc (\pi z))^{6}\right)=0}\tag{6.6.3.6}\]\\
in final form\\\\
\[\boxed{\psi(5,1-z) / \pi^{6}+\left(-16 \cdot(\cot (\pi z))^{4} \cdot(\csc (\pi z))^{2}-88 \cdot(\cot (\pi z))^{2} \cdot(\csc (\pi z))^{4}-16 \cdot(\csc (\pi z))^{6}\right)=0\tag{6.6.3.7}}\]\\
To solve such a transcendental equation, we will use Periodic Radicals, which means we will split the equation into 2 parts. We solve the 2 nd part because it's a trigonometric equation and we have,\\
\[\left(-16 \cdot(\cot (\pi z))^{4} \cdot(\csc (\pi z))^{2}-88 \cdot(\cot (\pi z))^{2} \cdot(\csc (\pi z))^{4}-16 \cdot(\csc (\pi z))^{6}\right)=0\]\\
whose has in general 12 solutions $\mathrm{j}:=\pm 1 ; \mathrm{m}:=1 \div 6 ;$ as below,\\\\
\[z:=\frac{1}{\pi}\left(j \star 2 \operatorname{ArcTan}\left[\sqrt{\operatorname{Root}\left[15+30 \#1 +17 \#1^{2} +(4+8u) \#1^{3} +17 \#1^{4} +30 \#1^{5}  +15 \#1^{6} \&, m \right]}\right]+2 \pi k\right)\]\\\\
A radical relation which gives us the inverse of the 2nd part of equation 6.6.3.7.The whole secret is based on the parameter k of the trigonometric equation, which in turn will give us multiple different roots, which solve our basic equation $\zeta (6,z)=0$. The total cases we consider are 2X6=12, and of these we keep the values where $0< Re(z) <1$. With these zeros we make the following table of the program below.
\newpage
\noindent Then the matrix program in mathematica witn infinity Periodic Radicals is,
\begin{figure}[h!]
    \centering
    \includegraphics[width=14cm, height=6.7cm]{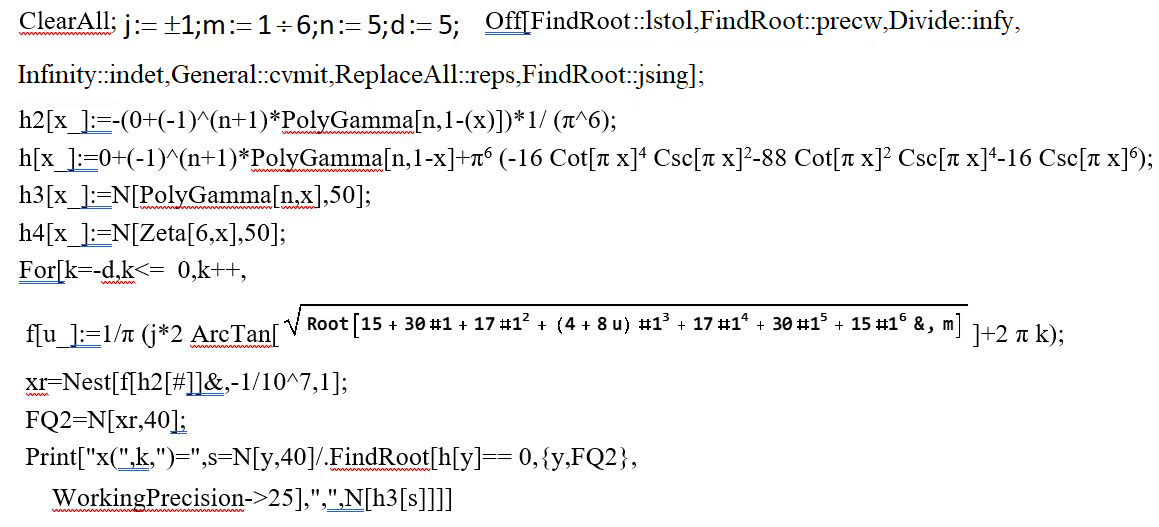}
\end{figure}\\
\hspace{-2.5ex}The final table of results is given below,
\begin{table}[h!]
\centering
\includegraphics[width=16.5cm, height=5.7cm]{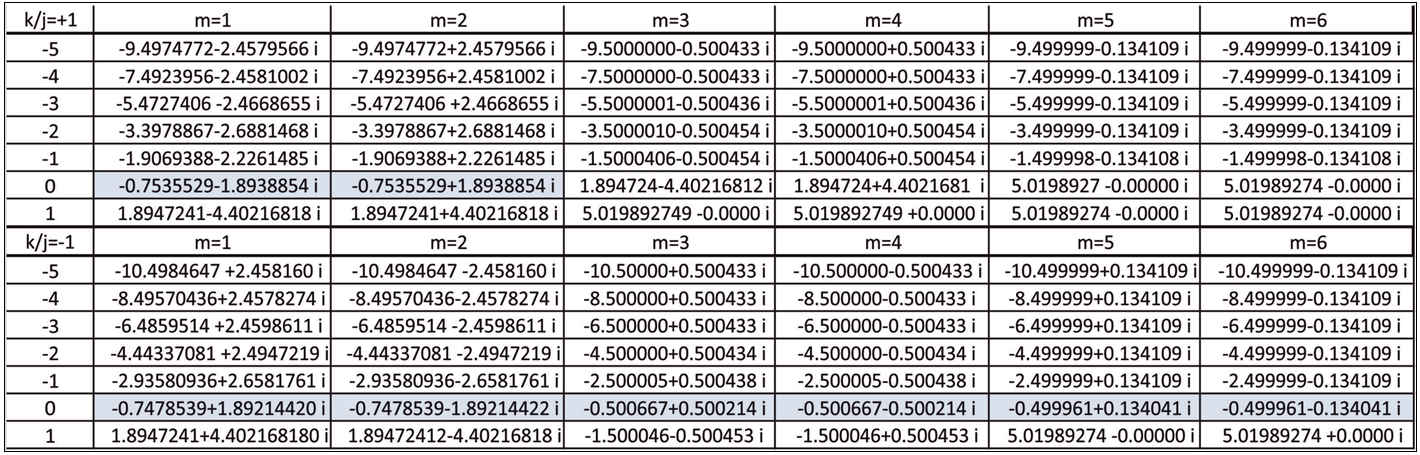}
\caption{Zeros of equation $\zeta(6,z)=0$} 
\end{table}\\
As we observe $\operatorname{Re}(\mathrm{z})<0$ and $\mathrm{k}<0$, for all roots. If we solve the equation for $\mathrm{k}>0$ we obtain exponential roots, i.e. large roots, were simply omitted. Therefore there is no solution for $Re(z)$ in the interval $(0,1)$ as a critical line.\\\\\\
\textbf{6.6.4. Direct proof of Epstein zeta equations $\boldsymbol{\zeta (Q(m,n),s)=0}$ and the contradiction of the deniers of R-H.}\\\\
In 1935 Potter and Titchmarsh for equation $\zeta(Q, s)$ proved has infinite zeros on the critical line and they also reported that a numerical calculation had shown that in the case, the Epstein function zeta has a zero on the critical strip but not on the critical line. This put an end to the speculation that the functional equation was the reason for the prevalence of a Riemann hypothesis. The deniers state "It is widely accepted that the generalized Riemann hypothesis is valid, with strong numerical proofs supporting it, but the proof remains long elusive." The Epstein zeta functions have the form of a double sum\\
\[\boxed{\zeta(\lambda, \mathrm{s})=\sum_{\mathrm{p}_{1}, \mathrm{p}_{2} \neq 0}^{\mathrm{n}} \frac{1}{\lambda\left(\mathrm{p}_{1}, \mathrm{p}_{2}\right)^{s}},\left\{\lambda\left(\mathrm{p}_{1}, \mathrm{p}_{2}\right)=\mathrm{a} \cdot \mathrm{p}_{1}^{2}+\mathrm{b} \cdot \mathrm{p}_{1} \mathrm{p}_{2}+\mathrm{c} \cdot \mathrm{p}_{2}^{2}, \mathrm{D}=\mathrm{b}^{2}-4 \cdot \mathrm{a} \cdot \mathrm{c}\right\},\left\{\mathrm{p}_{1}, \mathrm{p}_{2}, \mathrm{a}, \mathrm{b}, \mathrm{c} \in \mathrm{Z}+\right\}}\]\\\\
To complete the functions related to the $\zeta$-function we have to deal with the $\zeta(4, z), \zeta(5, z)$ which the deniers mention. We consider the sum:\\
\[\boxed{\zeta(\lambda, \mathrm{s})=\sum_{\mathrm{p}_{1}, \mathrm{p}_{2}}^{\mathrm{n}} \frac{1}{\left(\mathrm{p}_{1}^{2}+\mathrm{p}_{2}^{2} \cdot \lambda^{2}\right)^{s}}}\]\\
where the sum over the integers $\mathrm{p}_1$ and $\mathrm{p}_{2}$ runs over all integer pairs, apart from $(0,0)$, as indicated by the superscript prime. The quantity $\mathrm{Q}$ corresponds to the period ratio of the rectangular lattice, and $\mathrm{s}$ is an arbitrary complex number. For $\mathrm{Q}^{2}$ an integer, this is an Epstein zeta function, but for $\mathrm{Q}^{2}$ non-integer we will refer to it as a lattice sum over the rectangular lattice.\\
\[\mathcal{K}(n, m ; s ; \lambda)=\pi^{n} \sum_{p_{1}, p_{2}=1}^{\infty}\left(\frac{p_{2}^{s-1 / 2+n}}{p_{1}^{s-1 / 2-n}}\right) K_{s-1 / 2+m}\left(2 \pi p_{1} p_{2} \lambda\right).\]\\
For $\lambda>=1$ and the (possibly complex) numbers small in magnitude, such sums converge rapidly, facilitating numerical evaluations. (The sum gives accurate answers as soon as the argument of the MacDonald function exceeds the modulus of its order by a factor of $1.3$ or so.) The double sums satisfy the following symmetry relation, obtained by interchanging $\mathrm{p}_1$ and $\mathrm{p}_2$ in the definition:\\
\[\mathcal{K}(n,-m ; s ; \lambda)=\mathcal{K}(n, m ; 1-s ; \lambda).\]\\
Finally we get the relation for [25]\\
\[\begin{aligned}
&S_{0}(\lambda, s)=\zeta(Q, s)=\sum_{p_{1}, p_{2}}^{n} \frac{1}{\left(p_{1}^{2}+p_{2}^{2} \times \lambda^{2}\right)^{s}}, Q=p_{1}^{2}+p_{2}^{2} \cdot \lambda^{2} \\
&S_{0}(\lambda, s)=\frac{2 \zeta(2 s)}{\lambda^{2 s}}+2 \sqrt{\pi} \frac{\Gamma(s-1 / 2) \zeta(2 s-1)}{\Gamma(s) \lambda}+\frac{8 \pi^{s}}{\Gamma(s) \lambda^{s+1 / 2}} \mathcal{K}\left(0,0 ; s ; \frac{1}{\lambda}\right) .
\end{aligned}\]\\
From the general relations below we will derive the final relations of interest for the function $\mathrm{z}(\mathrm{Q}, \mathrm{s})$ :
\[\begin{aligned}
S_{0}(1, s) &=4 \zeta(s) L_{-4}(s), S_{0}(\sqrt{2}, s)=2 \zeta(s) L_{-8}(s) \\
S_{0}(\sqrt{3}, s) &=2\left(1-2^{1-2 s}\right) \zeta(s) L_{-3}(s), S_{0}(\sqrt{4}, s)=2\left(1-2^{-s}+2^{1-2 s}\right) \zeta(s) L_{-4}(s) \\
S_{0}(\sqrt{5}, s) &=\zeta(s) L_{-20}(s)+L_{-4}(s) L_{+5}(s), S_{0}(\sqrt{6}, s)=\zeta(s) L_{-24}(s)+L_{-3}(s) L_{+8}(s) \\
S_{0}(\sqrt{7}, s) &=2\left(1-2^{1-s}+2^{1-2 s}\right) \zeta(s) L_{-7}(s) .
\end{aligned}\]\\
From the zeroing of the first part of the transcendental representations for the cases where the values for $\lambda$
are $\lambda=\sqrt{3}, \lambda=\sqrt{4}, \lambda=\sqrt{7}$ and will we take:
\[\begin{gathered}
S_{0}(\sqrt{3}, s): s=\frac{1}{2}\left(1+\frac{(2 n+1) \pi i}{\ln 2}\right) \\
S_{0}(\sqrt{4}, s): s=\frac{1}{2} \pm \frac{i \arctan \sqrt{7}}{\ln 2}+\frac{2 n \pi i}{\ln 2} \\
S_{0}(\sqrt{7}, s): s=\frac{1}{2}+\frac{i \pi}{4 \ln 2}+\frac{2 n \pi i}{\ln 2}
\end{gathered}\]
We note that the critical line $1 / 2$ is technically constructed by solving the representations for these cases
and, of course, a number of infinite imaginary roots on the critical line. And in the case of the Epstein Function
we have points on the critical line. The only way to refute any negation for the $z(\mathrm{Q}, z)$ functions is to solve completely the 3 transcendental equations $\zeta(4, \mathrm{z})=0, \zeta(5, \mathrm{z})=0, \zeta(8, \mathrm{z})=0$ and check the results if any cases belong on the critical line.\\\\
Solving the other cases of Epstein functional equation approximately with the Periodic Radicals method.\\\\
If we assume that\\
$$
\boxed{
\begin{aligned}
&\begin{array}{lll}
Q(m, n)&=&a m^{2}+b m n+c n^{2} \\
&=&a(m+n \tau)(m+n \bar{\tau}) \\
&=&a|m+n \tau|^{2} 
\end{array}\\
&Z(s)=\sum_{m, n=-\infty}^{\infty} \frac{1}{a^{s}|m+n \tau|^{2 s}}, \sigma>1 .
\end{aligned}}
$$\\
in addition, the following relationships apply\\
$$
\begin{gathered}
\tau=\frac{b+\sqrt{b^{2}-4 a c}}{2 a} \\
\def\arraystretch{2.2}
\begin{array}{c}
x=\frac{b}{2 a}, \quad y=\frac{\sqrt{D}}{2 a}, \quad \tau=x+i y=\frac{b+i \sqrt{D}}{2 a} \\
\tau+\bar{\tau}=\frac{b}{a} \\
\tau \bar{\tau}=\frac{c}{a}
\end{array}
\end{gathered}
$$\\
separating the term with $\mathrm{n}=0$ and we have\\
\[
\boxed{
Z(s)=\frac{2}{a^{s}} \sum_{m=1}^{\infty} \frac{1}{m^{2 s}}+\frac{2}{a^{s}} \sum_{n=1}^{\infty} \sum_{m=-\infty}^{\infty} \frac{1}{|m+n \tau|^{2 s}}, \sigma>1}\tag{6.6.4.I}
\]\\
which is done\\
$$\sum_{m=-\infty}^{\infty} \frac{1}{|m+n \tau|^{2 s}}=\frac{\Gamma\left(s-\frac{1}{2}\right) \sqrt{\pi}}{n^{2 s-1} y^{2 s-1} \Gamma(s)}+\frac{4 \sqrt{\pi}}{n^{2 s-1} y^{2 s-1} \Gamma(s)} \sum_{m=1}^{\infty}(m n \pi y)^{s-\frac{1}{2}} \cos (2 m n \pi x) K_{s-\frac{1}{2}}(2 m n \pi y)$$\\
with use the Bessel function\\
$$K_{\nu}(y)=\frac{1}{\sqrt{\pi}}\left(\frac{2}{y}\right)^{\nu} \Gamma\left(\nu+\frac{1}{2}\right) \int_{0}^{\infty} \frac{\cos y t}{\left(1+t^{2}\right)^{\nu+\frac{1}{2}}} \mathrm{~d} t, \quad y>0, \quad \Re(v)>-\frac{1}{2},$$\\
and finally we get the equation for the Epstein function in the form\\\\
\resizebox{\textwidth}{!}{$\boxed{Z(s)=2a^{-s} \zeta(2 s)+2 a^{-s} y^{1-2 s} \frac{\Gamma\left(s-\frac{1}{2}\right) \sqrt{\pi}}{\Gamma(s)} \zeta(2 s-1)+\frac{8 a^{-s} y^{\frac{1}{2}-s} \pi^{s}}{\Gamma(s)} \sum_{n=1}^{\infty} n^{1-2 s} \sum_{m=1}^{\infty}(m n)^{s-\frac{1}{2}} \cos (2 m n \pi x) K_{s-\frac{1}{2}}(2 m n \pi y)}$}\\
\begin{flushright}
(6.6.4.II)
\end{flushright} 
It also turns out that is valid\\
\[\boxed{\left(\frac{\sqrt{D}}{2 \pi}\right)^{s} \Gamma(s) Z(s)=\left(\frac{\sqrt{D}}{2 \pi}\right)^{1-s} \Gamma(1-s) Z(1-s)}\tag{6.6.4.III}\]\\
which is the functional equation of the Epstein function. But as we proved [Theorem $1,2,3$ pages 1-10] with the function $\mathrm{\zeta}($ ) that it is $1-1$, the same is proved here by (6.6.4.3), which follows because $\mathrm{Z}(\mathrm{s})=\mathrm{Z}(1-\mathrm{s})$ i.e. $\mathrm{s}=1-\mathrm{s}=>\mathrm{s}=1 / 2$ for the Real part of the variable s. But we can also prove it from equation (6.6.4.II) if we consider the third part of the sum with a Bessel function as a small loss in the final result and because it is always imaginary, keeping the first 2 parts and obtain equation (6.6.4.IV) and work on it to find the roots.\\
\[\boxed{Z(s)=2 a^{-s} \zeta(2 s)+2 a^{-s} y^{1-2 s} \frac{\Gamma\left(s-\frac{1}{2}\right) \sqrt{\pi}}{\Gamma(s)} \zeta(2 s-1)}\tag{6.6.4.IV}\]\\
Thus we secure the real part which will show us the critical line and indeed uniquely. With the method we follow [24] (Generalized Theorem -Method Periodic Radicals) we find the generator (inverse function) of a part of the equation and solve the equation as a whole. For this equation a good generator is the term $\mathrm{y}^{1-2s}$. The analysis that will result after logarithm will give us multiple different roots, which is exactly what we have in such a case. Below we give the relationship of the generator,
\[
\def\arraystretch{2.2}\boxed{\begin{array}{l}
\mathrm{y}^{(1-2 \mathrm{x})}=\mathrm{u} \Rightarrow \mathrm{y}=\log (\mathrm{u})+2 \cdot \mathrm{k} \cdot \mathrm{i} \cdot \pi\\
\mathrm{x}=\frac{1}{2}-\left(\frac{\log (\mathrm{u})+2 \cdot \mathrm{k} \cdot \mathrm{i} \cdot \pi}{2 \cdot \log (\mathrm{y})}\right)
\end{array}}\tag{6.6.4.V}
\]\\
Then the matrix program in mathematica witn infinity Periodic Radicals is,
\begin{figure}[h!]
    \centering
    \includegraphics[width=14cm, height=6.7cm]{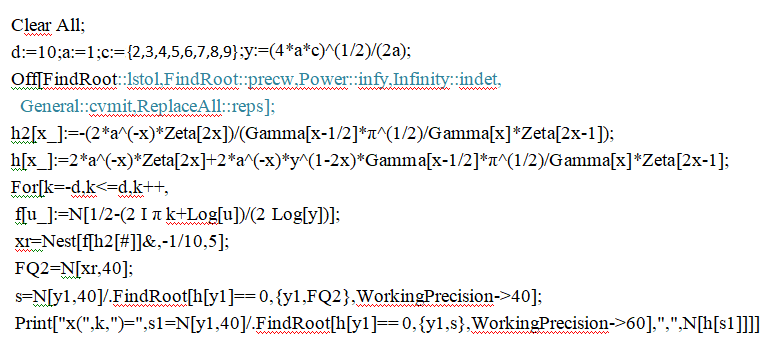}
\end{figure}\\
By running the program, all the roots with real part equal  1/2 are obtained. Here is a table with partial roots for $\lambda^{2} =2,3,4,5,6,7,8,9$.\\\\
Let's look in detail at the actual zeros per case in Table 10.
\newpage
\begin{table}[h!]
\centering
\includegraphics[width=12.98cm, height=8.89cm]{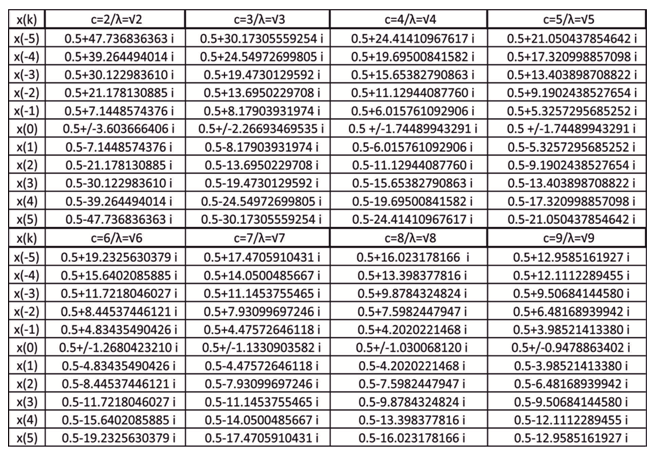}
\caption{Zeros of equation $Z(z)=0$ for $\lambda^{2}=2,3,4,5,6,7,8,9$} 
\end{table}
\vspace{2\baselineskip}
\noindent \textbf{6.6.5. Direct proof of $\boldsymbol{\zeta(4, z)=0}$.}\\\\
The whole procedure follows from chapter $6.3$ and equation $6.3 .3 .1$ It satisfies the Reflection relation:
\[\boxed{(-1)^{n} \psi(n, 1-z)-\psi(n, z)=\pi \cdot \frac{\partial^{n}}{\partial z^{n}} \cot (\pi \cdot z)}\tag{6.7.I}\]
this equation can be written making the first part positive as
\[\boxed{(-1)^{n+1} \psi(\mathrm{n}, 1-\mathrm{z})+\psi(\mathrm{n}, \mathrm{z})=-\pi \cdot \frac{\partial^{\mathrm{n}}}{\partial \mathrm{z}^{n}} \cot (\pi \cdot z)}\tag{6.7.II}\]
But since we associate the polygamma function and \underline{Hurwitz zeta function} arising
\[\boxed{\psi(\mathrm{n}, \mathrm{z})=(-1)^{\mathrm{n}+1} \mathrm{n} ! \zeta(\mathrm{n}+1, \mathrm{z})} \tag{6.7.III}\]
As a complementary relation we take the $\mathrm{n}$ -th derivative of $\cot (\pi z)$ which will be
\[\boxed{-\pi \cdot \frac{\partial^{\mathrm{n}}}{\partial \mathrm{z}^{\mathrm{n}}} \cot (\pi \cdot \mathrm{z})=\pi^{4}\left(4 \cdot \cot (\pi z)^{2} \csc (\pi \mathrm{z})^{2}+2 \csc (\pi z)^{4}\right)}\tag{6.7.IV}\]
We come therefore to the final system\\
\[\boxed{ 
\begin{array}{l}
(-1)^{n+1} \psi(\mathrm{n}, 1-z)+\psi(n, z)=\pi^{4}\left(4 \cdot \cot (\pi z)^{2} \csc (\pi z)^{2}+2 \csc (\pi z)^{4}\right)\\
\psi(n, z)=(-1)^{n+1} n ! \zeta(n+1, z)\\
\zeta(n+1, z)=0
\end{array}}  \tag{6.7.V}\]\\
Which takes the final form\\
\[ \boxed{\begin{array}{l}
(-1)^{n+1} \psi(n, 1-z)=\pi^{4}\left(4 \cdot \cot (\pi z)^{2} \csc (\pi z)^{2}+2 \csc (\pi z)^{4}\right)\\
\psi(n, z)=(-1)^{n+1} n ! \zeta(n+1, z)=0\\
\zeta(n+1, z)=0
\end{array}} \tag{6.7.VI}\]\\\\
The whole philosophy is therefore that if we admit that we zero $\zeta(\mathrm{n}+1, \mathrm{z})$ then automatically $\psi(\mathrm{n}, \mathrm{z})$ will also be zero. If we know a method that can solve $6.7$.II we can also solve $\zeta(\mathrm{n}+1, z)=0$. The method we follow here is the Infinite Periodic Radicals, which is considered here the \textbf{fastest and safest}, because we have high polynomial powers resulting from the multiple derivation of $\cot (\pi z)$ of the equation (6.7.II). In our case $\mathrm{n}=3$ to solve $\zeta(4, \mathrm{z})=0$.\\\\
We therefore have the following program,
\vspace{-0.5\baselineskip}
\begin{figure}[h!]
    \centering
    \includegraphics[width=14.5cm, height=9.5cm]{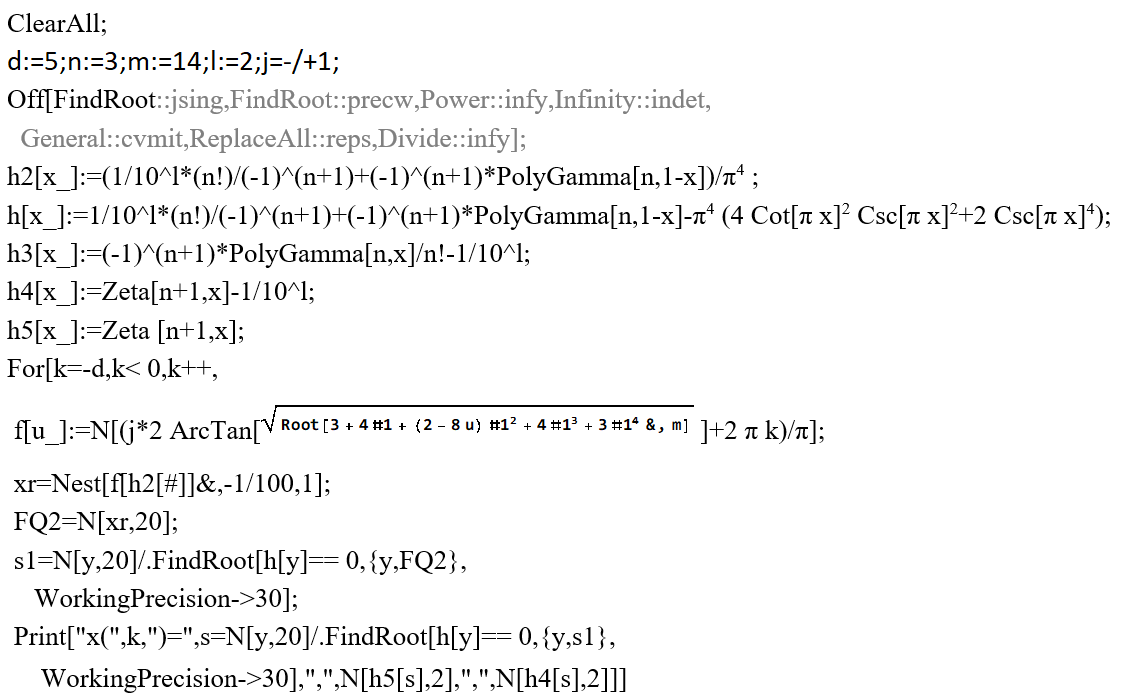}
\end{figure}
\vspace{-1\baselineskip}
\begin{table}[h!]
    \centering
    \includegraphics[width=16.5cm, height=5cm]{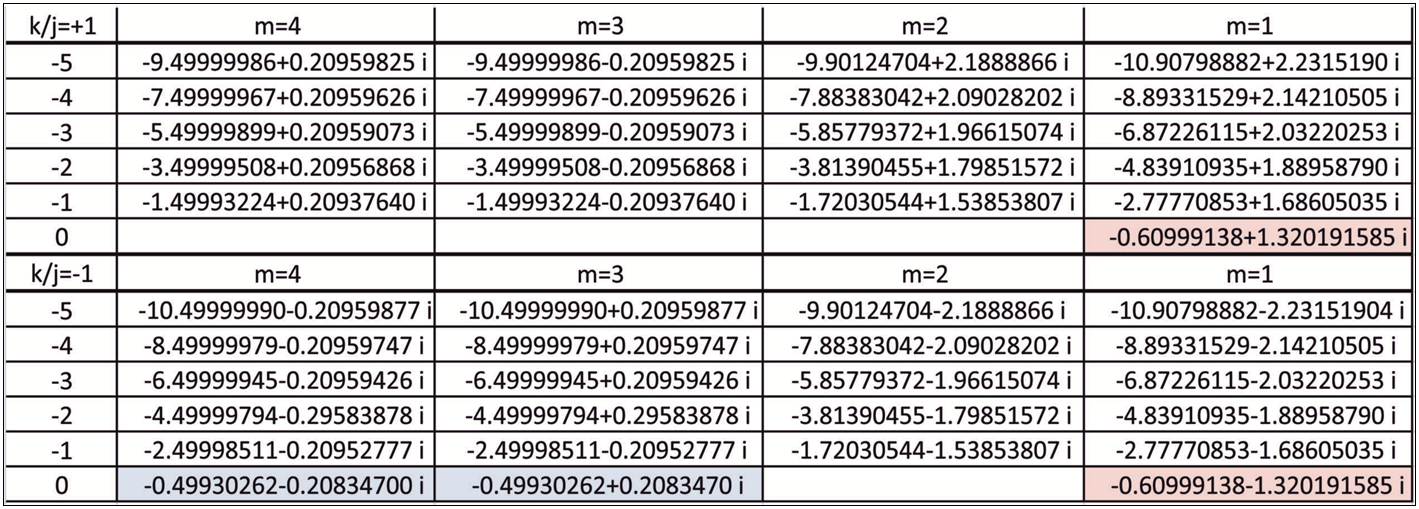}
    \caption{Zeros of equation $\zeta(4,z)=0$} 
\end{table}
\newpage
\noindent As we can see the program has $2$ variables $\mathrm{j}=+/-1$ and $\mathrm{m}=1 \div 4$ and gives us 8 cases that have been listed in the table 11, in detail for 5 consecutive roots. Moreover, as we can see, there is no critical line for these roots (because $\operatorname{Re}(\mathrm{z})<0)$, in the interval $(0,1)$ but not in the interval greater than 1. Also the cases of transcendental equations $\zeta(2, \mathrm{z})=0, \zeta(3, \mathrm{z})=0, \zeta(7, \mathrm{z})=0$ do not have the real part of complex zeros in the interval $(0,1)$ have been checked.\\\\\\
\textbf{6.6.6. The zeros of equations $\boldsymbol{\zeta(8,z)=0, \zeta(5,z)=0, \zeta(3,z)=0, \zeta(2,z)=0}$ and the tables of these.}\\\\
Based on the same theory and method of chapter $6.5$ we calculate the zeros around the interval $(0,1)$ below zero and above unity. If the values are not in this interval for the values of $Re(z)$ then it means that there are no zeros in the interval $(0,1)$ and therefore there are no critical lines other than $1/2$. From the given examination we see that there are all points outside the interval and per pair of points of conjugate complexes. They  have all in the interval $(-1,0)$ from $1$ pair upwards as q grows in the form $\zeta(q,z)$ that we consider. This means that as $q$ increases in value, the probability of finding a critical line in the interval $(0,1)$ is zero. The process of developing complex pairs is given on page 42 and they show a peculiarity around $-1/2$ and $-1$.\\\\
\phantom{ }\\
\textbf{I. Zeros of $\boldsymbol{\zeta(8,z)=0}$.}\\
\begin{table}[h!]
    \centering
    \includegraphics[width=16.5cm, height=5.7cm]{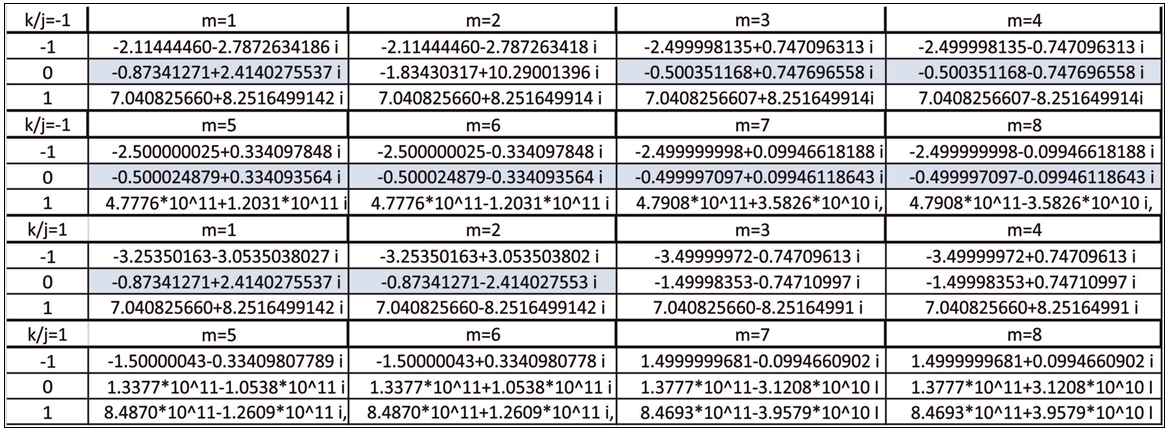}
    \caption{Zeros of equation $\zeta(8,z)=0$} 
\end{table}\\
\noindent\textbf{II. Zeros of $\boldsymbol{\zeta(5,z)=0}$.}\\
\begin{table}[h!]
    \centering
    \includegraphics[width=16.5cm, height=3cm]{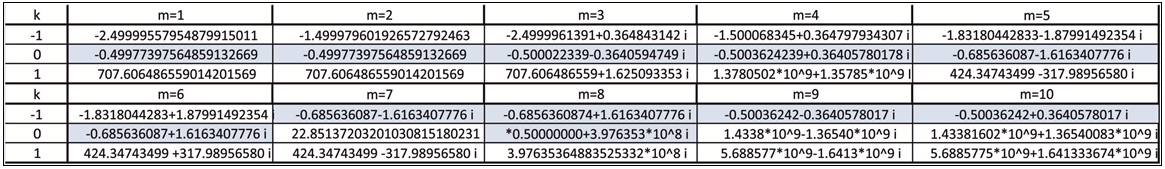}
    \caption{Zeros of equation $\zeta(5,z)=0$} 
\end{table}
\newpage
\begin{flushleft}
\textbf{III. Zeros of $\boldsymbol{\zeta(3,z)=0}$.}\\
\end{flushleft}
\begin{table}[h!]
    \centering
    \includegraphics[width=15.1cm, height=1.73cm]{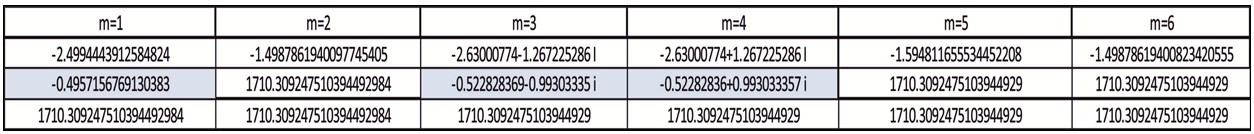}
    \caption{Zeros of equation $\zeta(3,z)=0$} 
\end{table}
\phantom{ }\\
\textbf{IV. Zeros of $\boldsymbol{\zeta(2,z)=0}$.}\\
\begin{table}[h!]
    \centering
    \includegraphics[width=13cm, height=1.72cm]{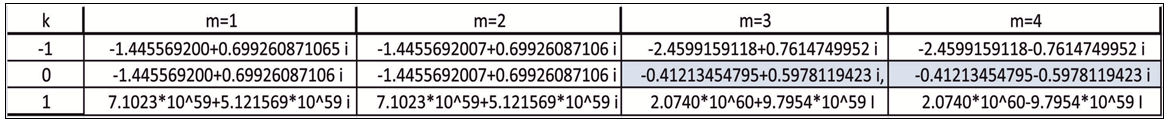}
    \caption{Zeros of equation $\zeta(2,z)=0$} 
\end{table}\\\\
Preceded by the given tables with zeros near the interval $(0,1)$ laterally in order of magnitude. The colored zeros are those closest to the interval $(0,1)$; notice immediately that they are divided into $3$ categories. Those where $Re(z)$ is around $-1/2$ and those where $Re(z)$ is close to $-1$. Following the logic of the matrices we construct the following Conclusion, which is general for the class of functional equations $z(q,z)$, with $q>1$ in $Z+$.\\\\\\\\
\textbf{6.6.7. Corollary 3.}\\\\
\textbf{The zeros of the Generalized Riemann Zeta function $\boldsymbol{\zeta(q,z)=0, q>1}$ in $\boldsymbol{Z+}$, do not belong to the critical strip of the interval $\boldsymbol{(0,1)}$, but are outside.}\\\\
According to the Generalized root finding theorem [24] and the Infinite Periodic Radicals method, by which we solved the whole series of equations with $q={2,3,4,5,6,8}$ we observe that the roots, are divided into $3$ categories, $2$ in those where is around $-1/2$ and which approach it as $q$ increases, and another one which moves away from $-1/2$ towards $-1$ as $q$ increases. Of course we are looking at these roots which are their real parts in the interval $(0,1)$.\\\\
We could also consider the case for $q$ in $R$ in general, and such a brief examination will be done in section 6.7 where we will consider the generalized equation $G-\zeta (z,q)$, i.e. we will also see some instances of $\zeta (q,z)$ with positive and negative $q$. Thus we will have a complete picture of the functions $\zeta (q,z)$ and $G-\zeta (z,q)$ which are closely related to Riemann's $\zeta (z)$.\\\\
Let's look in detail at the actual zeros per case in Table 16.
\newpage
\begin{table}[h!]
    \centering
    \includegraphics[width=13.5cm, height=6.5cm]{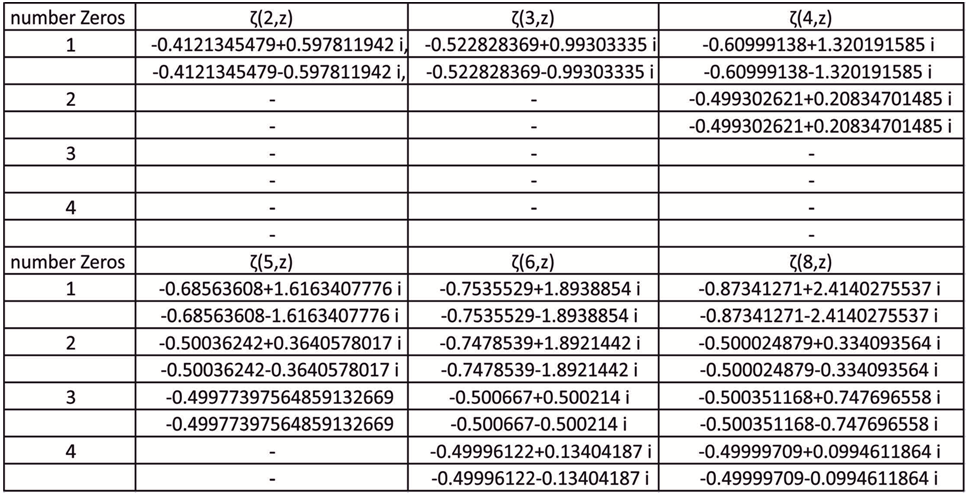}
    \caption{Zeros table for $q=\{2,3,4,5,6,8\}$} 
\end{table}
\phantom{ }\\
Suppose we have the zeros\\\\
\[\left\{z_{q_{i_1}^{1}}^{1}, z_{q_{i_2}^{1}}^{2}, z_{q_{i_3}^{1}}^{2}\right\},\left\{z_{q_{i_1}^{2}}^{1}, z_{q_{i_2}^{1}}^{2}, z_{q_{i_3}^{1}}^{2}\right\}, \ldots,\left\{z_{q_{i_{1}}^{k}}^{1}, z_{q_{i_2}^{k}}^{2}, z_{q_{i_3}^{k}}^{3}\right\}, \ldots,\left\{z_{q_{i_{1}}^{n}}^{1}, z_{q_{i_2}^{n}}^{2}, z_{q_{i_3}^{n}}^{3}\right\},\]\\
\[\text{ where } \left\{\kappa=1 \div n, i_{1}, i_{2}, i_{3} \in N\right\} \]\\\\
with ascending sequence of $q$, such that for the real parts of zeros apply the following\\\\
\[\boxed{\left\{\begin{array}{c}
\operatorname{Re}\left(z_{q_{i_1}^{k}}^{1}\right)=-\frac{1}{2}-\varepsilon_{q_{i_1}^{k}}^{1} \\
\operatorname{Re}\left(z_{q_{i_2}^{k}}^{2}\right)=-\frac{1}{2}+\varepsilon_{q_{i_2}^{k}}^{2} \\
\operatorname{Re}\left(z_{q_{i_3}^{k}}^{3}\right)=-1+\varepsilon_{q_{i_3}^{k}}^{3} \\
\kappa, n \in N, \kappa=1 \div n
\end{array}\right\} \quad \Leftrightarrow \quad\left\{\begin{array}{c}
\varepsilon_{q_{i_1}^{k}}^{1} \rightarrow 0 \\
\varepsilon_{q_{i_2}^{k}}^{2} \rightarrow 0 \\
\varepsilon_{q_{i_3}^{k}}^{2} \rightarrow 0 \\
\kappa>>2
\end{array}\right\} \wedge\left\{\begin{array}{l}
\operatorname{Re}\left(z_{q_{i_1}^{k}}^{1}\right) \rightarrow-\frac{1}{2} \\
\operatorname{Re}\left(z_{q_{i_2}^{k}}^{2}\right) \rightarrow-\frac{1}{2} \\
\operatorname{Re}\left(z_{q_{i_3}^{k}}^{3}\right) \rightarrow-1
\end{array}\right\}}\]\\\\
Therefore all zeros they have a negative real part of the complex zero, which shows that no zero is in the critical Strip of the interval $(0,1)$. So the \textbf{Generalized Riemann Zeta equation $\boldsymbol{\zeta(\mathrm{q}, \mathrm{z})=0, \mathrm{q}>1}$ in $\boldsymbol{\mathrm{Z}_{+}}$}, is outside of the critical line, but also has no critical line in the interval (0,1). For the Table 16, are used the functions Hurwitzzeta$[q, z]$ or $\operatorname{Zeta}(q, z)$ to check the zeros with the program mathematica 12, clearly with more
decimals to achieve approximations above $10^{-10}$.\\\\
The whole procedure for values of $q>1$, clearly shows that the function $\zeta (q,z)=0$ does not generate zeros even just in the interval $(0,1)$ which is the only critical line for the function $\zeta (z)=0$. Of course this does not exclude the possibility for other values of $q$, maybe it gives us some pairs or entire critical lines. This research will be completed below.
\newpage
\noindent\textbf{6.6.8. The zeros of equations $\boldsymbol{\zeta(-3, \mathbf{z})=0, \zeta(-2, \mathbf{z})=0, \zeta(-2 / 3, \mathbf{z})=0, \zeta(-1 / 2, \mathbf{z})=0,}$\\ $\boldsymbol{\zeta(-1 / 4, \mathbf{z})}$ $\boldsymbol{=0}$, and of equations $\boldsymbol{\zeta(1 / 4, \mathbf{z})=0, \zeta(1 / 2, \mathbf{z})=0, \zeta(2 / 3, z)=0}$ and the tables of these.}\\\\
To be able to close the important issue of the zeros of $\zeta(\mathrm{q}, \mathrm{z})=0$ and the positions of the values of $\operatorname{Re}(\mathrm{z})$ within the intervals $(0,1),(-1,0) \&(-\infty,-1)$ we need to solve several forms of transcendental equations and derive some characteristic matrices. The only method that has this possibility is the Periodic Radicals method, the same method that we solved the previous forms. the zeros we find will show the selective positions that are found and to which of these three intervals the $\operatorname{Re}(z)$ of each equation belongs.\\\\
\textbf{I. Control if $\boldsymbol{\operatorname{Re}(\mathbf{z})}$ of zeros in interval $\boldsymbol{(0,1)}$}\\\\
\textbf{I.1. Bernoulli polynomials}\\\\
Here there are 2 ways for integer values. Both the Infinity Periodic Radicals method and the Bernoulli polynomials. In this case we will prefer the Bernoulli polynomials. The Bernoulli polynomials are a generalization of the Bernoulli numbers. They have a variety of interesting properties, and will feature in our proof of the Euler-Maclaurin Summation Formula.
$$
\frac{x e^{x y}}{e^{x}-1}=\sum_{k=0}^{\infty} \frac{B_{k}(y) x^{k}}{k !}
$$
The generating function for the Bernoulli polynomials is the generating function for the Bernoulli numbers is multiplied by a term of $e^{x y}$. Our first observation of the Bernoulli polynomials is that the constant term of $\mathrm{B}_{k}(\mathrm{y})$ is in fact $\mathrm{B}$. If we set $\mathrm{y}=0$, then $\mathrm{B}_{k}(0)=\mathrm{B}_{k}$.
$$
\frac{x e^{x y}}{e^{x}-1}=\frac{x}{e^{x}-1}
$$
\textbf{Proposition.} The Bernoulli polynomials, $\mathrm{B}_{k}(\mathrm{y})$ satisfy the recurrence relation
$$
B_{k}(y)=\sum_{n=0}^{k}\left(\begin{array}{l}
k \\
n
\end{array}\right) B_{n} y^{k-n}
$$
$$
\hspace{8.5ex}\begin{aligned} \sum_{k=0}^{\infty} \frac{B_{k}(y) x^{k}}{k !} &=\frac{x e^{x y}}{e^{x}-1} \\ &=\frac{x}{e^{x}-1} \cdot e^{x y} \\ &=\sum_{k=0}^{\infty} \frac{B_{k} x^{k}}{k !} \cdot \sum_{k=0}^{\infty} \frac{(x y)^{k}}{k !} \end{aligned}$$\\
We find the Cauchy product by the same process as earlier to obtain
$$
\begin{aligned}
\sum_{k=0}^{\infty} \frac{B_{k}(y) x^{k}}{k !} &=\sum_{k=0}^{\infty} \sum_{n=0}^{k} \frac{(x y)^{k-n}}{(k-n) !} \cdot \frac{B_{n} x^{n}}{n !} \\
&=\sum_{k=0}^{\infty} \sum_{n=0}^{k} \frac{y^{k-n} B_{k}}{(k-n) ! n !} x^{k} \\
&=\sum_{k=0}^{\infty} \sum_{n=0}^{k}\left(\begin{array}{l}
k \\
n
\end{array}\right) y^{k-n} B_{n} \frac{x^{k}}{k !}
\end{aligned}
$$\\
Now, if we compare terms of $\mathrm{x}$ in the right and left-hand summations, we see the following:
$$
B_{k}(y)=\sum_{n=0}^{k}\left(\begin{array}{l}
k \\
n
\end{array}\right) B_{n} y^{k-n}
$$
as desired Using this recurrence, we can calculate the \_rst few Bernoulli polynomials:
$$
\begin{aligned}
&B_{0}(y)=1 \\
&B_{1}(y)=y-\frac{1}{2} \\
&B_{2}(y)=y^{2}-y+\frac{1}{6} \\
&B_{3}(y)=y^{3}-\frac{3}{2} y^{2}+\frac{1}{2} y \\
&B_{4}(y)=y^{4}-2 y^{3}+y^{2}-\frac{1}{30} \\
&B_{5}(y)=y^{5}-\frac{5}{2} y^{4}+\frac{5}{3} y^{3}-\frac{1}{6} y \\
&B_{6}(y)=y^{6}-3 y^{5}+\frac{5}{2} y^{4}-\frac{1}{2} y^{2}+\frac{1}{42}
\end{aligned}
$$\\
Notice our earlier result that the constant term of each polynomial is a Bernoulli number. Let us consider one of these polynomials, say $\text{B}_{5}(y)$, more closely. What if we differentiate it, or integrate it over 0 to $1?$
$$
\begin{aligned}
\frac{d}{d y} B_{5}(y) &=5 y^{4}-10 y^{3}+5 y^{2}-\frac{1}{6} \\
&=5\left(y^{4}-2 y^{3}+y^{2}-\frac{1}{30}\right) \\
&=5 B_{4}(y)
\end{aligned}
$$\\
so we observe that in this case,
$$
B_{k}^{\prime}(y)=k B_{k-1}(y)
$$\\
We also see
$$
\begin{aligned}
\int_{0}^{1} B_{5}(y) &=\frac{1}{6} y^{6}-\frac{1}{2} y^{5}+\frac{5}{12} y^{4}-\left.\frac{1}{12} y^{2}\right|_{0} ^{1} \\
&=0
\end{aligned}
$$\\
These two observations are, remarkably, true in general. In fact, they give us the following inductive definition for the Bernoulli numbers.\\\\
More generally about Bernoulli polynomials apply:\\
\begin{enumerate}
    \item $B_{0}(y)=1$
    \item $B_{k}^{\prime}(y)=k B_{k-1}(y)$;
    \item $\int_{0}^{1} B_{k}(y) d y=0$ for $k \geq 1$
\end{enumerate}
\newpage
\noindent\textbf{I.2. Zeros of $\boldsymbol{\zeta(-3, \mathbf{z})=0}$ and $\boldsymbol{\zeta(-2, \mathbf{z})=0}$.}\\\\\\
In special functions, the following applies to the Zeta function
$$
\boxed{\zeta(-q, z)=-\frac{B(q+1, z)}{q+1}}
$$
Therefore in these 2 cases we will have:\\\\
$$
\boxed{
\begin{aligned}
&\zeta(-2, z)=-\frac{B(3, z)}{4}=0 \Rightarrow \\
&\Rightarrow 1 / 3\left(-(z / 2)+\left(3 z^{\wedge} 2\right) / 2-z^{\wedge} 3\right)=0 \\
&\text { the zeros are : } z=0\|z=\frac{1}{2}\| z=1
\end{aligned}}
$$\\
And for the case $\zeta(-3, z)=0$ we have:\\\\
$$
\boxed{
\begin{aligned}
&\zeta(-3, z)=-\frac{B(4, z)}{4}=0 \Rightarrow\\
&\Rightarrow 1 / 4\left(1 / 30-z^{2}+2 z^{3}-z^{4}\right)=0 \\
&\text{the zeros are:} \\
&z=\frac{1}{2}\left(1-\sqrt{1-2 \sqrt{\frac{2}{15}}} \| z=\frac{1}{2}\left(1+\sqrt{1-2 \sqrt{\frac{2}{15}}} \|\right.\right.\\
&z=\frac{1}{2}\left(1-\sqrt{1+2 \sqrt{\frac{2}{15}}} \| z=\frac{1}{2}\left(1+\sqrt{1+2 \sqrt{\frac{2}{15}}}\right.\right. . \\
&\text{ie}\\
&z=0.240335\|z=0.759665\| \\
&z=-0.157704 \| z=1.1577
\end{aligned} }
$$\\\\
Solving other cases we come to the conclusion that we have 2 zeroes only real in the interval $(0,1)$. Therefore we have no critical lines in the interval $(0,1)$.\\\\\\\\
\textbf{I.3. Zeros of $\boldsymbol{\zeta(-2 / 3, z)=0}$, $\boldsymbol{\zeta(-1 / 2, z)=0}$, $\boldsymbol{\zeta(-1 / 4, z)=0}$}\\\\\\
This case can only be solved with Infinity Periodic Radicals with use Generalized Theorem[24]. The matrix of the program that generates all solutions is the following:
\newpage
The matrix program for $\zeta(-2/3,z)=0, \zeta(-1/2,z)=0, \zeta(-1/4,z)=0$  in mathematica with infinity Periodic Radicals is:\\
\phantom{ }
\begin{figure}[h!]
    \centering
    \includegraphics[width=13cm, height=8cm]{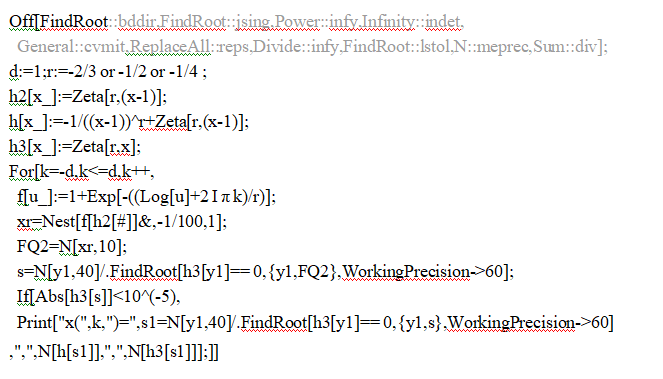}
\end{figure}
\phantom{ }\\
\vspace{\baselineskip}
\noindent 1. The zeros of $\zeta(-2 / 3, \mathbf{z})=0$ are:\\
$$
\boxed{
\begin{aligned}
&z_{1}=0.703226180818014 \\
&z_{2}=0.112695242215194
\end{aligned}}
$$\\
2. The zeros of $\zeta(-1 / 2, z)=0$ is:\\
$$
\begin{array}{|l|}
\hline z_{1}=0.066489124132138 \\
\hline
\end{array}
$$\\
3. The zeros of $\zeta(-1 / 4, \mathbf{z})=0$ is:\\
$$
\boxed{
\mathrm{z}_{1}=0.582245789118785}
$$\\
As we observe in all three cases we have only real roots although they are located in the interval $(0,1)$ they obviously neither create nor belong to any critical line.\\\\\\\\\\
\textbf{I.4. Zeros of $\boldsymbol{\zeta(1 / 4, z)=0}$, $\boldsymbol{\zeta(1 / 2, z)=0}$, $\boldsymbol{\zeta(2 / 3, z)=0}$}\\\\
Using the same program matrix we find all 3 zeros respectively. The matrix of the program that generates all solutions is the following:
\newpage
\begin{figure}[h!]
    \centering
    \includegraphics[width=13cm, height=8cm]{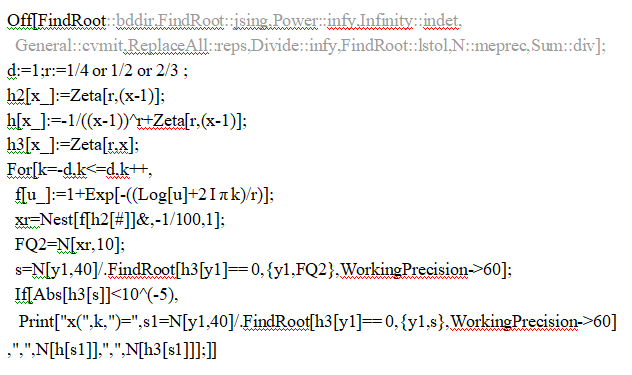}
\end{figure}
\vspace{\baselineskip}
\noindent 1. The zeros of $\zeta(1/4, \mathbf{z})=0$ are:\\
$$
\begin{array}{|l|}
\hline z_{1}=0.408064210108089 \\
\hline
\end{array}
$$\\
2. The zeros of $\zeta(1 / 2, z)=0$ is:\\
$$
\begin{array}{|l|}
\hline z_{1}=0.302721828598366 \\
\hline
\end{array}
$$\\
3. The zeros of $\zeta(2 / 3, \mathbf{z})=0$ is:\\
$$
\boxed{
\mathrm{z}_{1}=0.221547787091544}
$$\\
Therefore for all cases with $\mathrm{q}<0$ and $\mathrm{q}$ in $(0,1)$ while there are zeros in the set of real numbers in the interval $(0,1)$ there are no critical lines in this interval.\\\\\\\\

\noindent\textbf{6.6.9. Dirichlet L-Series and their relationship to the critical line.}\\\\
\textbf{I. Dirichlet L-series is a series} of the form [27]\\\\
\[\boxed{\mathrm{L}_{\mathrm{k}}(\mathrm{s}, \chi) \equiv \sum_{\mathrm{n}=0}^{\infty} \chi_{\kappa}(\mathrm{n}) \cdot \mathrm{n}^{-\mathrm{s}}}\]\\\\
where the number theoretic character $\chi_{k}(n)$ is an integer function with period $k$, are called Dirichlet L-series. These series are very important in additive number theory (they were used, for instance, to prove Dirichlet's theorem), and have a close connection with modular forms. Dirichlet L-series can be written as sums of Lerch transcendents with $z$ a power of $e^{2 \pi i / k}$. Dirichlet L-series is implemented in the Wolfram Language as DirichletL $[k, j, s]$ for the Dirichlet character $\chi(n)$ with modulus $k$ and index $j$. The generalized Riemann hypothesis conjectures that neither
the Riemann zeta function nor any Dirichlet L-series has a zero with real part larger than 1/2.\\\\ 
The Dirichlet lambda function also is,
\[\lambda(s) \equiv \sum_{n=0}^{\infty} \frac{1}{(2 n+1)^{s}}=\left(1-2^{-s}\right) \zeta(s)\]
The function Dirichlet $\left(\mathrm{k}, \mathrm{j}, \mathrm{s} \right)$, gives the Dirichlet L-function $\mathrm{L}(\chi, \mathrm{s})$ for the Dirichlet character $\chi(\mathrm{n})$ with modulus $k$ and index $j$.\\\\\\
\noindent
\textbf{Example(with mathematica):}\\\\
FindRoot [Dirichlet $[8,3, \mathrm{t}]==0,\{\mathrm{t}, 14I\}]$, Find $\mathrm{Root}[\mathrm{Dirichlet}[8,3, \mathrm{t}]==0,\{\mathrm{t}, 15 I\}]$\\
Out $:\{\mathrm{t} \rightarrow 0.5+12.9881 \mathrm{i}\},\{\mathrm{t} \rightarrow 0.5+16.3426 \mathrm{i}\}$\\\\
We therefore always \textbf{have as a critical line }$\mathbf{1} / \mathbf{2}$.\\\\\\
\textbf{II. The Dirichlet beta function} is defined by the sum\\
\[\boxed{\beta(s) \equiv \sum_{n=0}^{\infty} \frac{(-1)^{n}}{(2 n+1)^{s}}=2^{-s} \Phi\left(-1, s, \frac{1}{2}\right)}\]\\
where $\Phi(z, s, a)$ is the Lerch transcendent. The beta function can be written in terms of the Hurwitz zeta function $\zeta(x, a)$ by\\
\[\boxed{\beta(\mathrm{s}) \equiv \frac{1}{4^{\mathrm{s}}}\left(\zeta\left(\mathrm{x}, \frac{1}{4}\right)-\zeta\left(\mathrm{s}, \frac{3}{4}\right)\right)}\]\\
The beta function can be defined over the whole complex plane using analytic continuation,\\
\[\boxed{\beta(1-\mathrm{z})=\left(\frac{2}{\pi}\right)^{\mathrm{z}} \sin \left(\frac{1}{2} \cdot \pi \cdot \mathrm{z}\right) \cdot \Gamma(\mathrm{z}) \cdot \beta(\mathrm{z})}\]
where $\Gamma(z)$ is the gamma function.\\\\\\
\textbf{Example(with mathematica):}\\\\
Find Root[DirichletLBeta[s] $==0,\{\mathrm{t}, 14 I\}]$\\
Out : $\{t \rightarrow 0.5+12.9881 i\}$\\\\
We therefore always \textbf{have as a critical line} $\mathbf{1} / \mathbf{2}$.
\newpage
\noindent\textbf{III. Dirichlet Eta Function}\\\\
The Dirichlet eta function is the function $\eta(s)$ defined [28] by
\[\boxed{\eta(\mathrm{s}) \equiv \sum_{\mathrm{n}=1}^{\infty} \frac{(-1)^{\mathrm{n}-1}}{(\mathrm{n})^{\mathrm{s}}}=\left(1-2^{1-\mathrm{s}}\right) \cdot \zeta(\mathrm{s})}\]\\
where $\zeta(s)$ is the Riemann zeta function. Note that Borwein and Borwein (1987, p. 289) use the notation $\alpha(s)$ instead of $\eta(s)$. The function is also known as the alternating zeta function and denoted $\zeta^{*}(s)$, $\eta(0)=1 / 2$ is defined by setting $s=0$ in the right-hand side of $(2)$, while $\eta(1)=\ln 2$ (sometimes called the alternating harmonic series) is defined using the left-hand side. The function vanishes at each zero of $1-2^{1-s}$ except $s=1$. The eta function is related to the Riemann zeta function and Dirichlet lambda function by\\
\[\frac{\zeta(v)}{2^{v}}=\frac{\lambda(v)}{2^{v}-1}=\frac{\eta(v)}{2^{v}-2}\]\\
And $\zeta(v)+\eta(v)=2 \lambda(v)$. The eta function is also a special case of the polylogarithm function,\\
\[\eta(x)=-\operatorname{Li}_{x}(-1)\]
\vspace{0.5\baselineskip}
The derivative of the eta function is given by
\vspace{0.5\baselineskip}
\[\eta^{\prime}(x)=2^{1-x}(\ln 2) \zeta(x)+\left(1-2^{1-x}\right) \zeta^{\prime}(x)\]
\vspace{0.5\baselineskip}
Therefore, the natural logarithm of 2 is
\vspace{0.5\baselineskip}
\[\boxed{\eta(1) \equiv \sum_{n=1}^{\infty} \frac{(-1)^{n-1}}{(n)^{1}}=\ln (2)}\]\\\\
\textbf{Example(with mathematica):}\\\\
FindRoot [DirichletLEeta [s] = =0,\{t,14I\}],\\
Out : $\{\mathrm{t} \rightarrow 0.5+14.1347 \mathrm{i}\}$\\\\
We therefore always have as a critical line $1 / 2$.\\\\\\\\
\textbf{IV. Ramanujan z-Diriehlet series}\\\\
In this paper we consider the zeros of the Dirichlet series
\[L(s)=\sum_{n=1}^{\infty} \frac{\tau(n)}{n^{s}} \quad(\sigma=\operatorname{Re} s>13 / 2)\]
formed with Ramanujan's tau-function which may be defined by
\[\sum_{n=1}^{\infty} \tau(n) z^{n}=z \prod_{n=1}^{\infty}\left(1-z^{n}\right)^{24} \quad(|z|<1)\]
It is known that $L(s)$ has all of its non-real zeros in the strip $11/ 2<\sigma<13 / 2$ and it is conjectured that all these zeros are on $\sigma=6$ and are simple. It has been shown (by Hafner [29]) that a positive proportion of the zeros of $L(s)$ are of odd multiplicity and are on $\sigma=6$, but it is unknown whether any of the zeros are simple. Some other Dirichlet series are known to have infinitely many simple zeros. Heath-Brown [30] and Selberg independently observed that the work [31] of Levinson in 1974 implies that a positive proportion of the zeros of the Riemann zeta function are simple. Of course, the same result works for a Dirichlet $\mathrm{L}$-function.\\\\
\textbf{Example(with mathematica):}\\\\
FindRoot[RamanujanTauL[s]==0,\{t,5+9I\}],\\
Out : $\{\mathrm{t} \rightarrow 6+9.22238 \mathrm{i}\}$\\\\
We therefore always \textbf{have as a critical line in 6.}\\\\\\
\textbf{V. Zeros of the Davenport-Heilbronn counterexample}\\\\
We compute zeros off the critical line of a Dirichlet series considered by $\mathrm{H}$. Davenport and $\mathrm{H}$. Heilbronn. This computation is accomplished by deforming a Dirichlet series with a set of known zeros into the Davenport- Heilbronn series.\\\\
For $s=\sigma+i$ t with $\sigma>1$, let\\
\[\boxed{f_{1}(s)=1+\frac{\xi}{2^{s}}-\frac{\xi}{3^{s}}-\frac{1}{4^{s}}+\frac{0}{5^{s}}+\cdots}\tag{1}\]\\
be a Dirichlet series with periodic coefficients of period 5. Then $f_{1}(s)$ defines an entire function satisfying the following functional equation\\
\[\boxed{f(s)=T^{-s+\frac{1}{2}} \chi_{2}(s) f(1-s)}\tag{2}\]\\
with $T=5$ and\\
\[\boxed{\chi_{2}(s)=2(2 \pi)^{s-1} \Gamma(1-s) \cos \left(\frac{\pi s}{2}\right)}\tag{3}\]\\
In 1936, H. Davenport and H. Heilbronn (see [32]) showed that $f_{1}(s)$, as defined in (1), has zeros off the critical line $\sigma=1 / 2$. In 1994, R. Spira (see [33]) computed the following zeros of the Davenport-Heilbronn example:\\
\[\begin{aligned}
&\boldsymbol{0.808517+85.699348 i, 0.650830+114.163343 i }\\
&\boldsymbol{0.574356+166.479306 i, 0.724258+176.702461 i}
\end{aligned}\]
In this note we present a scheme for computing additional zeros of the Davenport- Heilbronn Dirichlet series. We certainly see isolated cases of pairs of complexes that individually define critical lines. We have seen this before in other cases with $\zeta$-functions and corresponding functional equations.\\\\\\
\textbf{"It is not a rare phenomenon but here we have to clarify which equations the Riemann hypothesis accepts as real cases as true cases that define the structure of the given hypothesis"}\\\\\\
\newpage
\section*{Part II. The paradox and the Rejection of the Hypothesis of\\ Riemann}
\noindent \textbf{The Riemann hypothesis incorrectly defines the non-trivial zeros of $\boldsymbol{\zeta(s)=0}$ as corresponding to the zeros of the sum of the series referred to by the function $\boldsymbol{\zeta()}$. This means that the non-trivial zeros of $\boldsymbol{z(s)=0}$ do not zero the sum represented by $\boldsymbol{\zeta()}$.}\\\\
\textbf{6.7 Introduction}\\\\
 What is in force to date. We know that the zeta function of Riemann is defined for complex s with real part greater than 1 by absolutely convergent series infinite order\\
\[\zeta(s)=\sum_{n=1}^{\infty} \frac{1}{n^s}=\frac{1}{1^s}+\frac{1}{2^s}+\cdots\tag{1}\]\\
In mathematics, the Riemann Hypothesis introduced by Riemann (1859) is the conjecture that the nontrivial roots of Riemann's zeta function all have real part $1 / 2$. In the strip $0<\operatorname{Re}(s)<1$ the zeta function satisfies the functional equation
\[\zeta(s)=\zeta(1-s) \cdot \Gamma(1-s) \cdot \sin \left(\frac{\pi s}{2}\right) \cdot 2^s \cdot \pi^{s-1}\tag{2}\]
One may then define $\zeta(s)$ for all remaining nonzero complex numbers $s(\operatorname{Re}(s) \leq 0$ and $s \neq 0)$ by applying this equation outside the strip, and letting $\zeta(s)$ equal the right-hand side of the equation whenever $s$ has non-positive real part (and $s \neq 0$ ). If $s$ is a negative even integer then $\zeta(s)=0$ because the factor $\sin (\pi s / 2)$ vanishes; these are the trivial zeros of the zeta function. (If $s$ is a positive even integer this argument does not apply because the zeros of the sine function are cancelled by the poles of the gamma function as it takes negative integer arguments.) The value $\zeta(0)=-1 / 2$ is not determined by the functional equation, but is the limiting value of $\zeta$ (s) as $s$ approaches zero. The functional equation also implies that the zeta function has no zeros with negative real part other than the trivial zeros, so all non-trivial zeros lie in the critical strip where $s$ has real part between 0 and 1.\\\\
\textbf{6.7.I. How to test the hypothesis and where is the contradiction and logical error in finding the non-trivial roots of the $z$ function and the sum.}\\\\
The more general form of the sum that includes the Riemann hypothesis has the form\\
\[{ }_{ \text{ }\text{ } m }^{\alpha, \beta} \sigma( s )=\sum_{n=1}^{\infty} \frac{1}{(\alpha \cdot n+\beta)^{m \cdot s}}=\frac{1}{(\alpha+\beta)^{m \cdot s}}+\frac{1}{(\alpha \cdot 2+\beta)^{m \cdot s}}+\cdots, n \in N^{+},\{\alpha, \beta, m \in R\}, s \in C \tag{3}\]\\
According to the above definitions, we will prove that the Riemann hypothesis, as defined, falls into contradictions, \textbf{firstly} because the non-trivial roots of the sum refer to the zeta function and are not the correct roots of the sum (1) when it is zeroed. They refer only to the Zeta function and when it is zeroed and to another sum. The known nontrivial roots of $\zeta( s )=0$ and $m=1$ refer to another sum for which the function $\zeta(s)$ is defined to be equal to a part of that sum, but which, if zeroed, then its roots are equivalent to the nontrivial roots of $\zeta(s)=0$. \textbf{Second}, for the generalized sum we are interested in (3) to have non-trivial roots, it must be true that $\alpha<0, \beta>0$ or $\alpha<0$ and $\beta<0$ and also in the generalized form of the sum when $m=1$, when we refer to the sum that the Riemann hypothesis deals with. When we go up integers values for $m$ then we can find non-trivial roots for $\alpha>0$. $R(s)=1 / 2$ is an adapted solution of an approximate formulation of the hypothesis, and concerns roots of $\zeta(s)=0$, as the $\zeta()$ function is defined, but which corresponds to a special case of the sum(3) not (1). This can be seen very simply if we substitute the complex values of the already found non-trivial solutions into the variable s of the sum(1), then we see that the sum does not go to zero, but instead for very high values it is almost infinite, i.e. gradually as the series grows it moves away from zero completely. This is demonstrated in other ways that we will see at the end of the chapter.
\newpage
\noindent\textbf{6.7.II. Generic forms of the Hurwitz zeta function}\\\\
The Hurwitz zeta function, is a generalization of the Riemann zeta function, that is also known as the generalized zeta function and is one of the fundamental transcendental functions and traditionally and defined [34] by the series\\
\[\zeta(s, \alpha)=\sum_{\kappa=0}^{\infty} \frac{1}{(\kappa+\alpha)^s}, R(s)>1, \alpha \notin Z_0^{-} \tag{4}\]\\
The series converges absolutely for $R(s)=\sigma>1$ and the convergence is uniform in every half-plane $\sigma>=1+\delta(\delta>0)$. Therefore, $\zeta(s, a)$ is analytic function of $s$ in the half-plane $R(s)=\sigma>1$. It satisfies the following integral representation\\
\[\zeta(s, \alpha)=\frac{1}{\Gamma(s)} \int_0^{\infty} \frac{x^{s-1} e^{-\alpha x}}{1-e^{-x}} d x, R(s)>1, R(\alpha)>0 \tag{5}\]\\
Using the Taylor expansion\\
\[\frac{1}{\left(1-e^x\right)}=-\frac{1}{x}+\frac{1}{2}+O(x)\text{ and splitting (2) into integrals}\]\\
\[\zeta(s, \alpha+1)=\frac{1}{\Gamma(s)} \int_0^{\infty} x^{s-1} e^{-\alpha x}\left(\frac{1}{1-e^x}-\frac{1}{x}+\frac{1}{2}\right) d x+\frac{1}{\Gamma(s)} \int_0^{\infty} x^{s-1} e^{-\alpha x}\left(\frac{1}{x}-\frac{1}{2}\right) d x\tag{6}\]\\
we can analytically continue $\zeta(s, a)$ into the strip $-1< R (s)<1$.\\\\
\noindent Evaluating the  integral, replacing $a$ by $a+1$ and making use of (Derivable from (4) by comparing its $n=0$ term with its sum over $n>=1$ )\\
\[\zeta(s, \alpha+1)=\zeta(s, \alpha)-\frac{1}{\alpha^s}, R(\alpha)>0\tag{7}\]\\
We obtain\\
\[\zeta(s, \alpha+1)=\frac{\alpha^{1-s}}{s-1}+\frac{\alpha^{-s}}{2}+\frac{1}{\Gamma(s)} \int_0^{\infty} x^{s-1} e^{-\alpha x}\left(\frac{1}{1-e^x}-\frac{1}{x}+\frac{1}{2}\right) d x, R(s)>-1, R (\alpha)>0\tag{8}\]\\
This could be further generalized\\
\[\begin{aligned}
\zeta( s , \alpha+1)&=\alpha^{- s }+\sum_{ k =0}^{ n } \frac{\Gamma( k + s -1)}{\Gamma( s )} \frac{B_{ k }}{ k } \alpha^{1-s- k }+ \\[5pt]
&+\frac{1}{\Gamma( s )} \int_0^{\infty} x ^{ s -1} e ^{-\alpha x}\left(\frac{1}{1- e ^{ x }}-\sum_{ k =0}^n \frac{B_{ k }}{ k } \alpha^{1- k }\right) dx,
\end{aligned}\]\\
\[R ( s )>-2[ n / 2]-1, R (\alpha)>0 \tag{9}\]\\
Where $B$ are the Bernoulli numbers. We can also analytically continue the Harwize zeta function to the whole complex s-plane (exept for a simple pole at $s=1$ )by means of the contour integral\\
\[\zeta(s, \alpha)=\frac{\Gamma(1-s)}{2 \pi i} \int_c \frac{x^{s-1} e^{\alpha x}}{1-e^{-x}} d x \tag{10}\]\\
Where the contour $C$ is a loop that starts from $-\infty$ along the lower side of the real axis, encircle the origin and then returns to $-\infty$ along the upper side of the real axis. From the integral (10) one come derive the Hurwitz series representation\\
\[\zeta( s , \alpha)=\frac{\Gamma(1- s ) i }{(2 \pi)^{1- s }}\left( e ^{-\pi i i / 2} \sum_{ n =1}^{\infty} \frac{ e ^{-2 \pi in \alpha}}{ n ^{1- s }}- e ^{\pi is / 2} \sum_{ n =1}^{\infty} \frac{ e ^{2 \pi i n \alpha}}{ n ^{1- s }}\right)\tag{11}\]\\
Valid in the half -palane $R ( s )<0$ and $0<\alpha<1$.If $a \neq 1$ this representation is also valid for $R((s)<1$. The Dirichlet series in (11) can be rewritten in terms of the polylogarithm $\operatorname{Li} m(z)$:\\
\[\zeta(s, \alpha)=\frac{\Gamma(1-s) i}{(2 \pi)^{1-s}}\left( e ^{-\pi is / 2} Li _{1-s}\left( e ^{2 \pi i \alpha}\right)- e ^{\pi is / 2} Li _{1-s}\left( e ^{-2 \pi i \alpha}\right)\right)\tag{12}\]\\
Or, equivalently, (for $R(s)>1)$\\
\[\zeta(1- s , \alpha)=\frac{\Gamma(1- s ) i }{(2 \pi)^{1-s}}\left( e ^{-\pi i s / 2} Li _{ s }\left( e ^{-2 \pi i \alpha}\right)+ e ^{\pi i i / 2} Li _{ s }\left( e ^{2 \pi i \alpha}\right)\right)\tag{13}\]\\
Setting $\alpha=1$, yields the functional equation for the Riemann zeta function\\
\[\zeta(1- s )=2 \frac{\Gamma( s )}{(2 \pi)^s} \cos \left(\frac{\pi s }{2}\right) \zeta( s )\tag{14}\]\\
The result holds for all admissible $s$ by analytic continuation.\\\\
\textbf{The function $\boldsymbol{\zeta()}$ in relation to the generalised sum} $\zeta(s)=\sum_{n=1}^{\infty} \frac{1}{n^s}=\frac{1}{1^s}+\frac{1}{2^s}+\cdots$.\\\\
The function $\zeta()$ is an approximate function that theoretically to date represents the sum $\sum_{n=1}^{\infty} \frac{1}{n^s}$
constructed by a particular method, but does not fully correspond to the sum. If we find the zeros of $\boldsymbol\zeta(s)=0$, i.e. the nontrivial complex roots, we observe that they do not zero the sum
Therefore, the strongest criterion for finding non-trivial roots $s$ of the sum is that the equation $\zeta(s)=0$ is zero and at the same time the sum itself is zero. The usefulness of $\zeta()$ is clearly great, and this follows both from its approximate relation to the sum and from the possibility provided by the very important relation (7) to find an infinite number of complex roots for many sums of generalized form (3) when they are equal zero. Finding non-trivial zeros for the \textbf{generalized series (3)} is achieved by solving transcendental equations directly related to the \textbf{Hurwitz zeta function}.\\\\\\
\textbf{6.7.III. Correlation of the Hurwitz Zeta function with the generalized sum $\boldsymbol{{ }_{ \text{ }\text{ } m }^{\alpha, \beta} \sigma( s )}$}\\\\
Our ultimate goal is to find a relation connecting the function $\zeta()$ and the generalized sum. The approximate correlation of $\zeta ()$ is a strong link that directly helps to determine the non-trivial roots of the generalized sum.\\
\[{ }_{ \text{ }\text{ } m }^{\alpha, \beta} \sigma( s )=\sum_{n=1}^{\infty} \frac{1}{(\alpha \cdot n+\beta)^{m \cdot s}} \tag{III.1}\]\\
is exactly what we were looking for.\\\\
Classically the Riemann zeta function or more generally, The Hurwitz zeta function, is defined on a half plane using a series and then it is analytically extended, with respect to s, to the whole plane except for a simple pole at $s=1$ with residue 1,\\
\[{ }_{ \text{ }\text{ }\mathsmaller{G}} \zeta(s, \beta)=\sum_{n=0}^{\infty} \frac{1}{(n+\beta)^s}, \text{ for } R(s)>1 \text{ and } 0<\beta \leq 1, m=1\tag{III.2}\]\\
If we have the sum of case (III.1) we have the relation in combination with (III.2)\\
\[{ }_{ \text{ }\text{ } \mathsmaller{G} }\zeta(s,(\beta+1))=\sum_{n=1}^{\infty} \frac{1}{(n+\beta)^s}, \text{ for } R(s)>1 \text{ and } 0<\beta+1 \leq 1 \text{ and } m=1, \beta \leq 0 \tag{III.3}\]\\
In the case of the generalized sum case (III.1), we obtain for the sum the relation\\
\[{ }_{ \text{ }\text{ } \mathsmaller{G} }\zeta(s,(\beta+\alpha) / \alpha)=\alpha^s \sum_{n=1}^{\infty} \frac{1}{(\alpha \cdot n+\beta)^s}=\alpha^{-s+s} \sum_{n=1}^{\infty} \frac{1}{(n+\beta / \alpha)^s}=\sum_{n=1}^{\infty} \frac{1}{(n+\beta / \alpha)^s}, \text{ for } R(s)>1\]\\
\[\text{ and } 0<\beta / \alpha+1 \leq 1 \text{ and } \beta / \alpha \leq 0, m=1\tag{III.4}\]\\
Also apply if $p=\dfrac{\alpha+\beta}{\alpha}$\\
\[\sum_{n=1}^{\infty} \frac{1}{(\alpha \cdot n+\beta)^s}=\alpha^{-s} \cdot \zeta(s,(\beta+\alpha) / \alpha), \text{ for } R(s)>1 \text{ and } 0<\beta / \alpha+1 \leq 1 \text{ and } \beta / \alpha \leq 0\tag{III.5}\]\\
From relation (7) of section II, we obtain\\
\[\zeta(s, \alpha+1)=\zeta(s, \alpha)-\frac{1}{\alpha^s} \Leftrightarrow \zeta(s, \alpha)=\zeta(s, \alpha+1)+\frac{1}{\alpha^s}, R(\alpha)>0\tag{III.6}\]\\
\[\sum_{n=1}^{\infty} \frac{1}{(\alpha \cdot n+\beta)^s}=\alpha^{-s} \cdot \zeta(s,(\beta+\alpha) / \alpha)=\alpha^{-s}\left(\frac{1}{p^s}+\zeta(s, 1+(\beta+\alpha) / \alpha)\right), \text { for } R(s)>1\]\\
\[ \text{and } 0<\beta / \alpha+1 \leq 1 \text{ and } \beta / \alpha \leq 0\tag{III.7}\]\\
If we now have the more generalized relation which is the numbers raised to power $m ^* s$ then we will get the generalized sum in terms of the \textbf{Hurwitz Zeta function} in the form\\
\[{ }_{ \text{ }\text{ } m }^{\alpha, \beta} \sigma( s )=\sum_{n=1}^{\infty} \frac{1}{(\alpha \cdot n+\beta)^{m \cdot s}}=\alpha^{-m \cdot s} \cdot \zeta(m \cdot s,(\beta+\alpha) / \alpha)=\alpha^{-m \cdot s}\left( \frac{1}{p^{m \cdot s}}+\zeta(m \cdot s,p)\right)=0, \]\\
\[\text { for } R(s)>1, \text{ where } 0<\beta / \alpha+1 \leq 1 \text{ and } \beta / \alpha \leq 0,\text{ } m \in R \tag{III.8}\]\\\\
Relation (III.8) is very useful when we generalize below for the Multiple Zeta function in relation to that ${ }_{ \text{ }\text{ } m }^{\alpha, \beta} \sigma( s )$, for which we know relation (3, chapter I)\\
\[{ }_{ \text{ }\text{ } m }^{\alpha, \beta} \sigma( s )=\sum_{n=1}^{\infty} \frac{1}{(\alpha \cdot n+\beta)^{m \cdot s}}=\frac{1}{(\alpha+\beta)^{m \cdot s}}+\frac{1}{(\alpha \cdot 2+\beta)^{m \cdot s}}+\cdots, n \in N^{+},\{\alpha, \beta, m \in R\}, s \in C.\]\\
The question of why it is useful is answered directly by the (GRIM) method[ ] which is a powerful procedural method for solving transcendental and polynomial equations using an iterative method.. We will see how this is done in practice in the next chapter IV.Finding the non-trivial complex roots resulting from the solution of the sum equation ${ }_{ \text{ }\text{ } m }^{\alpha, \beta} \sigma( s )=0$.
\newpage\noindent
In principle, we start from relation (III.8) and if we assume $p=1+\frac{\beta}{\alpha}$\\
\[{ }_{ \text{ }\text{ } m }^{\alpha, \beta} \sigma( s )=\sum_{n=1}^{\infty} \frac{1}{(\alpha \cdot n+\beta)^{m \cdot s}}=\alpha^{-m \cdot s} \cdot \zeta(m \cdot s,(\beta+\alpha) / \alpha)=\alpha^{-m \cdot s}\left(\frac{1}{p^{m \cdot s}}+\zeta(m \cdot s, 1+p)\right)=0\]\\
for $m=1, R(s)>1$ where $0<\beta / \alpha+1 \leq 1$ and $\beta / \alpha \leq 0 \quad m \in C$ \\\\
Corresponding to the generalised sum ${ }_{ \text{ }\text{ } m }^{\alpha, \beta} \sigma( s )=0$ for which we ask for the non-trivial roots when it is zeroed. As we observe, the relation\\
\[\begin{aligned}
&{ }_{ \text{ }\text{ } m }^{\alpha, \beta} \sigma( s )=\sum_{n=1}^{\infty} \frac{1}{(\alpha \cdot n+\beta)^{m \cdot s}}=\alpha^{-m \cdot s} \cdot \zeta(m \cdot s,(\beta+\alpha) / \alpha)=\alpha^{-m \cdot s}\left(\frac{1}{p^{m \cdot s}}+\zeta(m \cdot s, 1+p)\right)=0 \Leftrightarrow \\
&{ }_{ \text{ }\text{ }\mathsmaller{G}} {\zeta(m \cdot s,(\beta+\alpha) / \alpha)}=a^{-m\cdot s}\zeta(m \cdot s,(\beta+\alpha) / \alpha)=a^{-m\cdot s}\left(\frac{1}{p^{m \cdot s}}+\zeta(m \cdot s, 1+p)\right)=0 \Leftrightarrow \\
&\frac{1}{p^{m \cdot s}}+\zeta(m \cdot s, 1+p)=0
\end{aligned}\tag{IV.1}\]\\
If $i$ call $p =\frac{\alpha+\beta}{\alpha}, m =1$ therefore will have that:\\\\
${ }_{ \text{ }\text{ } \mathsmaller{G} }{\zeta( s , p )}=
a^{-s} \left( \frac{1}{p^s}+\zeta(s,1+p) \right)=0$ where ${ }_{ \text{ }\text{ } \mathsmaller{G} }\zeta( s , p )$ the generalized equation that zeroes the sum ${ }_{ \text{ }\text{ } m }^{\alpha, \beta} \sigma( s )$.
If $i$ replace with $\left(\frac{1}{ p }\right)^{ s }= u$
$\{ k \in Z, u \neq 0, p \neq 0\}$ then\\
\[s = f ( u )=\frac{2 \cdot i \cdot \pi \cdot k }{\log (1 / p )}+\frac{\log (u)}{\log (1 / p )} \tag{IV.2}\]\\
the resulting function $f ( u )$ and value of s, named "Generator function".\\\\
H Generator function finds, according to the GRIM[37] method with a fixed program, the non-trivial roots of case with different $\{ a , b , m \}$. We will see in another chapter, where we will analyze the cases when a,b take different values and different signs how many categories we have in total. Of all the categories we will mention about $\{ a , b \} 4$ in total only 2 are of interest to us. The calculation and the approximation is of course done with mathematica 12 with a high-powered programming language.\\\\\\
\textbf{6.7.IV.a. Program for finding non-trivial roots in mathematica language, method (GRIM[37]) :}\\\\
Off[FindRoot::bddir,FindRoot::jsing,FindRoot::precw,Power::infy,Infinity::indet,General::cvmit,General::st op,\\[5pt]
ReplaceAll::reps,Divide::infy,FindRoot::Istol,N::meprec];\\\\
\textbf{'program data'}\\\\
$d =$ Input["Range values of $k , d =$ "];\\[5pt]
$a=$ Input["Value of $a=$ " ]\\[5pt]
$b =$ Input["Value of $b =$ "];\\[5pt]
InAp=Input["Initial Approx of Iteration, $\operatorname{In} A p=$ " ];\\[5pt]
$R =$ Input["Replays of Iterations, $R =$ "];\\[5pt]
Asum=Input["Approach value of Sum, Asum="];\\[5pt]
$m =$ Input["Multiplicity of exponent of $m =$ "];\\\\
\textbf{'program execution'}\\\\
$p :=( a + b ) / a$;\\[5pt]
h2 $\left[ x _{-}\right]:=$ -Zeta $[ mx ,(1+ p )]$;\\[5pt]
$h \left[ x _{-}\right]:=\left((1 /( p ))^{\wedge}( mx )+\operatorname{Zeta}[ mx ,(1+ p )]\right) ;$\\[5pt]
$h 3\left[ x _{-}\right]:=\operatorname{Zeta}[ mx , p ] ; j :=0$;\\[5pt]
$h 4\left[ x _{-}\right]:=\sum_{ k =1}^{\infty} \frac{1}{( a \cdot k + b )^{\wedge}( mx )}$;\\[5pt]
For $[k=-d, k<=d, k++$,\\[5pt]
$f\left[u_{-}\right]:=(2$ I $\pi k+\log [u]) / \log [1 /(p)] ;$\\[5pt]
xr=Nest[f[h2[\#]] \&,InAp,R];\\[5pt]
$FQ 2= N [ xr , 10]$;\\[5pt]
$s = N [ y 1,40] /$. FindRoot $[ h 4[ y 1]==0,\{ y 1, FQ 2\}$, WorkingPrecision $\to 580]$;\\[5pt]
$s 2= N [ y 3,40] /$. FindRoot $[ h 4[ y 3]=0,\{ y 3, s \}$, WorkingPrecision- $\to 180]$;\\[5pt]
If $[0<\operatorname{Re}[ s ]<1$ \& \& Abs[h4[s2]] $<$ Asum \& \& Abs $[ h3 [ s 2]< A s u m , j = j +1$;\\[5pt]
$\operatorname{Print}[" x (", k , ")=$ ", $s 1= N [ y 2,40] /$. FindRoot $[ h 4[ y 2]==0,\{ y 2, s \}$,\\[5pt]
WorkingPrecision$\to$180],",",N[h[s1]],",",N[h4[s1]],",",N[h3[s1]]];,Loopback]]\\[5pt]
Print["Percent of Zeros\%=", $\left.N \left[ j /(2 d +1)^* 100\right]\right]$\\\\
The results of the program are very useful in terms of which roots we have for the sum data a,b and also in terms for checking the sum with respect to the Hurwitz Zeta function. We will have to look at 3 Categories and see in practice what happens depending on $\alpha$ and $\beta$. Of course the research will extend to generalized summation ${ }_{ \text{ }\text{ } m }^{\alpha, \beta} \sigma( s )=0$ and see what additional variations there are when we have multiplicity in the exponent $x$. We usually choose values for the approaches Asum $=10^{-10}$ and $\ln A p=\pm \frac{1}{1000}$.\\\\\\
\textbf{6.7.IV.b. Tow very important forms of $\boldsymbol{\zeta( s , a )}$}\\\\
Searching for better forms of equations for $\zeta( s , a )$ from a wide range of equations we choose 2 forms that best fit the GRIM solution method.The group of Kanemitsu, Tanigawa, Tsukada and Yoshimoto obtained a new proof of $\zeta( s , a )$. Their proof begins with the use of [10,Equation (see also [35, Equation (47)])\\
\[\begin{aligned}
& \zeta(s, a)=1 / 2 * p^{-s}-p^{1-s} /(s-1)+ \\
& +\sum_{n=1}^{\infty}\left(\frac{e^{-2 \pi \cdot i * n * p}}{(-2 \pi * i * n * p)^{\wedge}(1-s)} * \Gamma[1-s,-2 \pi * i * n * p]+\frac{e^{2 \pi * i * n * p}}{(2 \pi * i * n * p)^{\wedge}(1-s)} * \Gamma[1-s, 2 \pi * i * n * p]\right) 
\end{aligned}\tag{IV.3}\]\\
which is a special case of the Ueno-Nishizawa formula [24] and then invoking the Fourier series of the Dirac-delta function (s).The aim of this note is to give a yet another new proof of beginning with Hermite's well-known formula for $\zeta( s$, a) [36, p. 609, Formula 25.11.29], valid for $\operatorname{Re}( a )>0$\\
\[\zeta(s, a)=\left(1 / 2 *(1 /(p))^s-p^{1-s} /(1-s)\right)+2 * \int_0^{\infty} \frac{\operatorname{Sin}[s * \operatorname{ArcTan}[x / p]]}{\left(p^{\wedge} 2+x^{\wedge} 2\right)^{\wedge}(s / 2) *(\operatorname{Exp}[2 * \pi * x]-1)} d x \tag{IV.4}\]\\
With these forms, especially IV.4, we can form a program that will be more complete for calculating the non-trivial zeros of each sum.\\\\\\
\textbf{6.7.IV.b.1 Two programs for finding non-trivial roots in mathematica language, method (GRIM) :}\\\\
\textbf{6.7.IV.b.1.I. The first form of the kernel with inversion of $\boldsymbol{\frac{1}{2} \cdot p^x}$}\\\\
Off[FindRoot::bddir,NIntegrate::ncvb,NIntegrate::slwcon,FindRoot::jsing,Power:infy, Infinity:indet,\\[5pt]
General::cumit,ReplaceAll:reps,Divide::infy,FindRoot:1stol,N:meprec, General:munfl];\\\\
\textbf{'program data'}\\\\
$d =$ Input["Range values of $k , d =$ "];\\[5pt]
$a=\operatorname{Input}[$ "Value of $a=$ "];\\[5pt]
$b=$ Input["Value of $b=$ ""];\\[5pt]
InAp=Input["Initial Approx of Iteration, $\operatorname{InAp}=$ "];\\[5pt]
$R =$ Input["Replays of Iterations, $R =$ "];\\[5pt]
Asum=Input["Approach value of Sum, Asum="];\\\\
\textbf{'program execution'}\\\\
$\begin{aligned}
& p :=( a + b ) / a ; \\[5pt]
& h 2\left[ x _{-}\right]:= p ^{\wedge}(1- x ) /(1- x )- \\[5pt]
& \left.-2^* \text { NIntegrate }\left[\operatorname{Sin}\left[ x ^* \operatorname{ArcTan}[ w / p ]\right] /\left( p ^{\wedge} 2+ w ^{\wedge}\right)^{\wedge}( x / 2)^*\left(\operatorname{Exp}\left[2^*[ P 1]^* w \right]-1\right)\right),\{ w , 0,10000\}\right] ; \\[5pt]
& h \left[ x _{-}\right]:=\left(1 / 2^*(1 /( p ))^{\wedge} x - p ^{\wedge}(1- x ) /(1- x )\right)+ \\[5pt]
& +2^* \text { NIntegrate }\left[\operatorname{Sin}\left[x^* \operatorname{ArcTan}[ w / p ]\right] /\left(\left( p ^{\wedge} 2+ w ^{\wedge} 2\right)^{\wedge}( x / 2)^*\left(\operatorname{Exp}\left[2^*[ Pi ]^* w \right]-1\right)\right),\{ w , 0,10000\}\right] ; \\[5pt]
& h 3\left[ x _{-}\right]:=\operatorname{Zeta}[ x , p ] ; j :=0 \text {; } \\[5pt]
\end{aligned}$\\
$h 4\left[ x _{-}\right]:=\sum_{k=1}^{\infty} \frac{1}{(a \cdot k+b)^{\wedge}(x)}$\\[5pt]
For $[ k =- d , k < d , k ++$,\\[5pt]
$f \left[ u _{-}\right]:=(2 I \pi k +\log [2 u ]) / \log [ 1 /( p )]$;\\[5pt]
xr=Nest [f[h2[\#]]\&,InAp,R];\\[5pt]
$FQ 2= N [ xr , 10]$;\\[5pt]
$s = N [ y 1,40] /$ FindRoot $[ h 4[ y 1]=0,\{ y 1, FQ 2\}$, WorkingPrecision $\to$ 580];\\[5pt]
$s 2= N [ y 3,40] /$ FindRoot $[ h 4[ y 3]=0,\{ y 3, s \}$, WorkingPrecision $\to 180]$;\\[5pt]
If $[0<\operatorname{Re}[ s ]<1 \& \& \quad A b s[h 4[s 2]]<$ Asum \& \& $A b s[h 3[$ s 2$]<$ 
Asum, $j=j+1 ;$\\[5pt]
Print$[" x (", k , ")=\text {, }, s 1= N [ y 2,40] / \text { FindRoot }[ h 4[ y 2]=0,\{ y 2, s \}$,\\[5pt]
WorkingPrecision$\to$180],",",N[h[s1]],",",N[h4[s1]],",",N[h3[s1]]];,Loopback]] \\[5pt]
Print["Percent of Zeros $\%=$ ", N[ $\left.\left./(2 d +1)^* 100\right]\right]$
\newpage\noindent
\textbf{6.7.IV.b.1.II. The second form of the kernel with inversion of $\boldsymbol{\frac{p^{1-x}}{1-x}}$}\\\\
Off[FindRoot::bddir,NIntegrate:uncvb, NIntegrate::slwcon, FindRoot: jsing,Power:infy,Infinity:uindet,\\
General::cvmit,ReplaceAll: :reps,Divide::infy, FindRoot: $: 1$ stol,N: meprec, General:munfl];\\\\\\
\textbf{'program data'}\\\\
$d =$ Input["Range values of $k , d =$ "];\\[5pt]
$a=\operatorname{Input}[$ "Value of $a=$ "];\\[5pt]
$b =$ Input["Value of $b =$ "];\\[5pt]
In Ap=Input["Initial Approx of Iteration, InAp="];\\[5pt]
$R =$ Input["Replays of Iterations, $R =$ "];\\[5pt]
Asum=Input["Approach value of Sum, Asum="];\\\\\\
\textbf{'program execution'}\\\\
$p :=( a + b ) / a ;$\\[5pt]
h2$\left[ x _{-}\right]$ := $=1 / 2 * p ^{\wedge}(- x )+$\\[5pt]
$+2^* NIn$ tegrate $\left[\operatorname{Sin}\left[ x ^* \operatorname{ArcTan}[ w / p ]\right] /\left(\left( p ^{\wedge} 2+ w ^{\wedge} 2\right)^{\wedge}( x / 2)^*\left(\operatorname{Exp}\left[2 * \backslash[ Pi ]^* w \right]-1\right)\right),\{ w , 0,10000\}\right] ;$\\[5pt]
$h \left[ x _{-}\right]:=\left(1 / 2 *(1 /( p ))^{\wedge} x - p ^{\wedge}(1- x ) /(1- x )\right)+$\\[5pt]
$+2^*$ NIntegrate $\left[\operatorname{Sin}\left[ x ^* \operatorname{ArcTan}[ w / p ]\right] /\left(\left( p ^{\wedge} 2+ w ^{\wedge} 2\right)^{\wedge}( x / 2)^*\left(\operatorname{Exp}\left[2^* \backslash[ Pi ]^* w \right]-1\right)\right),\{ w , 0,10000\}\right] ;$\\[5pt]
$h 3\left[ x _{-}\right]:=\operatorname{Zeta}[ x , p ] ; j :=0$;\\[5pt]
$h 4\left[ x _{-}\right]:=\sum_{k=1}^{\infty} \frac{1}{( a \cdot k + b )^{\wedge}( x )}$\\[5pt]
For[k=-d,k<=d,k++\\[5pt]
$f \left[ x _{-}\right]:= 1 +($ ProductLog[k,-Log[p]/u] $) / L o g [ p ]$\\[5pt]
xr=Nest[f[h2[\#]]\&,InAp,R]\\[5pt]
$FQ 2= N [ xr , 10] ;$\\[5pt]
$s = N [ y 1,40] /$ FindRoot $[ h 4[ y 1]=0,\{ y 1, FQ 2\}$, WorkingPrecision $\to 580]$;\\[5pt]
$s 2= N [ y 3,40] /$ FindRoot $[ h 4[ y 3]=0,\{ y 3, s \}$, WorkingPrecision $\to 180]$;\\[5pt]
If $[0<\operatorname{Re}[ s ]<1 \& \&$ Abs $[ h 4[ s 2]]<$ Asum \& \& Abs[h3[s2]<Asum, $j = j +1 ;$\\[5pt]
$\operatorname{Print}[" x (", k , ")=", s 1= N [ y 2,40] /$ FindRoot $[ h 4[ y 2]=0,\{ y 2, s \}$,\\[5pt]
WorkingPrecision$\to$ 180],",",N[h[s1]],",",N[h4[s1]],",",N[h3[s1]]];,Loopback]]\\[5pt]
Print["Percent of Zeros $\left.\%=, N \left[ j /(2 d +1)^* 100\right]\right]$\\\\
Clearly it is a more sophisticated and more complete solution of non-trivial zeros of sums. The method we use is GRIM, and it turns out to be very simple and useful here. Personally, I think it is the only method of finding and calculating the non-trivial zeros of each sum mentioned.\\

\newpage
\noindent\textbf{6.7.V. Solve the generalized equation $\boldsymbol{{ }_G \zeta( m \cdot s , p )=0}$ of Hurwitz Zeta function.}\\\\
To get a general picture of all forms of functions $\zeta()$ we have to deal with the generalized function which we will call ${}_{\text{ }\text{ }\mathsmaller{G}}\zeta( s , p )=0$. If we define a complex $S =\sigma+\tau \cdot i$ and $D$ is the set of real $\sigma$, we have to investigate in the interval $0<\sigma<1$ i.e. if $\operatorname{Re}( s )$ is within the open interval $(0,1)$, for which values of $\sigma$ we have non-trivial roots and the generalized sum ${ }_{ m }^{\alpha, \beta} \sigma( s )=0$ tends to zero and at the same time the function Hurwitz Zeta to tends to zero .It is naive to consider only the function and not the sum. The form of the sum will now be,\\
\[{ }_{ \text{ }\text{ } G } \zeta( m \cdot s , p )={ }_{ m }^{\alpha, \beta} \sigma( s )=\sum_{ k =1}^{\infty} \frac{1}{(\alpha \cdot k +\beta)^{ m \cdot s }}, \text { where } p =\frac{\alpha+\beta}{\alpha},\{\alpha, \beta \in R \}, s \in C\]\\
Equivalent relations functions that hold with respect to the values $\alpha, \beta$ and $p$.\\
\[\begin{aligned}
&\text { If where } p =\frac{\alpha+\beta}{\alpha},\{\alpha, \beta \in R \}, s \in C , \text { apply per case, } \\
&{ }_{ \text{ }\text{ } \mathsmaller{G} } \zeta( m \cdot s , p )= \alpha^{-s \cdot m} \cdot\zeta( m \cdot s , p ) \ (1)\\
&\zeta( m \cdot s , p )=\left(\frac{1}{ p }\right)^{ m \cdot s }+ Zeta ( m \cdot s , p +1)  \ (2)
\end{aligned}\]\\\\\\
\noindent\textbf{6.7.V.I.a. The case where $\boldsymbol{\alpha<0, \beta>0}$ and $\boldsymbol{0<\beta<-\alpha \&|p|<1, m=1}$.}\\\\
The procedure used here is the same and we solve the basic equation (2) mentioned above after zeroing it. The method is GRIM [ ] which is iterative and involves one part of the equation ,because the inverse in general of the Zeta Hurwitz function is in general difficult .In some cases it is more easily achievable and this will be shown below in the proof of the rejection of the Riemann hypothesis.\\\\
So we have relations with the series:\\\\
Starting from the original Sum and we will have:\\
\[{ }_{ \mathsmaller{G} } \zeta( m \cdot s , p )=\alpha^{- m \cdot s } \cdot \operatorname{HurwitzZeta}( m \cdot s , p )=\alpha^{- m \cdot s } \cdot\left(\left(\frac{1}{ p }\right)^{ m \cdot s }+\operatorname{Zeta}( m \cdot s , p +1)\right)=0\]\\
with the method GRIM, we take Generator $\left(\frac{1}{ p }\right)^{ m \cdot s }= u \Leftrightarrow$\\
$\{ k \in Z , u \neq 0, p \neq 0\}$ and final
$s=f(u)=\frac{2 \cdot i \cdot \pi \cdot k}{m \cdot \log (1 / p)}+\frac{\log (u)}{m \cdot \log (1 / p)}$, is the Generator function.\\\\
Therefore, it is understood that the non-trivial roots generated by this program construction are infinite and in pairs define a critical line , i.e. infinite in the interval $(0,1)$. Let us make a example, using the program matrix IV.3 in mathematica with method (GRIM[]) by infinite periodic radicals.\\\\
\textbf{Example 1. $\boldsymbol{\alpha=-3 \& \beta=2,0<\beta<-\alpha \&|p|<1, m=1}$.}\\\\
The arrangement of Table 18 is given in ascending order of $\sigma$ for the non-trivial roots of complexes
$s =\sigma+\tau \cdot i$ and of course the approximation of \textbf{the Sum} ${ }_{ \text{ } \text{ } \text{ }1 }^{\alpha, \beta} \sigma( s )=\sum_{ n =1}^{\infty} \frac{1}{(\alpha \cdot n +\beta)^s}$ is given accordingly, so that it clearly appears to be correct and therefore only then are we consistent that these roots are the required ones.
\newpage\noindent
\begin{figure}[h!]
    \centering
    \includegraphics[width=15.8cm, height=10.2cm]{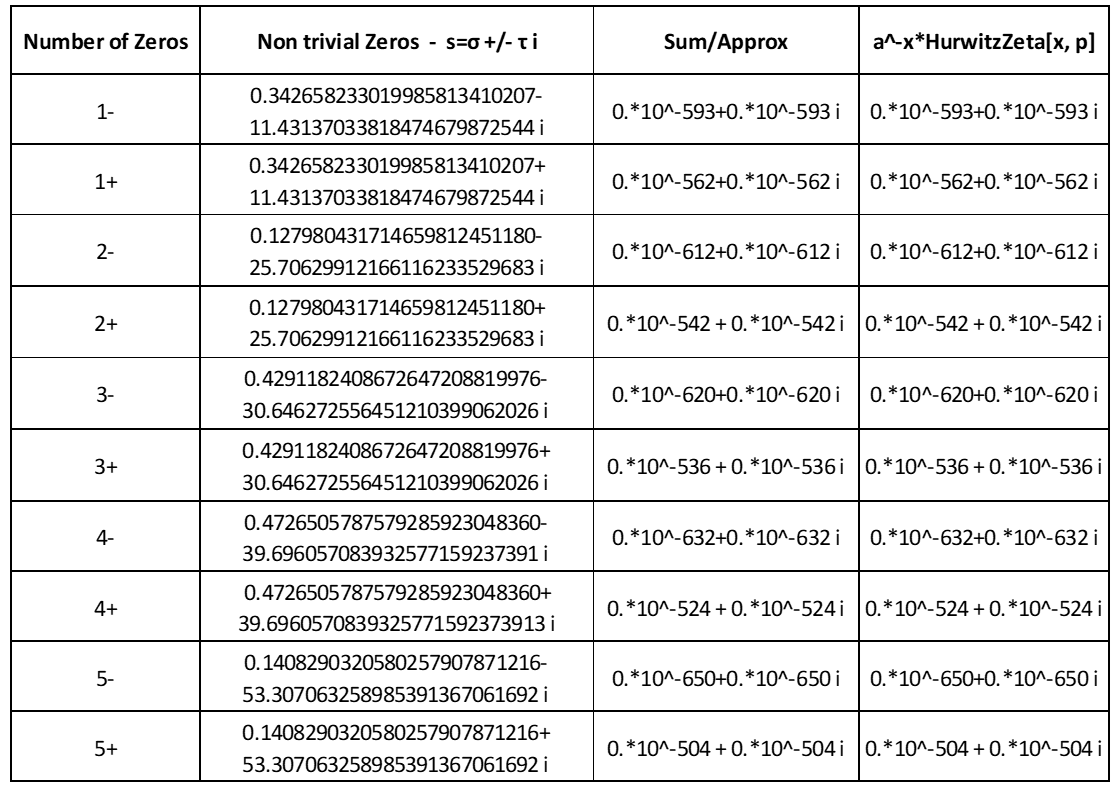}
\end{figure}\\
\noindent\textbf{Table 18. If $\boldsymbol{\alpha<0 \& \beta>>0,0<\beta<-\alpha \& | p |<1}$. Number 5 zeros for $\boldsymbol{\operatorname{Re}( s )}$ in interval $\boldsymbol{(0,1)}$}\\\\\\
\noindent Therefore, as we observe, we have a number of conjugate zeros in the interval $(0,1)$.This means that we will have scattered critical lines in pairs of complex zeros, which are certainly infinite. The percentages of conjugate zeros found by the program i.e. when d grows, increase. The sum is zero and also and \textbf{Hurwitz Zeta function}. We don't have it on the board but it's true that it tends to zero.\\\\
\textbf{Important observation}. In this case we should note the fact that the criterion applied to the specific function \textbf{Hurwitz Zeta function} for $a, b$ as we saw in chapter IV i.e. $\boldsymbol{0<\beta / \alpha+1<=1}$, when applied, a logical relation with the values that $a$ and $b$ should take arises. Therefore it equates to $\beta / \alpha>-1$ and $\beta / \alpha<=0$. It therefore follows that relation that should apply will be $0<\beta<-\alpha$. We reject the value for $\beta=0$ because we will consider it in the following. If we had $\beta=4$ this category is overturned and the finding of non-trivial roots is not valid for Hurwitz Zeta function function. But for the sum it is possible that. For values of $\beta$ outside the interval it does not mean that there are no solutions, just that they are bounded.\\\\\\
\textbf{6.7.V.I.b. The case where $\boldsymbol{\alpha<0, \beta>0}$ and $\boldsymbol{\beta>-\alpha \&|p|>0, m=1}$.}\\\\
We follow the same procedure by solving the basic equation (2) after zeroing it.Mentioned above with the GRIM [ ] method, which is iterative and concerns the one part of the equation where inversion is most feasible. We make another example, using the program matrix IV.3 in mathematica by infinite periodic radicals.
\newpage\noindent
In summary we have the relations we had before :\\\\
The generalised Sum is:\\\\
\[{ }_{ \text{ }\text{ } \mathsmaller{G} }  \zeta( m \cdot s , p )=\alpha^{- m \cdot s } \cdot \operatorname{HurwitzZeta}( m \cdot s , p )=\alpha^{- m \cdot s } \cdot\left(\left(\frac{1}{ p }\right)^{ m \cdot s }+ Zeta ( m \cdot s , p +1)\right)=0\]\\\\
with the method GRIM, we take Generator $\left(\frac{1}{p}\right)^{m \cdot s}=u \Leftrightarrow$\\\\\\
$\{ k \in Z , u \neq 0, p \neq 0\}$ and final\\\\\\
$s = f ( u )=\frac{2 \cdot i \cdot \pi \cdot k }{ m \cdot \log (1 / p )}+\frac{\log ( u )}{ m \cdot \log (1 / p )}$, is the Generator function.\\\\\\\\
\textbf{Example 2. $\boldsymbol{\alpha=-3 \& \beta=10, \beta>-\alpha \&|p|>1, m=1}$.}\\\\
In this case the order of Table 19 is given in ascending order of $\tau$ for the non-trivial roots of the complexes
$s =\sigma+\tau \cdot i$. Similarly the approximation of the sum ${ }_{ \text{ }\text{ } 1 }^{\alpha, \beta} \sigma( s )=\sum_{ n =1}^{\infty} \frac{1}{(\alpha \cdot n +\beta)^s}$ is given with a large exponent to clearly prove that it is correct. It is now obvious that the sum tends to zero. Also the column of the other form of the sum that partially tends to zero has been added using the Zeta Hurwitz function which is generally is non-zero.\\
\begin{figure}[h!]
    \centering
    \includegraphics[width=15.8cm, height=9.5cm]{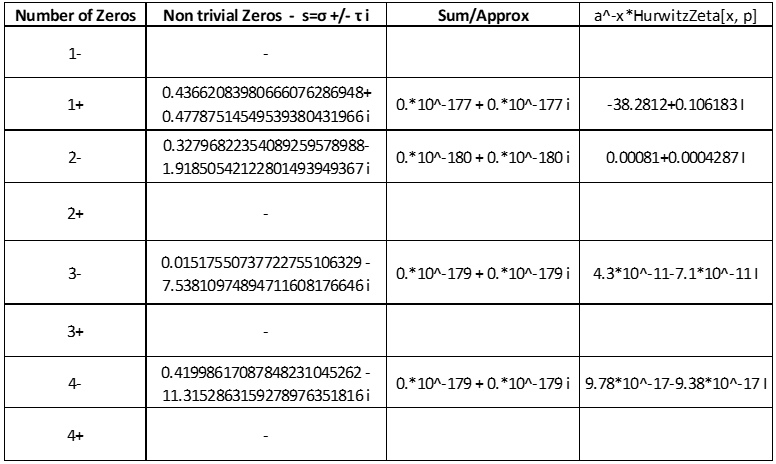}
\end{figure}\\
\vspace*{-0.5cm}
\begin{center}
\noindent\textbf{Table 19. If $\boldsymbol{\alpha<0 \& \beta>0, \beta>-\alpha \&|p|>1}$. Number 4 zeros for Re$\boldsymbol{(s)}$ in interval $\boldsymbol{(0,1)}$}    
\end{center}
\newpage\noindent
We also have a finite number of conjugate zeros in the interval (0,1). And this means that we will have scattered critical lines in pairs of conjugate zeros. We also note that the values of Im(z) that do not have opposite signs are what we call orphan complex roots.\\\\\\
\noindent\textbf{6.7.V.II. The case where $\boldsymbol{\alpha<0, \beta<0}$ and $\boldsymbol{m=1}$.}\\\\
This is the last category of the sum that can be solved. Using the same formalism and zeroing equation (2) we find infinite roots. The equations we use are the same and give us many great results. the detailed equations are again the following:\\\\
The generalised Sum is again:\\
\[{ }_{ \text{ }\text{ } \mathsmaller{G} }  \zeta( m \cdot s , p )=\alpha^{- m \cdot s } \cdot \text { HurwitzZeta }( m \cdot s , p )=\alpha^{- m \cdot s } \cdot\left(\left(\frac{1}{ p }\right)^{ m \cdot s }+ Zeta ( m \cdot s , p +1)\right)=0\]\\
with the method GRIM, we take Generator $\left(\frac{1}{p}\right)^{m \cdot s}= u \Leftrightarrow$\\\\
$\{ k \in Z , u \neq 0, p \neq 0\}$ and final
$s = f ( u )=\frac{2 \cdot i \cdot \pi \cdot k }{ m \cdot \log (1 / p )}+\frac{\log ( u )}{ m \cdot \log (1 / p )}$, is the Generator function.\\\\
We will make another typical example, using the matrix IV.3 program in mathematica with infinite periodic roots.It is clearly shown when \textbf{the Sum} is zero and also when the \textbf{Hurwitz Zeta function}.\\\\
\textbf{Example 3. $\boldsymbol{\alpha=-3 \& \beta=-1, m=1}$.}\\\\
In this case we have the Table 20 is given in ascending order of $\tau$ for the non-trivial roots of the complexes It is now obvious that both the method and the sum tend to zero. Also the column of the other form of \textbf{the Sum} which clearly tends to zero, has been added using \textbf{Hurwitz Zeta function} which is always near zero.
\begin{figure}[h!]
    \centering
    \includegraphics[width=15.8cm, height=8.5cm]{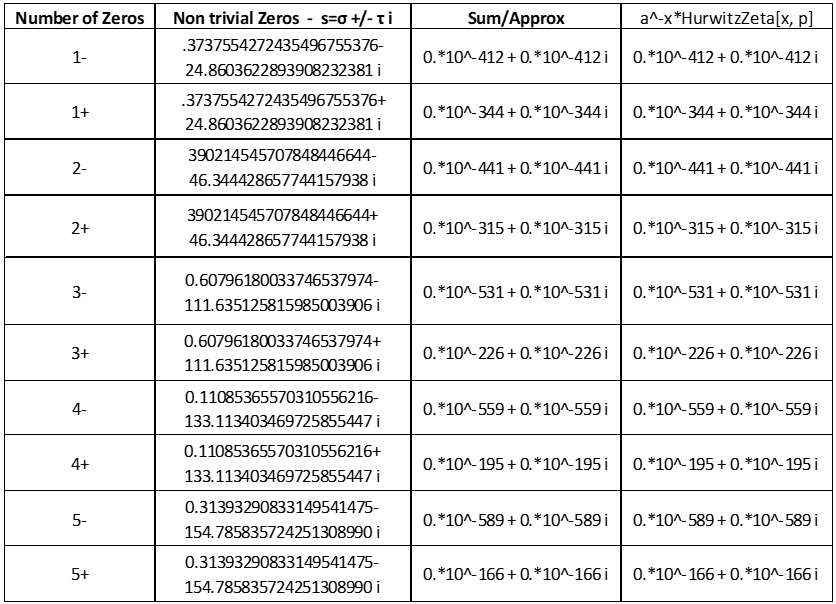}
\end{figure}\\
\vspace*{-0.5cm}
\begin{center}
 \textbf{Table 20. If $\boldsymbol{\alpha<0 \& \beta<0}$. Number 5 zeros for $\boldsymbol{\operatorname{Re}( s )}$ in interval $\boldsymbol{(0,1)}$}   
\end{center}
\newpage\noindent
We see that the zeros of this class are all with $\operatorname{Re}(z)$ within the interval $(0,1)$ per root. And this means that we will have scattered critical lines on pairs of conjugate zeros, which will be infinite.\\\\\\
\textbf{6.7.V.III. The case where $\boldsymbol{\alpha>0, \beta>0}$ or $\boldsymbol{\alpha>0, \beta<0}$ and $\boldsymbol{m=1}$.}\\\\
So we have the same procedure with respect to the equations, namely solving equation (2) by the GRIM method. Our main purpose is to approximate the solution according to the criterion that must exist and we should not forget this. That is, we need to find the non-trivial roots by applying 2 conditions: \textbf{The sum tends to zero when the series goes to infinity and the Hurwitz Zeta function also tends to zero}. This is the golden twin to accept the existence of the non-trivial roots. Thus we achieve the \textbf{Strong criterion}.\\\\\\
\textbf{6.7.V.III. $\boldsymbol\alpha$. The case where $\boldsymbol{\alpha>0, \beta>0, m=1}$}\\\\
Here there is no particular problem in finding the roots we follow the procedure of program IV.a. In this case we form Table 21 if we remove the criterion $0<\operatorname{Re}(s)<1$ from the logic check of the program. The table is given in ascending order $\tau$ for the non-trivial roots of the complex. As we will see we get values for $\operatorname{Re}(s)>1$ which is expected . So we do not have in this case $0<\operatorname{Re}(s)<1$.\\\\
\noindent\textbf{Example 3. $\boldsymbol{\alpha=5 \& \beta=6, m=1}$.}\\
\begin{figure}[h!]
    \centering
    \includegraphics[width=16.5cm, height=10.8cm]{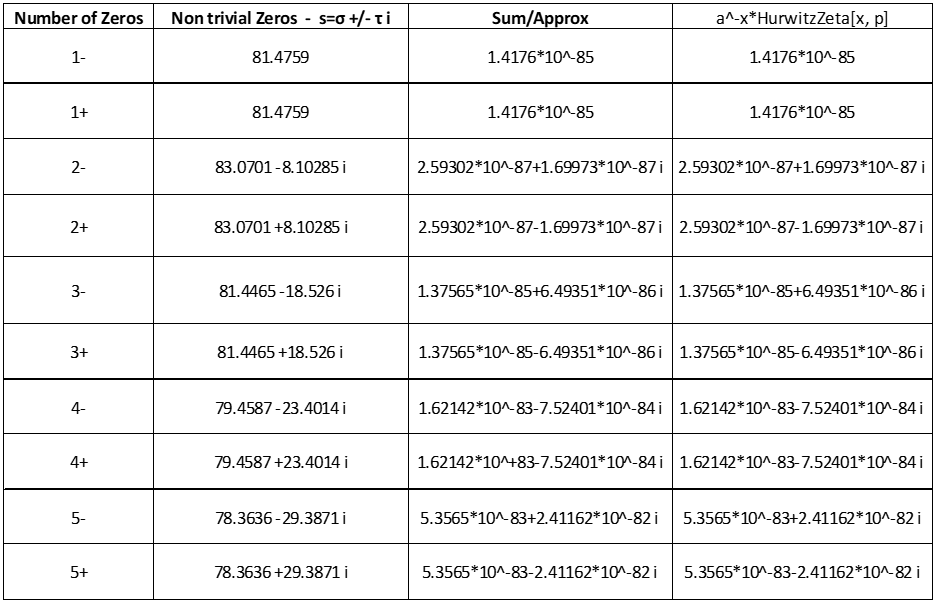}
\end{figure}\\
\vspace*{-0.5cm}
\begin{center}
  \textbf{Table 21. If $\boldsymbol{\alpha>0 \& \beta>0}$. Number 5 zeros for $\boldsymbol{\operatorname{Re}( s )}$ outside from interval $\boldsymbol{(0,1)}$}  
\end{center}
\newpage\noindent
What is clear is that all the roots are outside the interval $(0,1)$ with positive $\sigma$ in the complex value of $s$.\\\\
Also, as we can see, the values of the sum and the the \textbf{Hurwitz Zeta function} absolutely coincide.\\\\\\
\textbf{6.7.V.III. $\boldsymbol{\beta}$. The case where $\boldsymbol{\alpha>0, \beta <0, m =1}$}\\\\
Here we distinguish two cases according to Important observation we mentioned in example 1.Here we will have the scheme $0<\beta / \alpha+1<=1$ with $\alpha>0 \& \beta<0$. This results in the relation $-\alpha<\beta<0$. For values of $\beta$ outside the interval it does not mean that there are no solutions, just that they are bounded.\\\\\\
\noindent\textbf{Example 4. $\boldsymbol{\alpha=5 \& \beta=-1, m=1}$.}
\begin{figure}[h!]
    \centering
    \includegraphics[width=16.8cm, height=10.2cm]{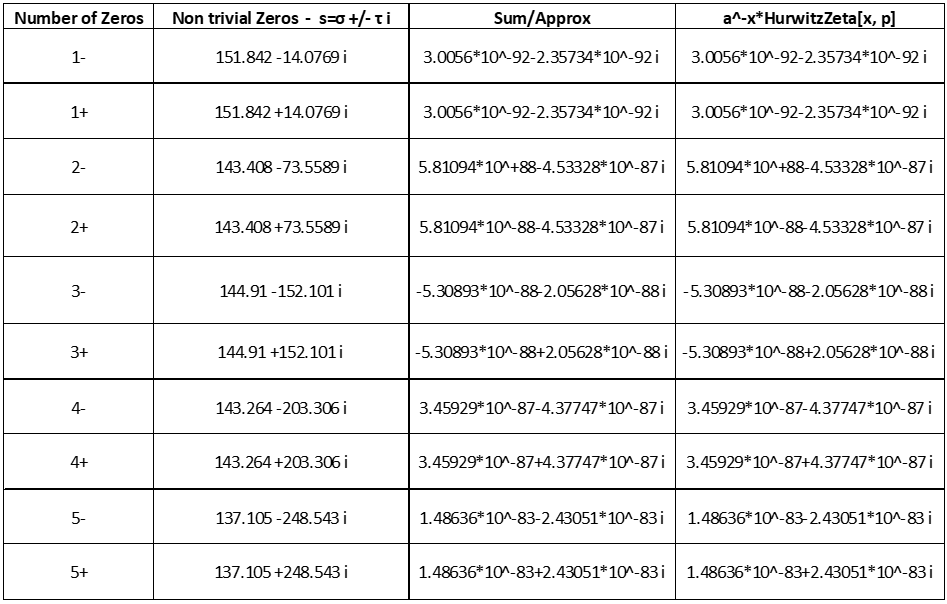}
\end{figure}\\
\vspace*{-0.5cm}
\begin{center}
    \textbf{Table 22. If $\boldsymbol{\alpha>0 \& \beta<0}$. Number 5 zeros for $\boldsymbol{\operatorname{Re}( s )}$ outside from interval $\boldsymbol{(0,1)}$}
\end{center}
\vspace{\baselineskip}
As we can see we are considering a value $\beta=-1$ that is within the range \textbf{we accept, i.e. $\boldsymbol{- \alpha < \beta < 0}$}. We do not accept zero for now because it will be considered later.It is the range that theoretically corresponds to the function $\zeta(1)$. It is also clear that all the roots are outside the interval and $(0,1)$ with positive $\sigma$ in the complex value of $s$. Also, the values of the sum and the Zeta Hurwitz function coincide perfectly. If we consider the case where $\beta$ is outside the interval $(-\alpha, 0)$ what we will observe is that we will have bounded values for the trivial zeros for the negative values of $k$ in the program.
\newpage\noindent
\textbf{6.7.IX. The paradox in the Hypothesis of Riemann is a strong pre-Hypothesis that raises the question that the Hypothesis of Riemann, it should be rejected by default.}\\\\
Let us come to the definition of the case: \textbf{In mathematics, the Riemann hypothesis is the conjecture that the Riemann zeta function has its zeros only at the negative even integers and complex numbers with real part $\boldsymbol{1 / 2}$}. But this directly implies that the same is true for the sum of the harmonic numbers. i.e. If the Riemann hypothesis holds, then both the Generalized Sum ${ }_m^{\alpha, \beta} \sigma(s)=\sum_{n=1}^{\infty} \frac{1}{(\alpha \cdot n+\beta)^{m \cdot s}} \rightarrow 0, s \in C$ and $\zeta(s) \rightarrow 0$, will they have the same non-trivial roots for $s$ where $\operatorname{Re}( s )=1 / 2$ and since $\{\beta=0, \alpha=1$ and $m=1\}$ for the \textbf{Generalized Sum} ${ }_{\text{ }\text{ }\text{ } 1}^{\alpha, \beta} \sigma(s)$.\\\\
There are 3 cases that raise the paradox that the $\zeta()$ function does not verify the finding of trivial roots at the same time with the sum ${ }_1^{1,0} \sigma(s)$ that they are supposed to be computationally identical. This means that the non-trivial zeros of $\zeta(z)=0$, the ones we already know, do not zero the sum of the harmonic series we mention and for which it is supposed to have been constructed. And therefore it will hold for the already known ZetaZeros the relation:\\\\
If $s_i=\zeta_{\text {zeecos }}(i)$ and $\zeta\left(s_i\right)=0$ or tending to zero $\Rightarrow{ }_1^{1,0} \sigma\left(s_i\right) \neq 0$ or not tending to zero by $s_i \in C$ and $\operatorname{Re}\left(s_i\right)=\frac{1}{2}$.\\\\\\
\textbf{6.7.IX.1. Sum comparison method and Zeta Hurwitz function.}\\\\
The first case, which is related to the Zeta Hurwitz function, concludes that the sum given $a=-2$ and $b=1$ is equal to a representation associated with $\zeta(s)$, so if $\zeta(s)$ is zero then the sum we are referring to will necessarily be zero. But if this is the case we will have the following relations\\\\
If $\zeta( s )=0 \Rightarrow{ }_{\text{ }\text{ }\text{ }\text{ }1}^{-2,1} \sigma( s )=\sum_{ n =1}^{\infty} \frac{1}{(-2 \cdot n +1)^{ s }} \rightarrow 0$ and ${ }_1^{1,0} \sigma( s )=\sum_{ n =1}^{\infty} \frac{1}{( n )^{ s }} \rightarrow 0, s \in C$ and $\operatorname{Re}( s )=\frac{1}{2}$\\\\
i.e two differential sums will have same the nontrivial zeros, which is a contradiction. Therefore, since it is proved by the criterion mentioned earlier that this sum has computationally nontrivial zeros of $\zeta(s)=0$, it follows that the nontrivial zeros of $\zeta(s)=0$ cannot be of the sum to which we know i.e of ${ }_1^{1,0} \sigma(s)$ it refers and therefore the zeros belong to the sum ${ }_{\text{ }\text{ }\text{ }\text{ }1}^{-2,1} \sigma( s )$.\\\\
\textbf{Lemma:} We will prove that ${ }_{\text{ }\text{ }\text{ }\text{ }1}^{-2,1} \sigma( s )=\sum_{n=1}^{\infty} \frac{1}{(-2 \cdot n+1)^s}=(-2)^s \cdot\left(-1+2^s\right) \cdot \zeta(s)$ has the zeros of $\zeta(s)=0$.\\\\\\
\textbf{Proof}\\\\
It follows from Chapter $V$ that\\
\[{ }_{ \mathsmaller{G} } \zeta( m \cdot s , p )={ }_{ m }^{\alpha, \beta} \sigma( s )=\sum_{ k =1}^{\infty} \frac{1}{(\alpha \cdot k +\beta)^{ m \cdot s }}, \text { where } p =\frac{\alpha+\beta}{\alpha},\{\alpha, \beta \in R \}, s \in C\]\\
We need to find functions of equivalent relations that hold with respect to the values $\alpha, \beta$ and $p$.\\\\
If where $p =\frac{\alpha+\beta}{\alpha},\{\alpha, \beta \in R \}, s \in C$, apply per case,\\
\[\begin{aligned}
&{ }_{ \mathsmaller{G}} \zeta( m \cdot s , p )=\alpha^{- m \cdot s } \cdot \zeta( m \cdot s , p ) \\[5pt]
&\zeta( m \cdot s , p )=\left(\frac{1}{ p }\right)^{ m \cdot s }+ Zeta ( m \cdot s , p +1)
\end{aligned}\tag{2}\]
\newpage\noindent
Knowing relation (2) and for our case where $a=-2$ and $b=1$ and $m=1$ we have the relation\\
\[\begin{aligned}
&{ }_{ \mathsmaller{G} } \zeta(1 \cdot s , p )=\alpha^{- s } \cdot \zeta( s , p )=\alpha^{- s } \cdot\left(\left(\frac{1}{ p }\right)^{ s }+\zeta( s , p +1)\right)=(-2)^{- s } \cdot\left((2)^{ s }+\zeta\left( s , \frac{3}{2}\right)\right) \Rightarrow \\[5pt]
&{ }_{ \mathsmaller{G} } \zeta\left( s , \frac{1}{2}\right)={ }^{-2,1} \sigma( s )=\sum_{ k =1}^{\infty} \frac{1}{(-2 \cdot k +1)^{ s }}=(-2)^{-s} \cdot\left((2)^{ s }+\zeta\left( s , \frac{3}{2}\right)\right)
\end{aligned}\tag{3}\]\\
But apply also\\
\[\zeta\left( s , \frac{3}{2}\right)=-2^{ s }+\left(-1+2^{ s }\right) \cdot \zeta( s )\tag{4}\]
From relations $(3 \& 4)$ it follows that\\
\[{ }_{\mathsmaller{G}} \zeta\left( s , \frac{1}{2}\right)={ }_1^{-2,1} \sigma( s )=\sum_{ k =1}^{\infty} \frac{1}{(-2 \cdot k +1)^{ s }}=(-2)^{-s} \cdot\left(2^s+-2^{ s }+\left(-1+2^{ s }\right) \cdot \zeta( s )\right)=(-2)^{- s } \cdot\left(-1+2^{ s }\right) \cdot \zeta( s )\tag{5}\]
It can therefore be concluded that\\
\[\sum_{ k =1}^{\infty} \frac{1}{(-2 \cdot k +1)^{ s }}=(-2)^{-s} \cdot\left(-1+2^{ s }\right) \cdot \zeta( s )\tag{6}\]\\
Because $(-2)^{-s}$ is not zeroed and because $s$ is not zero therefore, all that remains is $\zeta(s)$ becomes zero, the sum $\sum_{ k =1}^{\infty} \frac{1}{(-2 \cdot k +1)^s}$ will also become zero which means that the known sum $\sum_{ k =1}^{\infty} \frac{1}{( k )^s}$ is not reduced to zero. Otherwise we would have 2 completely different sums with the same non-trivial zeros which is impossible. In the following cases we will prove computationally that the sum $\sum_{ k =1}^{\infty} \frac{1}{( k )^{ s }}$ with zeros of $\zeta( s )=0$ is not zero.\\\\
The verification that the zeros of $\zeta(s)=0$ coincide with the zeros of the sum is done with program IV.a by checking directly back to the sum. In this case we use the the intermediate function\\
\[\zeta\left( s , \frac{3}{2}\right)=-2^s+\left(-1+2^s\right) \cdot \zeta(s) \text { which we transform to the form } \frac{\zeta\left( s , \frac{3}{2}\right)}{1-2^\sigma}=\frac{2^s}{\left(-1+2^s\right)}-\zeta( s )\]\\
A corresponding program that produces the ZetaZeros is the following in mathematica 12. The data are $d=10, a=-2, b=1, \operatorname{In} A p=-1 / 1000, R=1, A \operatorname{sum}=10^{\wedge}(-10), m=1$\\\\
Off[FindRoot::lstol,N::meprec,Sum::div,General::munfl,NIntegrate::slwcon,NIntegrate::ncvb]; \\\\
\textbf{'program data'}\\\\
$d =$ Input["Range values of $k , d =$ "];\\[5pt]
$a=$ Input["Value of $a=$"];\\[5pt]
$b=\operatorname{Input}[$ "Value of $b=$ "];\\[5pt]
InAp=Input["Initial Approx of Iteration, $\operatorname{In} A p=$ "];\\[5pt]
$R =$ Input["Replays of Iterations, $R =$ "];\\[5pt]
Asum=Input["Approach value of Sum, Asum="];\\[5pt]
$m =$ Input["Multiplicity of exponent of $m =$ "];
\newpage\noindent
\textbf{'program execution'}\\\\
$p :=( a + b ) / a$\\[5pt]
$h2\left[ x _{-}\right]$:=Zeta $[ x ]$;\\[3pt]
$h \left[ x _{-}\right]:=2^{ x } /\left(-1+2^{\wedge} x \right)-\operatorname{Zeta}[ x ]$;\\[3pt]
$h3 \left[x_{-}\right]:=(-2)^{\wedge}(-x)^* \operatorname{Zeta}[ x , p ] ; j:=0$;\\[3pt]
$h4 \left[ x _{-}\right]:= N \left[\sum_{ k =1}^{\infty} \frac{1}{( a \cdot k + b )^{\wedge}( mx )}\right]$;\\[3pt]
For $[ k =- d , k <= d , k ++$,\\[3pt]
$f\left[ u _{-}\right] :=N[(2$ I $\pi k ) / \log [2]+\log [u /(-1+u)] / \log [2]] ;$\\[3pt]
$x r =$ Nest $[ f [ h 2[\#]] \& , I n A p , R ] ;$\\[3pt]
$FQ 2= N [ xr , 10]$;\\[3pt]
$s = N [ y 1,40] /$.FindRoot $[ h 4[ y 1]==0,\{ y 1, FQ 2\}$, WorkingPrecision $\to 280]$;\\[3pt]
If $[0<\operatorname{Re}[ s ]< 1 \& \& A b s [ h 4[ s 2 ]]< A s u m , j = j +1$;\\[3pt]
Print["x(",k,")=",s1=N[y2,40]/.FindRoot[h4[y2]==0, $\{ y 2, s \}$, WorkingPrecision$\to 180]$,\\[3pt]
",", N[h[s1]],",",N[h4[s1]],",",N[h3[s1]]];\\[3pt]
$j = j +1 ;$;Loopback $]$ ]\\[3pt]
Print["Percent of Zeros $\%=$ ", $\left.N \left[ j /(2 d +1)^* 100\right]\right]$\\\\
The program itself simulates the equation to find the zero prime of the sum in two ways, i.e. using the Sum and using the function Zeta $(x, 1 / 2)$. Both must agree, since they are equal as we have proved.\\\\
\noindent\textbf{Example 4. $\boldsymbol{\alpha=-2 \& \beta=1}$.}\\\\
\textbf{Table 23} gives the known non-trivial zeros agreeing with the sum $\sum_{ k =1}^{\infty} \frac{1}{(-2 \cdot k +1)^s}$ and $\zeta( x )$ when zeroed.
\begin{figure}[h!]
    \centering
    \includegraphics[width=13cm, height=8.8cm]{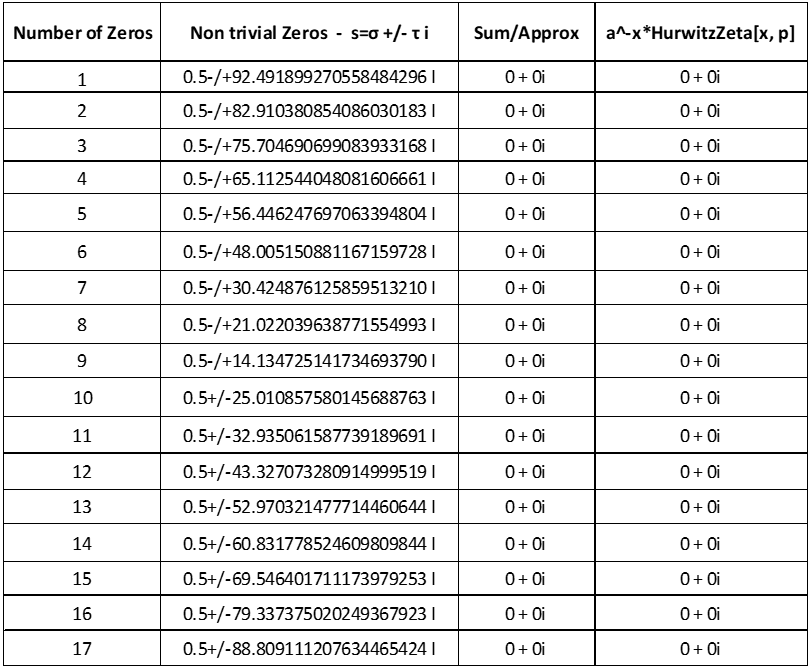}
\end{figure}\\
\vspace*{-0.5cm}
\begin{center}
    \textbf{Table 23. If $\boldsymbol{\alpha=-2 \& \beta=1}$ check the values for the Sum and $\boldsymbol{(-2)^{\wedge}(-s) * \zeta(x, 1 / 2)}$}
\end{center}
\newpage\noindent
We observe according to the program that the non-trivial roots that are found with program agree to be the Sum at infinity zeroed as well as $\zeta(x, 1 / 2)=0$. Therefore they agree perfectly . Of course, the sum will also be zeroed by the remaining non-trivial zeros of $\zeta(x)$.\\\\
\textbf{6.7.IX.II. Calculation method using the program.}\\\\
We will try to calculate the \textbf{non trivial} zeros with a different method. We will use the program VI.3 and with the  technique to have in the Sum the form. So we will see where the values for the \textbf{non trivial} roots go, if we add and subtract a constant number a very small. Both cases are useful and show us the kind of \textbf{non trivial roots}.
\[\sum_{ k =1}^{\infty} \frac{1}{\left( k \pm \frac{1}{\varepsilon}\right)^{ s }}, \varepsilon>>1\]
So we will see where the values for the trivial roots go, if we add and subtract a constant number a very small. Both cases are useful and show us the kind of trivial roots.\\\\
\textbf{1st Category } ${\sum_{ k =1}^{\infty} \dfrac{1}{\left( k +\frac{1}{\varepsilon}\right)^{ s }}, \varepsilon>>1}$\\\\
Using the program we will see for large $\varepsilon=5^* 10^{\wedge} 4$ which non-trivial zeros are produced. What we will see is disappointing because while the equation is zeroed and $\zeta(s, p)$, but is not zeroed the Sum. We also observe large values for $\operatorname{Re}( s )$ as well as $\operatorname{Im}( s )$ of the non-trivial roots and $\operatorname{Re}( s )<0$.\\\\
Off[FindRoot::bddir,FindRoot::jsing,Power::infy,Infinity::indet,General::cvmit,\\
ReplaceAll::reps,Divide::inf y,FindRoot::1stol,Sum::div, General::stop],General::munfl];\\\\
\textbf{'program data'}\\\\
$d =\operatorname{Input}\left[\right.$ "Range values of $\left.k , d =^{\prime \prime}\right]$\\[5pt]
$a =\operatorname{Input}[$ "Value of $a =$ "];\\[5pt]
$b =$ Input["Value of $b =$ "]\\[5pt]
InAp $=\operatorname{Input}[$ "Initial Approx of Iteration, $\operatorname{In} A p=$ "]\\[5pt]
$R =\operatorname{Input}[$ "Replays of Iterations, $R=$"];\\[5pt]
Asum=Input["Approach value of Sum, Asum="];\\[5pt]
$m =$ Input["Multiplicity of exponent of $m =$"]\\\\
\textbf{'program execution'}\\\\
$p :=( a + b ) / a$;\\[5pt]
$h2\left[ x _{-}\right]:=-Zeta[m*x,(1+p)]$;\\[5pt]
$h \left[ x _{-}\right]:= a ^{- m ^* x *}\left((1 /( p ))^{\wedge}\left( m ^* x \right)+\operatorname{Zeta}^2\left[ m ^* x ,(1+ p )\right]\right)$;\\[5pt]
h3[x $]:= a ^{- m ^* x *}\left(\right.$ HurwitzZeta $\left.\left[ m ^* x , p \right]\right) ; j :=0$;\\[5pt]
$h 4\left[ x _{-}\right]:= N \left[\sum_{k=1}^{\infty} \frac{1}{(a \cdot k+b)^{\wedge}(m x )}\right] ;$\\[5pt]
For $[ k =- d , k <= d , k ++$,\\[5pt]
$f [ u ]]:= N \left[(2 I \pi k +\log [ u ]) / \log \left[1 / p ^{\wedge} m \right]\right] ;$\\[5pt]
xr= $=$ Nest $[f[h 2[\#]] \&, 1 / 1000,1] ;$\\[5pt]
$FQ 2= N [ xr , 10]$;\\[5pt]
$s = N [ y 1,40] /$ FindRoot $[ h 4[ y 1]==0,\{ y 1, FQ 2\}$, WorkingPrecision $\to 180]$;\\[5pt]
Print $[$ "x $(", k, ")=", s 1=N[y 2,40] /$. FindRoot $[ h 4[ y 2]==0,\{ y 2, s \}$, WorkingPrecision $\to 120]$,\\[5pt]
",",N[h[s1]],",",N[h4[s1]],",",N[h3[s1]]];]\\[5pt]
Print[p,"Percent of Zeros $\left.\%=", N\left[j /(2 d +1)^* 100\right]\right]$\\\\
A corresponding program that produces the ZetaZeros is the following in mathematica 12. The data are $d =5, a =1, b =1 / 50.000, \operatorname{InA} p =1 / 1000, R =1, m =1$.\\
\begin{figure}[h!]
    \centering
    \includegraphics[width=12cm, height=5.2cm]{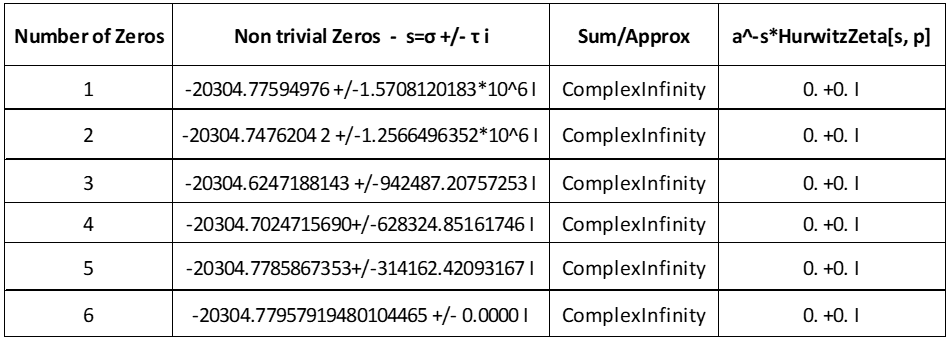}
\end{figure}\\
\vspace*{-0.5cm}
\begin{center}
    \textbf{Table 24. If $\boldsymbol{\alpha=1 \& \beta=1 / 5^* 10^{\wedge} 4}$ check the values for the Sum and $\boldsymbol{\zeta(s, 1+b / a)}$}
\end{center}
\vspace{0.6\baselineskip}
It is clear from the results in \textbf{Table 24} that the Sum $\sum_{ k =1}^{\infty} \frac{1}{( k )^{ s }}$ is not zero but tends to infinity in complex form. The non trivial zeros retain very large values and $\operatorname{Re}( s )$ is outside the interval $(0,1)$ and negative. This means it does not have the roots of $\zeta(s)=0$, which we all know.\\\\
\textbf{2st Category} $\sum_{ k =1}^{\infty} \frac{1}{\left( k -\frac{1}{\varepsilon}\right)^{ s }}, \varepsilon>>1$\\\\
Using again the previous program we will see for large $\varepsilon=5^* 10^{\wedge} 4$ which non-trivial zeros are produced. We will also see that while the equation is zeroed and $\zeta(s, p)$, the sum is not zeroed; we also observe large values for $\operatorname{Re}( s )$ as well as $\operatorname{Im}( s )$ of the non trivial roots and $\operatorname{Re}( s )<0$.\\
Also the same program that produces the ZetaZeros is the following in mathematica 12.The data are $d =5, a =1, b =-1 / 50.000, \operatorname{In} A p=1 / 1000, R =1, m =1$.\\
\begin{figure}[h!]
    \centering
    \includegraphics[width=11.5cm, height=5.3cm]{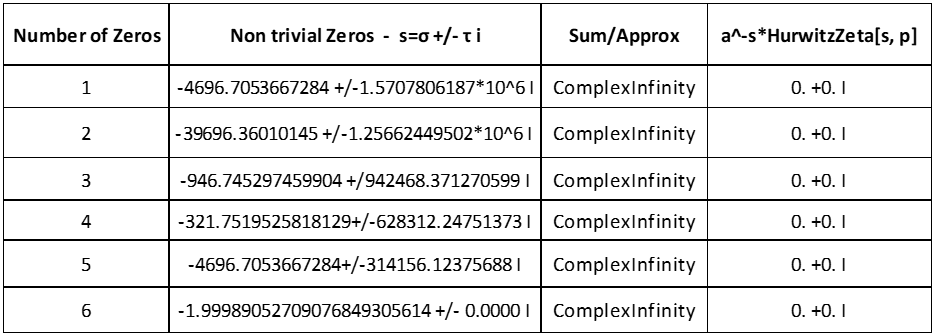}
\end{figure}\\
\vspace*{-0.5cm}
\begin{center}
\textbf{Table 25. If $\boldsymbol{\alpha=1 \& \beta=-1 / 5^* 10^{\wedge} 4}$ check the values for the Sum and $\boldsymbol{\zeta(s, 1+b / a)}$}
\end{center}
\newpage\noindent
\vspace{\baselineskip}
\noindent Here too it is obvious that the sum $\displaystyle\sum_{ k =1}^{\infty} \frac{1}{( k )^s}$ is not zero but tends to infinity in complex form. The trivial zeros retain large values and $\operatorname{Re}( s )$ is outside the interval $(0,1)$ and negative. We reaffirm that it does not have the roots $\zeta(s)=0$, which we all know, and that final up the sum lacks trivial zeros.\\\\\\\\
\textbf{6.7.IX.III. Approximate Check by substituting the non-trivial roots, in the infinite sums $\boldsymbol{\displaystyle\sum_{ k =1}^{\infty} \frac{1}{( k )^s}}$ and $\boldsymbol{\displaystyle\sum_{ k =1}^{\infty} \frac{1}{(-2 \cdot k +1)^s}}$}\\\\\\\\
We will do the final check with the sum in mathematica 12. We note that in the sum $\displaystyle\sum_{ k =1}^{\infty} \frac{1}{(-2 \cdot k +1)^s}$ the approximation for a series of known trivial zeros approaches close to zero while for the sum of the $\displaystyle\sum_{ k =1}^{\infty} \frac{1}{( k )^s}$ for the non trivial zeros it does not zero, on the contrary, the result is infinite. We construct the program below and check for the first 10 non-trivial zeros. The results in Table 26 show exactly what we mentioned\\\\
So we have the corresponding program that produces the ZetaZeros of $\zeta( s )=0$ and we find the corresponding values for Sums. The data are $a 1=-2, b 1=1, a 2=1, b 2=0, m =1$.\\\\
Off[FindRoot::1stol,N::meprec,Sum::div,General::munfl,NIntegrate::slwcon,NIntegrate::ncvb]; \\\\
\textbf{'program data'}\\\\
al=Input["Value of a $=$ "']\\[5pt]
b1=Input["Value of bl="];\\[5pt]
a2 $=$ Input["Value of a2="];\\[5pt]
b2=Input["Value of b2="];\\[5pt]
$m =$ Input["Multiplicity of exponent of $m =$ "']\\\\
\textbf{'program execution'}\\\\
$h 1\left[ x _{-}\right]:= N \left[\sum_{k=1}^{\infty} \frac{1}{\left( a _1 \cdot k+b_1\right)^{\wedge}(m \cdot x)}, 10\right] ;($Sum1$)$\\[5pt]
$h 2\left[ x _{-}\right]:= N \left[\sum_{ k =1}^{\infty} \frac{1}{\left( a _2 \cdot k + b _2\right)^{\wedge}( m \cdot x )}, 10\right] ;($Sum2$)$\\[5pt]
For $[ i =1, i <=15, i ++$,\\[5pt]
$xn = N [$ ZetaZero[i],550];\\[5pt]
Print[i,",",N[xn,10],h1[xn],h2[xn]]]
\newpage\noindent
\begin{figure}[h!]
    \centering
    \includegraphics[width=14.5cm, height=7.3cm]{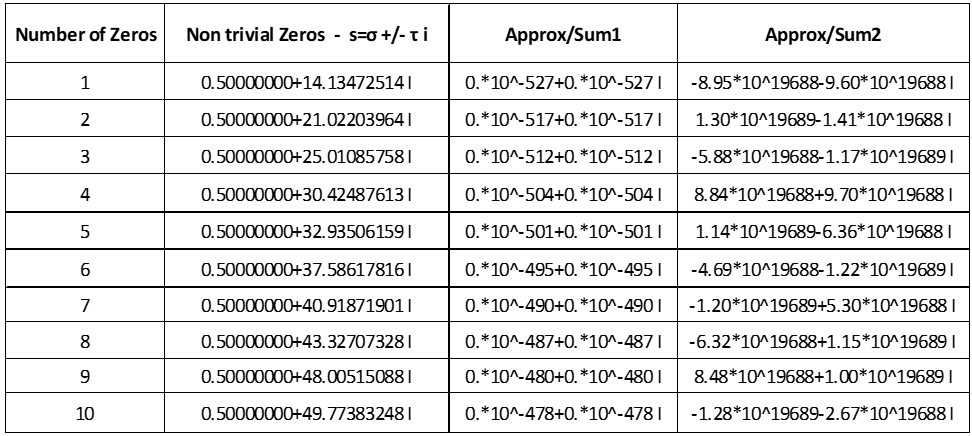}
\end{figure}\\
\textbf{Table 26. If $\boldsymbol{(\alpha 1=-2 \& \beta 1=1)}$ and $\boldsymbol{(\alpha 2=1 \& \beta 2=0)}$ approximation for the Sum1 and Sum2}\\\\\\
As we see the results of approximate sums that Sum1 $=\sum_{k=1}^{\infty} \frac{1}{(-2 \cdot k+1)^x}$ tends to zero while the Sum2 $=\sum_{k=1}^{\infty} \frac{1}{k^x}$ it tends to infinity. For Sum1 it is as expected. But the Sum2 makes the overturn and the unexpected.
The question is not simply whether Riemann's Hypothesis is true or not, but which functions $\zeta()$ we accept to be included in the Hypothesis, as Riemann defined it, and whether the function $\zeta(s)$ with its non trivial zeros represents the sum $\sum_{ k =1}^{\infty} \frac{1}{( k )^5}$ or some other sum if it is zeroed. Also what is the set of existing functions $\zeta()$ to be included in the Hypothesis. If we take as a general definition of Riemann's hypothesis all similar $\zeta()$ functions then it is not fully valid as a test of functions with s zeros satisfying the constraint $0<\operatorname{Re}[ s ]<1$ and especially when it constrains us to be $\operatorname{Re}(s)=1 / 2$ and also what is the sum represented by the hypothesis and if it is what is implied by the simple harmonic series. This I believe is the complete conundrum of the definition of the Riemann hypothesis that has revolutionized the mathematical community for over 160 years.\\\\\\
\textbf{6.7.X. Rejection of the Riemann hypothesis. Reconsideration of the sum corresponding to the $\boldsymbol\zeta$ function.}\\\\
In an attempt to overturn the paradox it will be necessary to examine more deeply and computationally all the cases where we have $p=1 / 2$ or $p=1$. First we will investigate one more sum that is very important and will play an important role in understanding which this sum is that ultimately represents the $\zeta()$ function.\\\\\\
\textbf{X.I Generalisation of the sum} ${ }_{ \quad\text{ } m }^{-2 \beta, \beta} \sigma( s )=\frac{1}{\beta^{m \cdot s}} \sum_{n=1}^{\infty} \frac{1}{(-2 \cdot n+1)^{m \cdot s}}$\\\\
The Hurwitz zeta function the generalized by\\
\[{ }_{ \text{ }\text{ } \mathsmaller{G} }  \zeta(m \cdot s, p)={ }_m^{\alpha, \beta} \sigma(s)=\sum_{k=1}^{\infty} \frac{1}{(\alpha \cdot k+\beta)^{m \cdot s}}, \text { where } p=\frac{\alpha+\beta}{\alpha},\{\alpha, \beta \in R-0\}, s \in C\]
\newpage\noindent
We want to have the form for $p=(\alpha+\beta) / \alpha=1 / 2$ i.e. $\alpha=-2 \beta$ and therefore\\
\[{ }_{ \quad\text{ } m }^{-2 \beta, \beta} \sigma( s )=\frac{1}{\beta^{m \cdot s}} \sum_{n=1}^{\infty} \frac{1}{(-2 \cdot n+1)^{m \cdot s}}=\frac{1}{(-\beta)^{m \cdot s}} \cdot(2)^{-m \cdot s} \cdot\left(-1+2^{m \cdot s}\right) \cdot \zeta(m \cdot s)\tag{7}\]\\
If $m =1$,$ \beta>0 $ we have\\
\[{ }_{ \quad\text{ } m }^{-2 \beta, \beta} \sigma( s )=\frac{1}{\beta^s} \sum_{n=1}^{\infty} \frac{1}{(-2 \cdot n+1)^s}=\frac{1}{(-\beta)^s} \cdot(2)^{-s} \cdot\left(-1+2^s\right) \cdot \zeta(s)\tag{8}\]\\
Alternatively we can use the digamma function\\
\[\psi(q)=\frac{d}{d q} \log \Gamma(q)\tag{9}\]\\
Appears in Laurent expansion of \\
\[\zeta(z, q)=\frac{1}{z-1}-\psi(q)+O(z-1)\tag{10}\]\\
And the polygamma function defined by\\
\[\begin{aligned}
&\psi^m(q)=\frac{d^m}{d q^m} \psi(q), m \in N \\[10pt]
&\psi^m(q+1)=\psi^m(q)+\frac{(-1)^m m !}{q^{m+1}}
\end{aligned}\tag{11}\]
2 very useful relations apply\\
\[\begin{aligned}
&\psi\left( n +\frac{1}{2}\right)=\sum_{k=1}^{n-1} \frac{1}{ k }+\sum_{ k = n }^{2 n +1} \frac{2}{ k }-\log (4)-\gamma, n \in N \\[10pt]
&\psi\left(\frac{1}{2}\right)=-\gamma-\log (4)
\end{aligned}\tag{12}\]\\\\
And final we have if $s=1$ will be take for Sum\\
\[\begin{aligned}
&{ }_{ \quad\text{ } 1 }^{-2 \beta, 2\beta} \sigma( s )=\frac{1}{\beta^1} \sum_{n=1}^{\infty} \frac{1}{(-2 \cdot n+1)^1}=\lim _{n \rightarrow \infty} \frac{1}{2 \cdot \beta}(-\gamma-\log (4)-\operatorname{PolyGamma}(0,1 / 2+n)) \rightarrow-\infty\\[10pt]
& \gamma = \text{ Eulergamma }
\end{aligned}\tag{13}\]\\\\\\
\textbf{X.II Generalisation of the sum} ${ }_{ \quad\text{ }\text{ }\text{ }\text{ } m }^{-2 \beta, -\beta} \sigma( s )=\frac{1}{\beta^{m \cdot s}} \sum_{n=0}^{\infty} \frac{1}{(-2 \cdot n-2)^{m \cdot s}}$\\\\
The Hurwitz zeta function the generalized by\\
\[{ }_{ \text{ }\text{ } \mathsmaller{G} }  \zeta( m \cdot s , p )={ }_{ \quad\text{ } m }^{\alpha, -2 \beta} \sigma( s )=\sum_{ k =0}^{\infty} \frac{1}{(\alpha \cdot k -2 \beta)^{ m \cdot s }} \text {, where } p =-\frac{2\beta}{\alpha},\{\alpha, \beta \in R -0\}, s \in C\]\\
We want to have the form for $p=(-2 \beta) / \alpha=1$ i.e. $\alpha=-2 \beta$ and therefore\\
\[{ }_{ \qquad\text{ } m }^{-2 \beta, -2 \beta} \sigma( s )=\frac{1}{(-\beta)^{m \cdot s}} \sum_{n=0}^{\infty} \frac{1}{(2 \cdot n+2)^{m \cdot s}}=\frac{1}{(-\beta)^{m \cdot s}} \cdot(2)^{-m \cdot s} \cdot \zeta(m \cdot s)\tag{14}\]\\
If $m =1$ we have\\
\[{ }_{ \qquad\text{ } 1 }^{-2 \beta, -2\beta} \sigma( s )=\frac{1}{(-\beta)^s} \sum_{n=0}^{\infty} \frac{1}{(n+1)^s}=\frac{1}{(-\beta)^s} \cdot(2)^{-s} \cdot \zeta(s)\tag{15}\]\\
1 very useful relations apply\\
\[\psi(n+2)=-\gamma-\frac{1}{2+n}+(2+n) \sum_{k=1}^{\infty} \frac{{ }_2 \tilde{F}_1\left(1,2 ; 2 ;-\frac{2+n}{k}\right)}{k^2}\tag{16}\]\\
where $2+\operatorname{Re}(n)>0$\\\\
And final we have if $s=1$ and $m=1$, $\beta>0$ will be take for Sum\\
\[\begin{aligned}
&{ }_{ \qquad\text{ } 1 }^{-2 \beta, -2\beta} \sigma( s )=\frac{1}{(-2 \beta)^1} \sum_{n=1}^{\infty} \frac{1}{(n+1)^1}=\lim _{n \rightarrow \infty} \frac{1}{2 \cdot \beta}(-\gamma-\text { PolyGamma }(0,2+n)) \rightarrow-\infty \\[10pt]
&\gamma=\text { Eulergamma }
\end{aligned}\tag{17}\]\\\\
\textbf{XIII. Generalisation of the sum}\\\\
We define the sum in relation to the $\zeta()$ function as follows\\
\[{ }_{ \quad m }^{-1,0} \sigma( s )=\sum_{n=1}^{\infty} \frac{1}{(-n)^{m \cdot s}}=\sum_{n=1}^{\infty} \frac{1}{(-2 \cdot n+1)^{m \cdot s}}+\sum_{n=0}^{\infty} \frac{1}{(-2 \cdot n-2)^{m \cdot s}}\tag{18}\]\\
If $\beta=1$ \& $m=1$ then if we replace relations (13\&18) we get that the sum\\
\[{ }_{ \quad 1 }^{-1,0} \sigma( s )=\sum_{n=1}^{\infty} \frac{1}{(-n)^{1 \cdot s}}=\sum_{n=1}^{\infty} \frac{1}{(-2 \cdot n+1)^{1 \cdot s}}+\sum_{n=0}^{\infty} \frac{1}{(-2 \cdot n-2)^{1 \cdot s}} \rightarrow-\infty+(-\infty) \rightarrow-\infty\tag{19}\]\\
Therefore from relation (17,19) the sum\\
\[{ }_{ \quad 1 }^{-1,0} \sigma( s )=\sum_{n=1}^{\infty} \frac{1}{(-n)^{1 \cdot s}} \rightarrow-\infty\]\\\\
From relations (8\&15) we conclude that for $m=1$\\
\[\begin{aligned}
{ }_{ \quad 1 }^{-1,0} \sigma( s )&=\sum_{n=1}^{\infty} \frac{1}{(-n)^{1 \cdot s}}=\sum_{n=1}^{\infty} \frac{1}{(-2 \cdot n+1)^{1 \cdot s}}+\sum_{n=0}^{\infty} \frac{1}{(-2 \cdot n-2)^{1 \cdot s}}= \\[10pt]
&=(-1)^{-s}(-2)^{-s} \cdot \zeta(s)+(-1)^{-s}(-2)^{-s} \cdot\left(-1+2^s\right) \cdot \zeta(s)=(-1)^{-2s} \cdot \zeta(s)
\end{aligned}\tag{20}\]\\\\
if $s=\sigma+\tau i$ and $\sigma=1 / 2$ and $\zeta( s )=0$, $m=1$ then the sum will be\\
\[{ }_{ \quad 1 }^{-1,0} \sigma( s )=\sum_{n=1}^{\infty} \frac{1}{(-n)^{s}}=(-1)^{-2s} \cdot \zeta(s)=-e^{2\cdot\pi \cdot \tau} \cdot \zeta(s)=-e^{2\cdot\pi \cdot \tau} \cdot 0=0\tag{21}\]\\\\
But from relation (20,21) if $\zeta(s)$ is zeroed then if we zero and the sum ${ }_{ \quad 1 }^{-1,0} \sigma( s )=\displaystyle \sum_{n=1}^{\infty} \frac{1}{(-n)^{s}}=0$ will have nontrivial zeros same to $\zeta( s )=0$. Moreover, if we compute this sum with trivial zeros of $\zeta( s )=0$ we will find for the first 10 consecutive zeros with conjugate roots the values according to the following table 27.\\\\\\
\textbf{XIV. Calculation with program.}\\\\
So we have the corresponding program that produces the ZetaZeros of $\zeta( s )=0$ and we find the corresponding values for Sums. The data are $a 1=-2$, $b 1=1$, $a 2=-2$, $b 2=-2$, $m =1$.\\\\
Off[FindRoot::1stol,N::meprec,Sum::div,General::munfl,NIntegrate::slwcon,NIntegrate::ncvb];\\\\\\
\textbf{'program data'}\\\\
a1=Input["Value of a1="];\\[8pt]
b1=Input["Value of b1="];\\[8pt]
a2 $=\operatorname{Input["Value~of~a2="];~}$\\[8pt]
b2=Input["Value of b2="];\\[8pt]
$m =$ Input["Multiplicity of exponent of $m =$ "];\\\\\\
\textbf{'program execution'}\\\\
$\begin{aligned}
& h 1\left[ x _{-}\right]:= N \left[\sum_{ k =1}^{\infty} \frac{1}{\left( a _2 \cdot k + b _2\right)^{\wedge}( m \cdot x )}, 10\right] ;(\text {Sum1}) \\[8pt]
& h 2\left[ x _{-}\right]:= N \left[\sum_{ k =0}^{\infty} \frac{1}{\left( a _2 \cdot k + b _2\right)^{\wedge}( m \cdot x )}, 10\right] ;(\text {Sum}2) \\[8pt]
&\text { For }[ i =1, i <=15, i ++\text {, } \\[8pt]
& xn = N [\text { ZetaZero[i],550]; } \\[8pt]
&\text { Print[i,",",N[xn,10],h1[xn],h2[xn]]]}
\end{aligned}$
\newpage\noindent
\begin{figure}[h!]
    \centering
    \includegraphics[width=14.5cm, height=7.3cm]{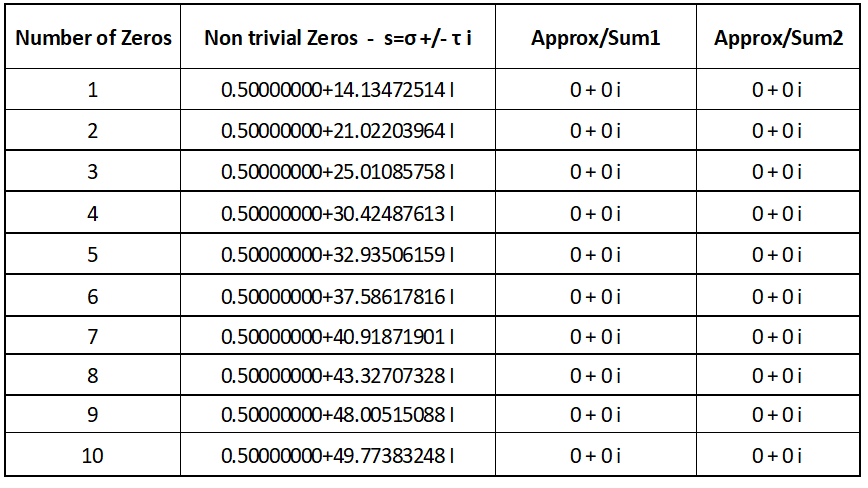}
\end{figure}\\
\begin{center}
  \textbf{Table 27. If $\boldsymbol{(\alpha 1=-2 \& \beta 1=1) \&(\alpha 2=-2 \& \beta 2=-2)}$ approx. for the Sum1 and Sum2}  
\end{center}
\vspace{\baselineskip}
As we see the results of approximate sums Sum1 $=\sum_{k=1}^{\infty} \frac{1}{(-2 \cdot k+1)^x}$ tends to complex zero and the Sum2 $=\sum_{k=0}^{\infty} \frac{1}{(-2 \cdot k-2)^x}$ it tends to complex zero. It is now obvious and computationally that the sum that is zeroed by the non-trivial zeros of $\zeta( s )=0$ is not the known sum $\sum_{k=1}^{\infty} \frac{1}{(k)^x}$ but the negative sum $\sum_{k=1}^{\infty} \frac{1}{(-k)^x}$ that we have shown ends in a common factor $\zeta( s )$.\\\\\\
\textbf{The Riemann hypothesis is therefore overturned on 2 points.}\\\\
\textbf{First.}Is the wrong choice of the Sum associated with the function $\zeta()$(Relation (p. 66-75)) and \\
\textbf{Second.}There are other generalized functions in combination with Hurwitz zeta function that have zeros with real part in the interval $(0,1)$ but also critical lines in the interval. (p. 60-65)
\newpage\noindent
\begin{center}
\textbf{Epilogue to Chapter 6}
\end{center}
\vspace*{\baselineskip}
\textbf{I. The zeros of $\boldsymbol{\zeta (z)=0}$}\\\\
The fact that the roots of $\zeta(\mathrm{z})=0$ all lie on the critical line, follows from 5 proofs, both basic and
elementary.\\\\
1.Using direct method of inverse of $\zeta ()$ on functional equations, chapter 6.1.1.\\\\
2.Solutions of the Equation $\zeta(\mathrm{s}) / \zeta(1-\mathrm{s})=1, \mathrm{~s}$ is Complex number, chapter $6.1 .3$.\\\\
3.With general inverse random function theory, chapter $6.2$.\\\\
4.Statistical evidence method chapter $6.3$.\\\\
5.With the define the function $\lambda(\mathrm{N})$, chapter $6.4 .1$.\\\\
Therefore, the proof that all zeros of $\zeta(\mathrm{z})=0$ belong on the critical line is complete, since the function $\zeta(\mathrm{z})=0$ as proved many times has zeros only on the critical line which means that  $\operatorname{Re}(\mathrm{z})=1 / 2$ and only this one.\\\\\\
Instances I-V of the characteristic generalized equation ${ }_{ \text{ }\text{ } \mathsmaller{G} } \zeta( s , q )=0, \alpha<0, \alpha, \beta \in R ( p .52-63)$ as well as the case $\zeta( q * s )=0$ with $q >1 / 2$ (p. 30-32) are subcases of the generalized mentioned above, and some cases of the periodic Dirichlet Davenport-Heilbronn series (page 51) rejectthe Riemann hypothesis. Certainly \textbf{the Paradox is the rejection by definition} (p.66-76) and lay the foundation for the definitive rejection of Riemann's Hypothesis. Therefore the Hypothesis is latent and the simple reason that the $\zeta()$ function does not correspond to the sum assumed has been constructed for to be $\zeta( s )=0$. But if it is formulated as follows: "The zeros of only $\zeta( s )=0$ and partial associated functions coincide in the critical line for $\operatorname{Re}( s )=1 / 2$ is correct and the Hypothesis holds", as has been shown. Obviously this happened because the equations of other related functions i.e. other related equations arising from related functions with the function $\zeta()$ were not fully or adequately solved or not solved at all and there was a large knowledge gap about where their zeros are located with respect mainly to the interval $(0,1)$. These equations are forms of transcendental equations quite difficult and their solution requires a powerful method (GRIM) and the corresponding generalized theorem for solving such equations.
\newpage
\section*{Part III. Distribution of primes numbers and Zeros}
\textbf{7. Introduction}\\\\
One way to get a feeling for where the primes are within the natural numbers is to count the number of primes $\pi (n)$ less than or equal to a given number $n$. There is no real hope for finding a simple rule for $\pi (n)$, but the Prime Number Theorem1 (PNT) says that $\pi (n)$ can be approximated by the logarithmic integral function. The \underline{distribution of prime numbers} is most simply expressed as the (discontinuous) step function $\pi (x)$, where $\pi (x)$ (Fig.9), is the number of primes less than or equal to $x$.
\phantom{ }\\
\flushleft
\justifying
\begin{figure}[h!]
\centering
\includegraphics[scale=0.55]{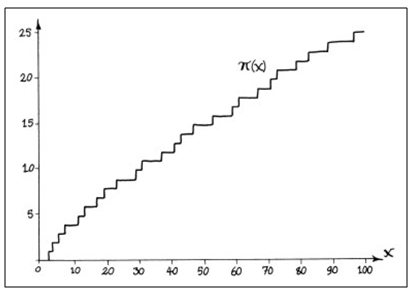}
\caption{The primes of function $\pi (x)$}
\end{figure}
The function $\pi (x)$ in relation to the random int. $x$. It turns out that $\pi (x)$ can be expressed exactly as the limit of a sequence of smooth functions $R_{n}(x)$. To define $R_{n}(x)$ we first introduce the \underline{logarithmic integral} function $Li(x)$,\\
\[Li(x)=\int_{2}^{x}\dfrac{du}{\log u}\]\\
which appears throughout the analytical theory of the prime distribution:\\\\       
This is a smooth function which simply gives the area under the curve of the function $1/\log u$ in the interval $[2,x]$. Don Zagier explains the reasoning behind the function $Li$ in his excellent introductory article ``The first 50 million prime numbers" [14], 
based on his inaugural lecture held at Bonn University, May 5, 1975: ``A good approximation to $\pi (x)$", which was first given by Gauss is obtained by taking as starting point the empirical fact that the frequency of prime numbers near a very large number $x$ is almost exactly $1/\log x$. From this, the number of prime numbers up to x should be approximately given by the \underline{logarithmic sum} $Li(x)=1/\log 2 + 1/\log 3 + 1/\log x$ or, what is essentially the same, by the logarithmic integral.                             
\[Li(x)=\int_{2}^{x}\dfrac{du}{\log u}\]\\
Using $Li(x)$ we then define another smooth function, $R(x)$, first introduced by Riemann in his original eight\\\-page paper, and given by\\
\[R(x)=\sum_{n=1}^{\infty}\dfrac{\mu (n)}{n}Li(x^{1/n})\]\\
Riemann's research on prime numbers suggests that the probability for a large number $x$ to be prime should be even closer to $1/\log x$ if one counted not only the prime numbers but also the powers of primes, counting the square of a prime as half a prime, the cube of a prime as a third, etc. This leads to the approximation:\\
\[\pi(x)+\dfrac{1}{2}\pi(x^{1/2})+\dfrac{1}{3}\pi(x^{1/3})+...\approx Li(x)\]\\
or, equivalently [\textit{by means of the Möbius inversion formula}]\\\\
\[\boxed{\pi(x)\approx Li(x)-\dfrac{1}{2}Li(x^{1/2})-\dfrac{1}{3}Li(x^{1/3})-...}\]\\\\
The function on the right side of this formula is denoted by $R(x)$, in honour of Riemann. It represents an amazingly good approximation to $\pi (x)$. For those in the audience who know a little function theory, perhaps I might add that $R(x)$ is an entire function of $\log x$, given by the rapidly converging power series:\\
\[R(x)=1+\sum_{k=1}^{\infty}\dfrac{(\ln x)^{k}}{kk!\zeta (k+1)}\]\\
where $\zeta (k+1)$ is the Riemann zeta. 
\phantom{ }\\
\flushleft
\justifying
\begin{figure}[h!]
\centering
\includegraphics[scale=0.55]{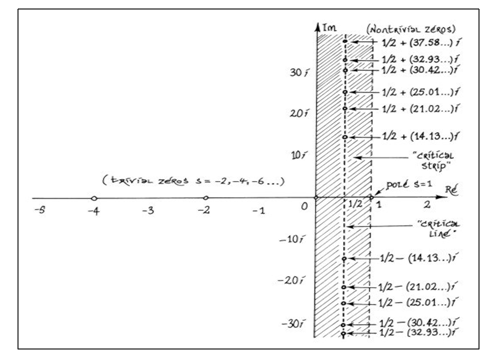}
\caption{The non-trivial zeros in the Critical Line}
\end{figure}
\justifying{
Here we see the zeros of the Riemann zeta function in the complex plane. These fall into two categories, \textit{trivial} and \textit{nontrivialzeros} (Fig.10). Here are some tables on nontrivial zeros compliled by Andrew Odlyzko[4]. The trivial zeros are simply the negative even integers. The nontrivial zeros are known to all \textbf{lie in the critical strip} that is $0 < Re[s] < 1$, and always come in complex conjugate pairs. All known nontrivial zeros lie on the \textbf{\underline{critical line}} $Re[s] = 1/2$. The \underline{Riemann Hypothesis} states that they \textit{all} lie on this line. The difference between the prime counting function and its ``amazingly good approximation" $R(x)$, i.e. the \textit{fluctuations} in the distribution of primes, can be expressed in terms of the entire set of zeros of zeta, which we shall represent by $\rho$, via the function $R$ itself:  }           
\[R(x)-\pi (x)=\sum_{\rho}R\left(x^{\rho}\right)\]    
Obviously some of the $x^{\rho}$ are complex values, so here R is the analytic continuation of the real-valued function R defined previously. This was mentioned above by Zagier, and is known as the \underline{Gram Series} expansion:
\[\sum_{k=1}^{\infty}\left[ R\left(x^{\rho_{k}}\right)+R\left(x^{\rho_{-k}}\right)\right]\]      
The contributions from the complex-conjugate pairs $\rho_{k}$ and $\rho_{-k}=\overline{\rho_{k}}$ cancel each others' imaginary parts, so
\[\boxed{\pi(x)=R(x)-\sum_{m=1}^{\infty}R\left(x^{-2m}\right)+\sum_{k=1}^{\infty}T_{k
}(x)}\]
Where $T_{k}(x)=-R\left(x^{\rho_{k}}\right)-R\left(x^{\rho_{-k}}\right)$ are real-valued. We can now define the sequence of functions $R_{n}(x)$ which approximate $\pi(x)$ in limit:
\[R_{n}(x)=R(x)-\sum_{m=1}^{\infty}R\left(x^{-2m}\right)+\sum_{k=1}^{n}T_{k}(x)\tag{1}\]
\[f(x)=-2Re\left(\sum_{k=1}^{\infty}Ei\left(\rho_{k}\log (x)\right)\right)+\int_{x}^{\infty}\dfrac{1}{(t^{3}-1)\log (t)}dt-\log (2)\tag{2}\]
where $\rho_{\kappa}$ is the $\kappa^{\text{th}}$ complex zero of the zeta function. In this formula, $Ei(z)$ (\textit{Mathematica's} built-in function ExpIntegral $Ei(z)$ is the generalization of the logarithmic integral to complex numbers. These equations come from references [13], [14], and [15]. First, let $M$ be the smallest integer such that $x^{1 / M}<2$. We need to add only the first M-1 terms (that is, $n=1,2,..,M$) in the sum in equation (1). For each of these values of $N$, we use equation (2) to compute the value of $f\left(x^{1 / n}\right)$. However, we will add only the first $N$ terms (that is, $\kappa=1,2,..,N$) in the sum in equation (2). Because the purpose of this Demonstration is to show how the jumps in the step function $\pi(x)$ can be closely approximated by adding to $R(x)$ a correction term that involves zeta zeros, we ignore the integral and
the $\log 2$ in second equation; this speeds up the computation and will not noticeably affect the graphs, especially for $x$ more than about 5. The more zeros we use, the closer we can approximate $\pi(x)$. For larger $x$, the correction term must include more zeros in order to accurately approximate $\pi(x)$.\\\\
\textbf{8. The GUE hypothesis.}\\\\
While many attempts to prove the $\mathrm{RH}$ had been made, a
few amount of work has been devoted to the study of the
distribution of zeros of the Zeta function. A major step
has been done toward a detailed study of the distribution
of zeros of the Zeta function by Hugh Montgomery [18],
with the Montgomery pair correlation conjecture.
Expressed in terms of the normalized spacing. 
$$\boxed{\delta_{n}=\left(\gamma_{n+1}-\gamma_{n}\right) \dfrac{\log \left(\gamma_{n} /(2 \pi)\right)}{2 \pi}}$$
this conjecture is that, for $M\rightarrow \infty$
$$
\begin{aligned}
\dfrac{1}{M} \#\left\{n: N+1 \leq n \leq N+M, \delta_{n} \in[\alpha, \beta]\right\} 
\sim \int_{a}^{b} p(0, u) d u \sim \int_{a}^{b} 1-\left(\dfrac{\sin \pi u}{\pi u}\right)^{2} d u
\end{aligned}
$$
In other words, the density of normalized spacing between non-necessarily consecutive zeros is $1-\left(\dfrac{\sin \pi u}{\pi u}\right)^{2}$. It was first noted by the Freeman Dyson, a quantum physicist, during a now-legendary short teatime exchange with Hugh Montgomery[18], that this is precisely the pair correlation function of eigenvalues of random hermitian matrices with independent normal distribution in (Fig.11, 12) of its coefficients. Such random hermitian matrices are called the Gauss unitary ensemble (GUE). As referred by Odlyzko in [16] for example, this motivates the GUE hypothesis which is the conjecture that the distribution of the normalized spacing between zeros of the Zeta function is asymptotically equal to the distribution of the GUE eigen values. Where
$p(0,u)$ is a certain probability density function, quite complicated to obtain for an expression of it). As reported by Odlyzko in [16], we have the Taylor expansion around zero. $p(0, u)=\dfrac{\pi^{2}}{3} u^{2}-2 \dfrac{\pi^{4}}{45} u^{4}+\ldots$
which under the GUE hypothesis entails that the proportion of $\delta$ n less than a given small value $\delta$ is asymptotic to $\dfrac{\pi^{2}}{9} \delta^{3}+O\left(\delta^{5} \right)$. Thus very close pair of zeros are rare.\\\\
\phantom{ }\\
\begin{figure}[h!]
\centering
\includegraphics[scale=0.5]{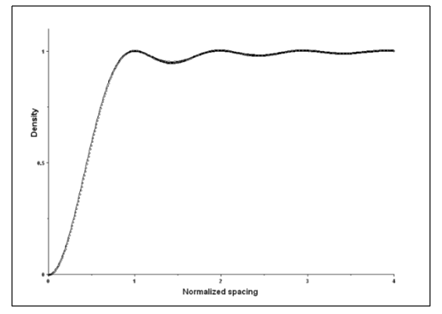}
\caption{Probability density of the normalized spacing }
\end{figure}
\flushleft
\justifying
between non-necessarily consecutive zeros and the GUE prediction[20]
\phantom{ }\\
\begin{figure}[h!]
\centering
\includegraphics[scale=0.5]{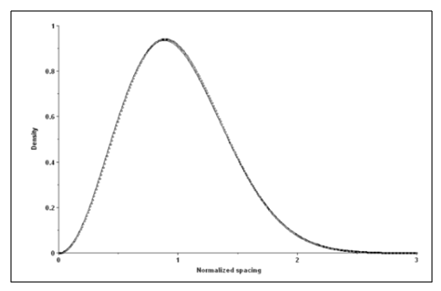}
\caption{Probability density of the normalized spacing $\delta n$ and the GUE prediction[20].}
\end{figure}\\\\
\textbf{9. Gaps between zeros}\\
The table below lists the minimum and maximal values
of normalized spacing between zeros $\delta n$ and of
$\delta_{n}+\delta_{n+1}$, and compares this with what is expected
under the GUE hypothesis. It can be proved that $p(0,t)$ have the following Taylor expansion around 0, $p(0, u)=\dfrac{\pi^{2}}{3} u^{2}-2 \dfrac{\pi^{4}}{45} u^{4}+\cdots$
so in particular, for small delta
$$
\boxed{Prob\left(\delta_{n}<\delta\right)=\int_{0}^{\delta} p(0, u) d u \sim \dfrac{\pi^{2}}{9} \delta^{3}}
$$
so that the probability that the smallest $\delta n$ are less than $\delta$ for $M$ consecutive values of $\delta n$ is about.
$$
1-\left(1-\dfrac{\pi^{2}}{9} \delta^{3}\right)^{M} \simeq 1-\exp \left(-\dfrac{\pi^{2}}{9} \delta^{3} M\right)
$$
This was the value used in the sixth column of the table. The result can be also obtained for the $\delta_{n}+\delta_{n+1} \ldots$
$$
\boxed{Prob\left(\delta_{n}+\delta_{n+1}<\delta\right) \sim \dfrac{\pi^{6}}{32400} \delta^{8}}
$$
from which we deduce the value of the last column (Table 28).
\setcounter{table}{27}
\phantom{ }\\
\begin{table}[h!]
\centering
\includegraphics[width =11cm, height=5cm]{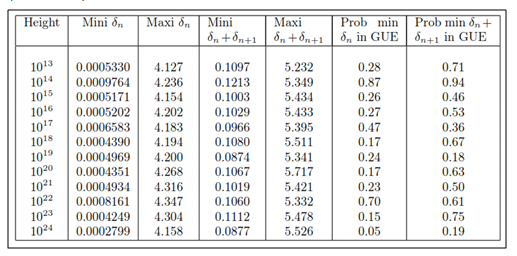}
\caption{The values of $\delta_{n}+\delta_{n+1}$} 
\end{table}\\
\noindent For very large spacing in the GUE, as reported by
Odlyzko in $[16]$, des Cloizeaux and Mehta $[17]$ have
proved that
$$
\log p(0, t) \sim-\pi^{2} t^{2} / 8 \quad(t \rightarrow \infty)
$$
which suggests that
$$
\max _{N+1 \leq n \leq N+M} \delta_{n} \sim \dfrac{(8 \log M)^{1 / 2}}{\pi}
$$
\phantom{ }\\
\begin{center}
\textbf{Statistics of False zetazeros, $\boldsymbol{\delta}$-intervals and count of primes.}
\end{center}
\textbf{10. General equations.}\\\\
With this statistic we find the crowd of individual parts intervals defined by consecutive Zetazeros in a fixed integer interval. Here we use $\delta=1000$. We have the general Equation $\boldsymbol{\delta+}$ \textbf{Zetazero (kin) - Zetazero (kf)} $\boldsymbol{=0}$ and after given initial value in kin, i calculate the \textbf{kf usually by the Newton method}. From the initial and final value of \textbf{kin, kf} and by performing the process of successive intervals, i calculate the number(count) of the primes ones that are within the intervals. In this way will have consecutive intervals.
$S_{k}=\{$ Zetazero $(k+1)<p<$ Zetazero $(k), k \in(k i n, k f)\}$
from where resulting the number of the primes in the given consecutive interval. The sum of the primes in the given interval $\delta=1000$ will be obvious.
$$S_{\delta}=\bigcup_{i=k_{i n}}^{k_{f}} S_{i}$$
The False intervals will be in a normalized form $\mathrm{F}_{\delta}=1000-\mathrm{S}_{\delta}$ over $N>=4000$.\\\\
\textbf{11. The statistics}\\\\
The First Statistic looking for ways to show us which function is the most ideal to get closer the points of interest, uses the NonlinearModelFit method of the
function $y=(a+b \cdot x) /(\log (d+c \cdot x)$ and therefore after determining the variables $\{a, b, c, d\}$ we are able to make statistical and probable prediction at higher levels of numbers . This function is directly related to the function $\pi(x)=x / \log x$ which was reported in the introduction that is, defining the number of primes
numbers relative to $x$. We will do a double statistic of the intervals $\boldsymbol{\lbrace\delta }$, \textbf{False intervals} $\boldsymbol{\rbrace }$ and $\boldsymbol{\lbrace \delta}$, \textbf{number of primes} $\boldsymbol{\rbrace}$ that we are ultimately interested in the statistics mainly the count of the primes inside at $\delta$ -intervals.\\\\
\textbf{11.1}\textbf{ The first statistic} is about count of \textbf{False intervals}
and the \textbf{count of the primes} that corresponding in them. The table below gives it in aggregate until the count of primes equal $p_{\delta}=29$ corresponding to the count of False $F_{\delta}=971$ and the Total range Integers for a interval $D=1000$ to $2 \cdot 10^{14}+1000$, (Table 29).
\setcounter{table}{28}
\flushleft
\justifying
\begin{table}[htb]
\centering
\includegraphics[width  = 7cm, height =5cm]{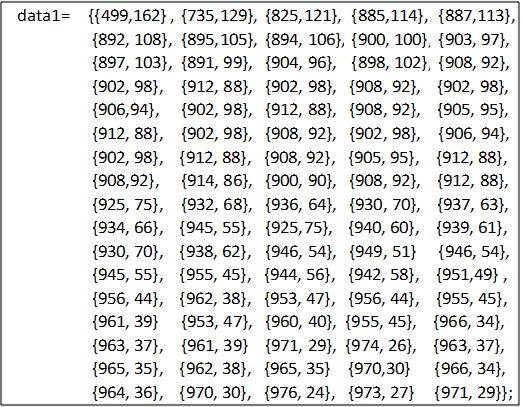}
\caption{The data.1} 
\end{table}
\phantom{ }\\
The diagram given below (Fig.13) shows the arrangement of the points at the level $(x, y)$ according to the data1[19].
\begin{figure}[htb]
\centering
\includegraphics[scale=0.6]{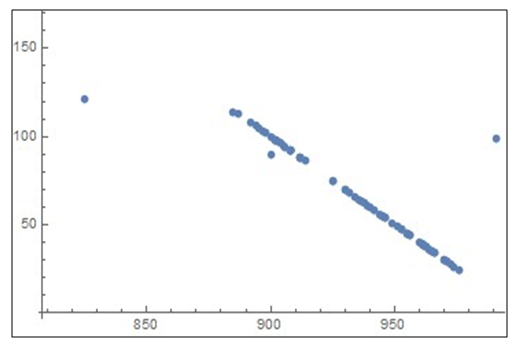}
\caption{The depiction of the intervals False and the primes ones (file:data1) contained in the intervals  $\boldsymbol{N+\delta-N = 1000}$.}
\end{figure}\\
Even more macroscopically (Fig.14), the points are shown by the line $\boldsymbol{y = 988.709 - 0.988372*x}$, which, as they appear, are stacked on its lower right,
\begin{figure}[h!]
\centering
\includegraphics[scale=0.64]{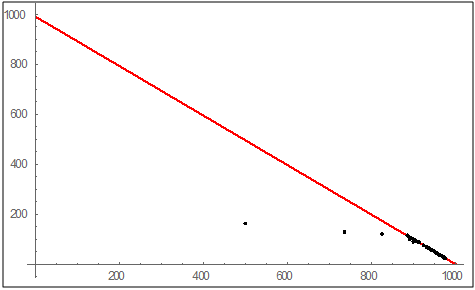}
\caption{Graphic depiction of the line $y = 988.709 - 0.988372*x$  with the archive(data1) points.}
\end{figure}\\
It is obvious that for $\{x = 999, y = 1.32\}$ and for $\{y = 1000, x = 0.33\}$, a value located below the unit and means that it can in such a interval, and we are talking about this is for high order integers intervals so there is not one prime within the interval $\delta = 1000$, chosen at random.\\\\
Basic thinking usually in statistics is to calculate at constant shows the number of data we need respectively, finding the frequencies that appear. This method is a real percentage of realistic Espigue isolates the properties showing a group of numbers.\\\\
\textbf{11.2 The second statistic} refers to the count of the primes ones that located in the intervals $[N, N+\delta]$ with $N=1000$ up to $2 \cdot 10^{14}$, range $\delta=1000$ and the \textbf{count of the primes within successive ZetaZeros}. After we found the roots kf of general Equation $\boldsymbol{\delta+}$ \textbf{Zetazero (kin) -Zetazero (kf)} $\boldsymbol{=0}$ of given initial value kin, usually by the Newton method. The data 19 (Table 30), that we have met with the above method are.
\phantom{ }\\
\begin{table}[htb]
\centering
\includegraphics[scale=0.66]{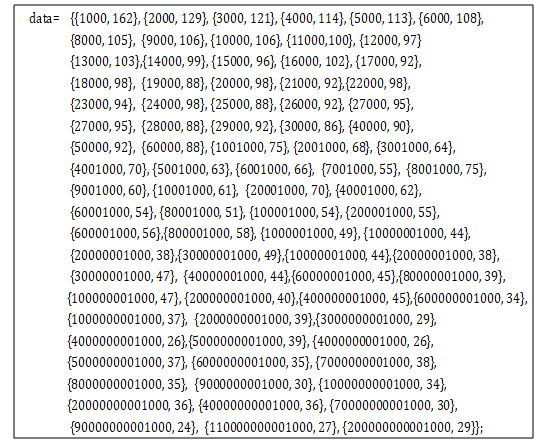}
\caption{The data} 
\end{table}\\
\noindent Using the NonlinearModelFit[19] method of the function $y=(a+c*x)/(\log[d+b*x])$ and after specifying the variables $\{a, b, c, d\}$ the following result the function will be:
$$y=\dfrac{1108.254246288494-1.116325646187134 \times 10^{-12} x}{\log [-4503.90336177023+5.425635171403081 x]}$$
with very good approach and value performance in each pairing.The diagram (Fig.15) of the above equation is shown in more detail below with gravity in the latest data.
\begin{figure}[htb]
\centering
\includegraphics[scale=0.6]{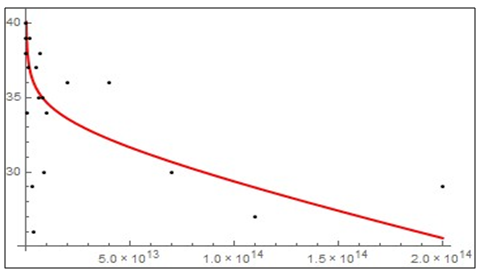}
\caption{The graphic depiction of the function}
\end{figure}\\
$y=(a+c * x) /(\log [d+b * x])$. Thus a test value for $x= 2 \cdot 10^{14}$ gives us $y=25.6$, close to 29 we took with the analysis. They follow \textbf{statistics ANOVA} and t-Statistic $[19,20]$ and we get the (Tables 31, 32).\\\\
Analysis of variance (ANOVA) is an analysis tool used in statistics that splits an observed aggregate variability found inside a data set into two parts: systematic factors and random factors. The systematic factors have a statistical influence on the given data set, while the random factors do not. An ANOVA test is a way to find out if survey or experiment results are significant. In other words, they help you to figure out if you need to reject the null hypothesis or accept the alternate hypothesis.
\phantom{ }\\
\flushleft
\justifying
\begin{table}[htb]
\centering
\textbf{Analysis-ANOVA}\\
\includegraphics[scale=0.65]{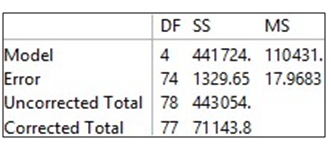}
\caption{} 
\end{table}
\begin{table}[htb]
\centering
\textbf{Analysis-T-statisrtic}\\
\includegraphics[scale=0.65]{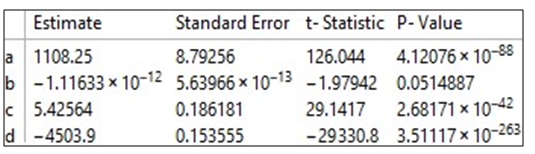}
\caption{} 
\end{table}
As we can see from the results, we have an important Standard Error only for \textbf{the variable a}. The other variables are observed to have a low statistical error and have good compatibility. By adapting the method as we see, we associate two lists of results of the number of primes and long intervals, which although disproportionately, work together with impeccable and good contact.\\\\
\underline{Ronald Fisher} introduced the term variance and proposed its formal analysis in a 1918 article The Correlation\\
Between Relatives \underline{on the Supposition of Mendelian Inheritance}. His first application of the analysis of variance was published in 1921. Analysis of variance became widely known after being included in Fisher's 1925 book \underline{\textit{Statistical Methods for Research Workers}}.\\\\
One of the attributes of ANOVA that ensured its early popularity was computational elegance. The structure of the additive model allows solution for the additive coefficients by simple algebra rather than by matrix calculations. In the era of mechanical calculators this simplicity was critical. The determination of statistical significance also required access to tables of the F function which were supplied by early statistics texts.\\\\
\textbf{12. Standard–Error and Confidence–Interval}\\\\
The \textbf{standard error (SE)} of a statistic (usually an estimate of a parameter) is the standard deviation of its sampling distribution or an estimate of that standard deviation. If the statistic is the sample mean, it is called the \textbf{standard error of the mean (SEM)}. The sampling distribution of a population mean is generated by repeated sampling and recording of the means obtained. This forms a distribution of different means, and this distribution has its own mean and variance. Mathematically, the variance of the sampling distribution obtained is equal to the variance of the population divided by the sample size. This is because as the sample size increases, sample means cluster more closely around the population mean. Therefore, the relationship between the standard error of the mean and the standard deviation is such that, for a given sample size, the standard error of the mean equals the standard deviation divided by the square root of the sample size. In other words, the standard error of the mean is a measure of the dispersion
of sample means around the population mean. Where we observed, in addition to the initial measurements \#11-\#12, the Standard - Error also the confidence interval stabilized at good and acceptable values with central average values around the lateral expected values. The values are given in two consecutive tables that we add below and refer to the whole values of range (Table 33).
\flushleft
\justifying
\begin{table}[htb]
\centering
\begin{minipage}{0.5\textwidth}
\includegraphics[width=8cm, height=13cm]{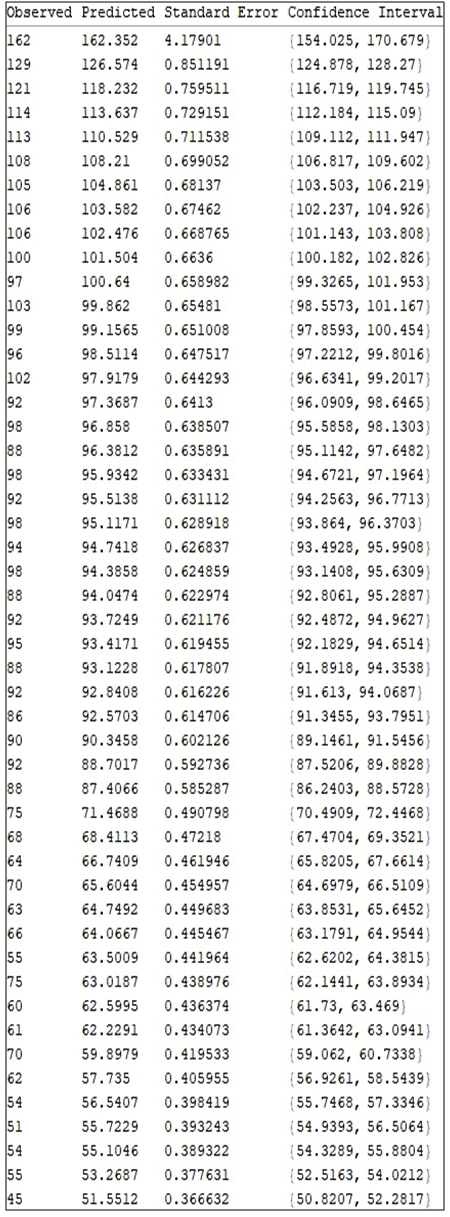}
\end{minipage}\hfill
\begin{minipage}{0.5\textwidth}
\vspace*{-13.49\baselineskip}
\includegraphics[width=8cm, height=7.28cm]{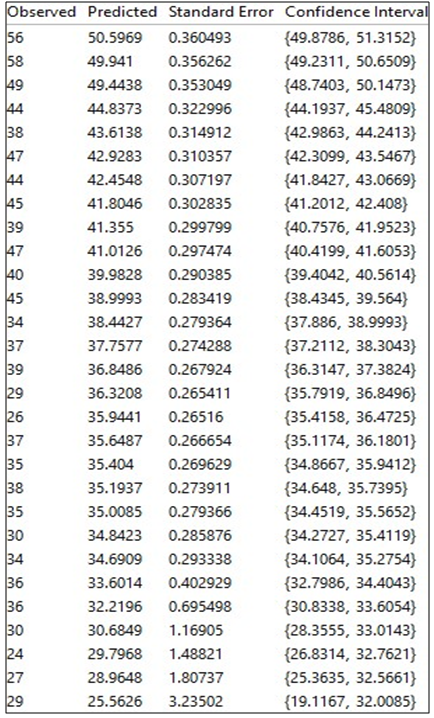}
\end{minipage}
\caption{} 
\end{table}
\textbf{13. Statistical comparison}\\\\
With statistical comparison it appears that as long as we moving to a higher order of integer-size that are within a given interval, we have chosen $\delta=1000$, the count of the primes diminishes to disappear or to there are $1-2$ primes at high order intervals of more than $10^{20}$. This also agrees with of \textbf{Gram's law}. In particular, the problem of distribution of the differences $t_{n+1}-t_{n}$ \textbf{(that is of difference of ZetaZeros)} is considered [21]. If we accept Gram's law then the order of this difference does not exceed the quantity  $t_{n+1}-t_{n} \approx \frac{2 \pi}{\ln (n)} \rightarrow 0, n \rightarrow \infty$ for much larger integers $[20,21]$ then mean and their mean value is close to,
$$t_{n+1}-t_{n} \approx \frac{\ln \ln (n)}{\ln (n)} \rightarrow 0, n \rightarrow \infty$$
And as the above analysis of-\textbf{False intervals}- as we have shown, \textbf{it is compatible with this result of law Gram's}.\\\\
\flushleft
\justifying
\textbf{14. From above analysis they arise 3 big conclusions:}\\\\
$\boldsymbol{1^{\text{st}}}$. The number(Count) of primes located within $\delta$ intervals is gradually \underline{decreasing} with a higher order size $\mathrm{n}$ of $10^{n}$.\\
$\boldsymbol{2^{\text{nd}}}$. The distribution of the number(Count) of the primes, within interval $\delta$ follows a \underline{Nonlinear correlation} by the function $y=(a+b \cdot x) /(\log (d+c \cdot x)$ similar to the number (Count) of the primes $\pi(x) \approx x / \log (x)$, which is apply approximating for large numbers.\\
$\boldsymbol{3^{\text{rd}}}$. The measurement(Count) of false intervals and the count of the primes follows a \underline{linear correlation} and it increases if the order of magnitude size of the integer increases.\\\\
\begin{center}
\textbf{\underline{Epilogue}}
\end{center}
If we would like to consider Riemann's Hypothesis without dealing with the sum corresponding to the function $\zeta()$ the question is not simply whether Riemann's hypothesis is valid or not, but which functions $\zeta()$ make it admissible, so that as defined by Riemann for the function $\zeta()$ alone, it can be true. Thus, if we consider as a general case of the hypothesis all parallel functions $\zeta()$, then it is not $100 \%$ true. Of course there is the Paradox in the Riemann Hypothesis for which we have \textbf{non-trivial zeros that zero out a Sum} that is different from what the function $\zeta( s )$ represents and not what it should represent. This is a strong pre-hypothesis that raises the very serious matter that the Riemann's Hypothesis \textbf{should be rejected by definition} or its definition should be revisited on another basis in a more general form that encompasses all cases. Of course, it is also true that the Riemann hypothesis has multiple connections with other sciences and gives directions appropriate and grouped according to the \textbf{function} we accept to apply to any particular application case. If we accept generalized functions $\zeta()$, this error with respect to only real answer to the conjecture of the Riemann hypothesis generalizes and corrects the hypothesis. A hypothesis that has occupied the mathematical community for more than 160 years I believe that with the existing explanations it is clarified and takes the correct form it he should had from the beginning. Certainly the function $\zeta()$ is very interesting and helps in the right direction for solving this complicated conjecture according to the structure and analysis of the correct sum, as proven in this paper.
\newpage\noindent
\phantomsection
\addcontentsline{toc}{section}{References}
\thispagestyle{plain}
\renewcommand\refname{References}


\begin{thebibliography}{26}
\bibitem{b19} The “data1 \& data \& NonlinearModelFit Statistic”, have been done with the program mathematica.
\end{thebibliography}
\end{document}